\documentclass[12pt]{amsart}
\usepackage{graphicx}

\usepackage{amsmath}
\usepackage{amsthm}
\usepackage{amssymb}
\usepackage[mathscr]{eucal}
\begin{document}

\title{Elliptic minuscule pairs and splitting abelian varieties}%
\author{V.~Kumar Murty and Ying Zong}
\address{Department of Mathematics\\University of
Toronto} \email{murty@math.toronto.edu, zongying@math.toronto.edu}


\maketitle
\section{Introduction : why elliptic minuscule pairs}
\smallskip

The motivation of this article was to study the following question (cf.
\cite{murty}, 1.1) through a monodromy
approach.

\smallskip

{\bf Question 1.1.} \emph{Let $A_K$ be an absolutely simple abelian variety over a
number field $K$. Does there exist a finite extension $L$ of $K$
such that the base change of $A_K$ to each finite extension of $L$
has simple specializations at a set of places of
positive density}?
\smallskip

Let us recall some notions before we formulate this question in more
precise terms and impose a natural hypothesis on $A_K$.
\smallskip

Let $t=\mathrm{Spec}(K)$, $\overline{t}$ a geometric point of $t$ and
$S$ a dense open sub-scheme of the normalization of
$\mathrm{Spec}(\mathbf{Z})$ in $t$ such that $A_t=A_K$ extends to an
abelian scheme $A$ over $S$.
\smallskip

We call an arbitrary $S$-fiber of $A$ a specialization of $A_t$. A
specialization $A_s=A\times_Ss$ at a point $s$ of $S$ is said to be simple if it is a
simple object in the category of $s$-abelian varieties up to
isogenies, that is, if
$\mathrm{End}_s(A_s)\otimes_{\mathbf{Z}}\mathbf{Q}$ is a
$\mathbf{Q}$-division algebra. And, a specialization $A_s$ is absolutely
simple if $A_s\times_s\overline{s}$ is simple for some geometric point $\overline{s}$ of $s$.
\smallskip

Recall that a subset $\Xi$ of $S\backslash\{t\}$ has (natural) density $d$ (\cite{serre_adic}, I--7), $0\leq d\leq 1$, if asymptotically in
$N\in\mathbf{R}$,
\[\mathrm{Card}(\{s\in\Xi, \mathrm{Card}(k(s))\leq N\})=d\
\frac{N}{\mathrm{log}\ N}+o(\frac{N}{\mathrm{log}\ N}).\]

As a fundamental example, the set \[\{s\in S\backslash\{t\}, k(s)\
\mathrm{is\ a\ prime\ field}\}\] has density $1$.
\smallskip

In general, the density of $\Xi$ is taken here to be the supremum of the densities of its ``measurable'' subsets.
\smallskip

What we asked above is whether there exists some finite extension
$L$ of $K$ such that for each finite extension $K'$ of $L$, if $S'$
denotes the normalization of $S$ in $t'=\mathrm{Spec}(K')$, the set
\[\{s'\in S'\backslash\{t'\}, A\times_Ss'\ \mathrm{is\ simple}\},\]
or what amounts to the same, the subset
\[\{s'\in S'\backslash\{t'\}, k(s')\
\mathrm{is\ a\ prime\ field}, A\times_Ss'\ \mathrm{is\ simple}\}\]
has positive density.
\smallskip

Note that
\smallskip

\smallskip

{\bf Lemma 1.2.} \emph{The algebra $\mathrm{End}_{t'}(A_{t'})\otimes_{\mathbf{Z}}\mathbf{Q}$ must be a field
if $A_{t'}$ has at least one simple specialization
$A_{s'}$ at a point $s'$ with
finite prime residue field.}

\begin{proof} For, $\mathrm{End}_{s'}(A_{s'})\otimes_{\mathbf{Z}}\mathbf{Q}$ is a field at such a point $s'$
(\cite{tate}, p.~98, line 1) and the specialization homomorphism
\[ sp: \mathrm{End}_{t'}(A_{t'})=\mathrm{End}_{S'}(A_{S'})\hookrightarrow
\mathrm{End}_{s'}(A_{s'})\] is injective.
\end{proof}

In particular, our question has a negative answer unless
\[E:=\mathrm{End}_{\overline{t}}(A_{\overline{t}})\otimes_{\mathbf{Z}}\mathbf{Q}\]
is a field, as was predicted in \cite{murty} and known to J. Achter in a less precise way (\cite{achter}, Theorem B).
\smallskip

One can ask if the hypothesis that
$E=\mathrm{End}_{\overline{t}}(A_{\overline{t}})\otimes_{\mathbf{Z}}\mathbf{Q}$ be a field is sufficient for the question to have a positive answer (cf. 1.5, 1.7).
\smallskip

For this, enlarge if necessary $K$ to a finite extension so that
\[\mathrm{End}_t(A_t)=\mathrm{End}_{\overline{t}}(A_{\overline{t}}).\]

Let $\ell$ be a prime number and let $\mathfrak{l}$ be a place of $E$ above $\ell$. Replacing $S$ by its open sub-scheme $S[1/\ell]$ if necessary, we assume that $\ell$ is prime to the residue characteristics of $S$.
Choose for each closed point $s\in S$ a geometric point
$\overline{s}$ located at $s$ and a ``chemin'' $ch_s$ connecting
$\overline{s}$ to $\overline{t}$ (SGA 1, \'{E}xpos\'{e} V, 7). Let
$F_s\in\pi_1(s,\overline{s})$ be the geometric Frobenius and
$F_s^*$ the image of $F_s$ under the composition
\[ \pi_1(s,\overline{s})\to
\pi_1(S,\overline{s})\stackrel{ch_s}{\longrightarrow}\pi_1(S,\overline{t})\stackrel{\rho_{\ell,\overline{t}}}{\longrightarrow}\mathrm{GL}_E(H^1(A_{\overline{t}},\mathbf{Q}_{\ell})),\]
where $\rho_{\ell,\overline{t}}$ is the $\ell$-adic monodromy
representation associated with the abelian scheme $A$. Let
$M_{\ell}=\mathrm{Im}(\rho_{\ell,\overline{t}})$ be the monodromy and
$M_{\ell}^{\mathrm{Zar}}$ its Zariski closure in
$\mathrm{GL}_E(H^1(A_{\overline{t}},\mathbf{Q}_{\ell}))$. Further enlarging
$K$ to a finite extension if necessary, one may assume that $M_{\ell}^{\mathrm{Zar}}$ is connected.
\smallskip

The group $M_{\ell}^{\mathrm{Zar}}$ is then by Faltings (\cite{faltings}, satz 3) reductive and (\emph{loc.cit.}, satz 4)
\[\mathrm{End}_t(A_t)\otimes_{\mathbf{Z}}\mathbf{Q}_{\ell}\ \widetilde{\to}\ \mathrm{End}_{M_{\ell}^{\mathrm{Zar}}}(H^1(A_{\overline{t}}, \mathbf{Q}_{\ell}))^{\mathrm{opposite}}.\]

If $V_{\mathfrak{l}}$ denotes the $E_{\mathfrak{l}}$-component of $H^1(A_{\overline{t}}, \mathbf{Q}_{\ell})$ and if $M_{\mathfrak{l}}^{\mathrm{Zar}}$ is the image of $M_{\ell}^{\mathrm{Zar}}$ in $\mathrm{GL}_{E_{\mathfrak{l}}}(V_{\mathfrak{l}})$, one has
\[E_{\mathfrak{l}}=\mathrm{End}_{M_{\mathfrak{l}}^{\mathrm{Zar}}}(V_{\mathfrak{l}}),\] which amounts to the absolute irreducibility of $V_{\mathfrak{l}}$ as a $E_{\mathfrak{l}}$-linear representation of $M_{\mathfrak{l}}^{\mathrm{Zar}}$. The image $M_{\mathfrak{l}}$ of $M_{\ell}$ in $M_{\mathfrak{l}}^{\mathrm{Zar}}(E_{\mathfrak{l}})$ is an open analytic subgroup by Bogomolov \cite{bogomolov}.
\smallskip

At each closed point $s$ of $S$, the commutant of $F_s^*$ on $V_{\mathfrak{l}}$ is
\[(\mathrm{End}_s(A_s)\otimes_{\mathbf{Z}}\mathbf{Q})^{\mathrm{opposite}}\otimes_EE_{\mathfrak{l}},\] as by Tate \cite{tate_invent}
\[\mathrm{End}_s(A_s)\otimes_{\mathbf{Z}}\mathbf{Q}_{\ell}\ \widetilde{\to}\ \mathrm{End}_{F_s^*}(H^1(A_{\overline{t}}, \mathbf{Q}_{\ell}))^{\mathrm{opposite}}.\]

Recall that if at one point $s\in S\backslash\{t\}$ with prime residue field $A_s$ is simple, then
\[\mathrm{End}_s(A_s)\otimes_{\mathbf{Z}}\mathbf{Q}=\mathbf{Q}(F_s^*)\] is a field of degree $2g$ over $\mathbf{Q}$, where $g=\mathrm{dim}(A_t)$. This Frobenius $F_s^*$ has all distinct eigenvalues on $H^1(A_{\overline{t}}, \mathbf{Q}_{\ell})$ and $F_s^*$ lies in a unique maximal torus of $M_{\ell}^{\mathrm{Zar}}$. So
\smallskip

\smallskip

{\bf Lemma 1.3.} \emph{If $A_t$ has at least one simple specialization at a point $s$ with finite prime residue field, then some, hence every, maximal torus of $M_{\ell}^{\mathrm{Zar}}$ acts on $H^1(A_{\overline{t}}, \mathbf{Q}_{\ell})$ without multiple weights.}
\smallskip

\smallskip

This weight multiplicity free condition has the following immediate implication :
\smallskip

\smallskip

{\bf Proposition 1.4.} \emph{Suppose that $M_{\ell}^{\mathrm{Zar}}$ is connected and that the monodromy representation $H^1(A_{\overline{t}}, \mathbf{Q}_{\ell})$ has no multiple weights. Then in a density $1$ set $\Sigma$ of points $s\in S\backslash\{t\}$ every positive power of $F_s^*$ has all distinct eigenvalues on $H^1(A_{\overline{t}}, \mathbf{Q}_{\ell})$. In particular, the simple factors of each $A_s$, $s\in\Sigma$, are absolutely simple mutually non-isogenous over $\overline{s}$ and $\mathrm{End}_{\overline{s}}(A_{\overline{s}})$ is commutative. A specialization at a point $s\in\Sigma$ is thus absolutely simple if it is simple.}

\begin{proof} Let $s$ be a point of $S\backslash\{t\}$. The Frobenius $F_s^{*}$ being semi-simple on
$H^1(A_{\overline{t}},\mathbf{Q}_{\ell})$ lies in a maximal torus
$\mathfrak{T}(s)$ of $M_{\ell}^{\mathrm{Zar}}$, as $M_{\ell}^{\mathrm{Zar}}$ is connected. And, $F_s^*$ has eigenvalues
$\chi_i(F_s^{*})$, where $\chi_i$ are the weights
of $H^1(A_{\overline{t}},\mathbf{Q}_{\ell})$ relative to
$\mathfrak{T}(s)$.
\smallskip

These eigenvalues generate over $\mathbf{Q}$ an extension of degree bounded by a constant, as the characteristic polynomial of $F_s^*$ has coefficients in $\mathbf{Z}$ (Weil). Thus,
if some ratio $\chi_i(F_s^{*})/\chi_j(F_s^{*})$ is a root of unity, its order divides an integer $N(g)>1$ depending only on $g=\mathrm{dim}(A_t)$.
\smallskip

The following subset of $M_{\ell}$
\[\{u\in M_{\ell}, u^{N(g)}\ \mathrm{has\ all\ distinct\ eigenvalues\
on}\ H^1(A_{\overline{t}},\mathbf{Q}_{\ell})\}\] is Zariski open and stable under conjugation. Its volume in the normalized Haar measure of $M_{\ell}$ is by Cebotarev's density theorem (\cite{serre_adic}, I--8, Corollary 2) the density of the set
\[\{s\in S\backslash\{t\}, (F_s^{*})^{N(g)}\ \mathrm{has\ all\ distinct\ eigenvalues}\}\] or the density of the set
\[\Sigma=\{s\in S\backslash\{t\}, (F_s^*)^N\ \mathrm{has\ all\ distinct\ eigenvalues}, \forall\ N\geq 1\}.\]

This volume and this density are $1$ because the characters $\chi_i$ are all distinct by assumption.
\smallskip

Consider an integer $N\geq 1$ and a finite extension $k'$ of $k(s)$ of degree $N$, where $s\in\Sigma$. Put $s'=\mathrm{Spec}(k')$ and $A_{s'}=A_s\times_ss'$. As $(F^*_s)^N$ has all distinct eigenvalues on $H^1(A_{\overline{t}}, \mathbf{Q}_{\ell})$, the ring
$\mathrm{End}_{s'}(A_{s'})$ is commutative, for by Tate
\[
\mathrm{End}_{s'}(A_{s'})\otimes_{\mathbf{Z}}\mathbf{Q}_{\ell}\
\widetilde{\to}\
\mathrm{End}_{(F^*_s)^N}(H^1(A_{\overline{t}},\mathbf{Q}_{\ell}))^{\mathrm{opposite}}.\]

Now $A_s$ is isogenous to a product of simple abelian varieties
$A_i$, $i\in I$. If one factor appears with multiplicity $>1$, or if
$A_i\times_s s'$ is not simple, or if $A_i\times_s s'$ and
$A_j\times_s s'$ are isogenous for $i\neq j$, then
$\mathrm{End}_{s'}(A_{s'})$ is not commutative. So these factors
$A_i$ of $A_s$ are absolutely simple mutually non-isogenous over $\overline{s}$, $\forall\ s\in\Sigma$.

\end{proof}

And, this weight multiplicity free condition means (\cite{howe}, 4.6.3) that the tensor components of each $V_{\mathfrak{l}}$, as a $E_{\mathfrak{l}}$-linear representation of the derived group of $M_{\mathfrak{l}}^{\mathrm{Zar}}$, are
\smallskip

--- either minuscule
\smallskip

--- or of the types $(A_n, r\omega_1)$, $(A_n, r\omega_n)$, $(B_n, \omega_1)$, $(C_3, \omega_3)$, $(G_2, \omega_1)$ for some integers $n, r>1$.
\smallskip

Recall that a minuscule representation is a highest weight representation all whose weights have the same length.
\smallskip

To seek a positive answer we now assume that some $V_{\mathfrak{l}}$ is minuscule. (The types $(A_n, r\omega_1)$, $(A_n, r\omega_n)$ are not self-dual and thus do not occur in $V_{\mathfrak{l}}$ if $E$ is totally real. The non-minuscule types might after all be ruled out by elementary means.)
\smallskip

We assume that $V_{\mathfrak{l}}$ is even \emph{elliptic minuscule}, namely, that the derived group $G_{\mathfrak{l}}$ of $M_{\mathfrak{l}}^{\mathrm{Zar}}$ admits at least one maximal torus acting irreducibly on $V_{\mathfrak{l}}$.
Such a torus has a nonempty Zariski open set of $E_{\mathfrak{l}}$-points acting irreducibly on $V_{\mathfrak{l}}$.
\smallskip

The subset of the compact analytic group $M_{\mathfrak{l}}$ consisting of those elements acting irreducibly on $V_{\mathfrak{l}}$ is a union of conjugacy classes and is open by Krasner's lemma (\cite{serre_corps locaux}, II, Exercice 2).
For elliptic minuscule $V_{\mathfrak{l}}$, this subset is nonempty whose nonzero volume in the normalized Haar measure of $M_{\mathfrak{l}}$ is by Cebotarev's density theorem the density of the set
\[\{s\in S\backslash\{t\}, F_s^*\ \mathrm{acts\ irreducibly\ on}\ V_{\mathfrak{l}}\},\] or equivalently the density of the set
\[\{s\in S\backslash\{t\}, (\mathrm{End}_s(A_s)\otimes_{\mathbf{Z}}\mathbf{Q})\otimes_EE_{\mathfrak{l}}\ \mathrm{is\ a\ division\ algebra}\},\] which is $\leq$ the density of the set
\[\{s\in S\backslash\{t\}, \mathrm{End}_s(A_s)\otimes_{\mathbf{Z}}\mathbf{Q}\ \mathrm{is\ a\ division\ algebra}\},\] or that of
\[\{s\in S\backslash\{t\}, k(s)\ \mathrm{is\ a\ prime\ field}, A_s\ \mathrm{is\ simple}\}.\]

So one has the following partial answer :
\smallskip

\smallskip

{\bf Theorem 1.5.} \emph{Let $\ell$ be a prime number. Suppose that $E:=\mathrm{End}_t(A_t)\otimes_{\mathbf{Z}}\mathbf{Q}=\mathrm{End}_{\overline{t}}(A_{\overline{t}})\otimes_{\mathbf{Z}}\mathbf{Q}$ is a field, that $M_{\ell}^{\mathrm{Zar}}$ is connected and that the monodromy representation $H^1(A_{\overline{t}}, \mathbf{Q}_{\ell})$ admits an elliptic minuscule factor $V_{\mathfrak{l}}$ for a place $\mathfrak{l}$ of $E$ above $\ell$.}
\smallskip

\emph{Then, for every prime $l$, $H^1(A_{\overline{t}}, \mathbf{Q}_l)$ has no multiple weights as a representation of the identity component of $M_l^{\mathrm{Zar}}$, and $A_t$ specializes to absolutely simple abelian varieties at a set of places of positive density.}

\smallskip

\smallskip

To provide substance to this answer, our goal is to classify elliptic minuscule representations, namely, to solve the problem below :
\smallskip

\smallskip

{\bf Question 1.6.} \emph{Let $G$ be a semi-simple algebraic group over the spectrum
$\eta$ of a finite extension of $\mathbf{Q}_{\ell}$ and $\rho_V:
G\to\mathrm{GL}(V)$ an absolutely irreducible $\eta$-linear algebraic
representation with finite kernel. Does $G$ admit some maximal torus acting irreducibly on $V$}?

\smallskip

One can assume $G$ to be simply connected. Let $\overline{\eta}$ be a geometric point of $\eta$. Notice that a maximal torus $\mathfrak{T}$ acts
irreducibly on $V$ if and only if the weights of
$V_{\overline{\eta}}$ relative to $\mathfrak{T}_{\overline{\eta}}$
are permuted transitively by $\pi_1(\eta,\overline{\eta})$. So if
such a torus exists, all the weights have the same length, that is,
$V_{\overline{\eta}}$ is minuscule.
\smallskip

Let $D_{\overline{\eta}}$ be the Dynkin diagram of
$G_{\overline{\eta}}$ and
$\rho_D:\pi_1(\eta,\overline{\eta})\to\mathrm{Aut}(D_{\overline{\eta}})$
the index. Let $\alpha_i$, $i=1,\cdots,r$, be the
$\pi_1(\eta,\overline{\eta})$-orbits in $D_{\overline{\eta}}$
consisting of minuscule vertices corresponding to a minuscule
representation $V=V_1\otimes_{\eta}\cdots\otimes_{\eta}V_r$ of
$G=G_1\times_{\eta}\cdots\times_{\eta} G_r$, $G_i$ being the simple
factors. Put $D=(D_{\overline{\eta}},\rho_D)$,
$\alpha_V=\sum\alpha_i$.
\smallskip

Whether or not $G$ has a maximal torus acting irreducibly on $V$
depends in fact only on $(D,\alpha_V)$ (2.3, 3.1). If $G$
admits such a torus, we call $(D,\alpha_V)$ an \emph{elliptic minuscule
pair} (2.2). The elliptic minuscule pairs with connected Dynkin
diagrams are enumerated in (3.2).
\smallskip

{\bf Remark 1.7.} Suppose that $\overline{t}$ has values in $\mathbf{C}$, that $\mathrm{End}_{\overline{t}}(A_{\overline{t}})\otimes_{\mathbf{Z}}\mathbf{Q}$ is a field and that the Mumford--Tate group of the Hodge structure on $H^1(A_{\overline{t}}^{an},\mathbf{Q})$ is definable by absolute Hodge cycles rational over $t$ (\cite{deligne_absolute}, 2.11, 2.9). It is possible that then $A_t$ has absolutely simple specializations at a set of places of positive density.

\smallskip

\emph{Acknowledgement ---} We owe much to M. Aschbacher for answering questions of Ying Zong. We thank J. Achter for kind correspondences. We thank especially C. Chai and
S. Lu for their critical reading of the manuscript and providing
comprehensive advice on writings.

\smallskip

\smallskip

\smallskip


\section{Elliptic minuscule pairs}
\smallskip

2.1. A Dynkin diagram is a finite set $D$ equipped with the
structure of a function $l:D\to \{1,2,3\}$ (``longueurs'') and of a
binary relation $L$ (``liaisons'') on $D$ such that $L$ is disjoint
with the diagonal of $D\times D$.
\smallskip

Every root system has its Dynkin diagram with its connected components
labeled according to types as $A, B,\cdots, G_2$ (\cite{bourbaki},
Chapitre VI, Th\'{e}or\`{e}me 3, p. 197).
\smallskip

Let $S$ be a scheme. An $S$-Dynkin diagram is a sheaf of sets $D$ on
$S$ for the \'{e}tale topology which is locally constant constructible and
is equipped with the structure of a morphism $l:D\to \{1,2,3\}_S$ and
of a sheaf of $S$-relations $L\subset D\times_S D$, $L$ locally
constant constructible on $S$, such that for every geometric point
$s$ of $S$ the fibre $D_s$ with the function $l_s$ and the
relation $L_s$ is a Dynkin diagram.
\smallskip

For every $S$-scheme $S'$, $D\times_SS'$ is an $S'$-Dynkin diagram
and every descent datum on $D$ relative to $S$ for the \'{e}tale
topology is effective.
\smallskip

The monodromy representation
\[
\rho_{D,s}: \pi_1(S,s)\to \mathrm{Aut}(D_s,l_s,L_s)
\]
associated with an $S$-Dynkin diagram $D$ at a geometric point $s\to S$ is said to be the
\emph{index} of $D$ at $s$ (cf. \cite{tits}, 2.3).
\smallskip

One defines $\pi_0(D)$ to be the quotient of $D$ by the equivalence
relation generated by $L$. Notice that $D$ is a $\pi_0(D)$-Dynkin diagram.
\smallskip

Every reductive $S$-group scheme has its $S$-Dynkin diagram which is
functorial with respect to isomorphisms and is compatible with every base
change (SGA 3, \'{E}xpos\'{e} XXIV, 3.3).
\smallskip

Given an $S$-Dynkin diagram $D$, if at every geometric point $s$ of
$S$ the components of the fibre $D_s$ are of the types $A,B,\cdots,
G_2$, then there is a quasi-\'{e}pingl\'{e} semi-simple simply
connected $S$-group scheme which has $D$ as its $S$-Dynkin diagram (SGA
3, \'{E}xpos\'{e} XXIV, Th\'{e}or\`{e}me 3.11).
\smallskip

And, for each semi-simple simply connected $S$-group scheme $G$,
there exists up to unique isomorphisms a unique pair $(Q,u)$
which consists of a quasi-\'{e}pingl\'{e} semi-simple simply connected
$S$-group scheme $Q$ and of an ``isomorphisme ext\'{e}rieur'' $u\in
\mathrm{Isom.ext}_S(Q,G)$ (SGA 3, \'{E}xpos\'{e} XXIV, Corollaire
3.12). The existence of $u$ enables the identification of the
$S$-Dynkin diagram $D$ of $Q$ with that of $G$ and permits one to
define the $S$-scheme of ``isomorphismes int\'{e}rieurs''
\begin{displaymath}
\underline{\mathrm{Isom.int}}_S (Q,G),
\end{displaymath}
which is a left torsor under the adjoint group of $G$ and a right
torsor under the adjoint group of $Q$.
\smallskip

Let $T\subset B$ be the canonical maximal torus and Borel subgroup
of $Q$, $U$ the unipotent radical of $B$, $N$ the normalizer
of $T$ in $Q$ and $W=N/T$ the Weyl group. Let
\[
\pi: X\to S
\]
denote the $S$-scheme $Q/B$, which is projective smooth with
geometrically connected fibres over $S$.
\smallskip

Suppose that
\[
\omega: T\to \mathbf{G}_{m,S}
\]
is a weight of $Q$ with respect to $T$ that is dominant relative to the
notion of positivity defined by $B$. Let
\[
\omega_B: B\to B/U=T\stackrel{\omega}{\longrightarrow}
\mathbf{G}_{m,S}
\]
be the composition. This character $\omega_B$, when twisted by the $B_X$-torsor
\[
Q\to Q/B=X,
\]
provides a $\mathbf{G}_{m,X}$-torsor
\[
Q\stackrel{B_X}{\wedge} \mathbf{G}_{m,X}
\]
and an invertible $\mathcal{O}_X$-module
\[
L_{\omega}=Q\stackrel{B_X}{\wedge}\mathbf{G}_{m,X}\stackrel{\mathbf{G}_{m,X}}{\wedge}\mathcal{O}_X.
\]

Recall that $E_{\omega}=\pi_{*}L_{\omega}$ is a representation of
$Q$ on a locally free $\mathcal{O}_S$-module of finite rank whose
formation is compatible with every base change $S'\to S$. And when $S$ is
the spectrum of an algebraically closed field of characteristic
zero, $E_{\omega}$ is irreducible with highest weight $\omega$.
\smallskip

In particular, to each section $\alpha\in D(S)$ of the $S$-Dynkin
diagram $D$, there corresponds a fundamental representation
$E_{\alpha}$ of $Q$ of fundamental weight $\omega_{\alpha}$.
\smallskip

We say that a section $\alpha\in D(S)$ is \emph{minuscule} if the Weyl orbit
\[
W\omega_{\alpha}\subset\underline{\mathrm{Hom}}_S(T,\mathbf{G}_{m,S})\]
is the sheaf of weights of $E_{\alpha}$ relative to $T$.
\smallskip

More generally, $\alpha=\sum_{i=1}^r\alpha_i$, $\alpha_i\in D(S)$,
is said to be \emph{minuscule} if each $\alpha_i$ is minuscule and if, for every
geometric point $s$ of $S$, $\alpha_{i,s}$ lie in distinct
components of $D_s$. Let
$W\omega_{\alpha}:=W\omega_{\alpha_1}\times_S\cdots\times_S
W\omega_{\alpha_r}$.
\smallskip

\smallskip

{\bf Definition 2.2.} \emph{Suppose that $S$ is connected and
that $\alpha=\sum^r_{i=1} \alpha_i$ is minuscule. The pair $(D,\alpha)$ is
said to be an elliptic minuscule pair or simply
elliptic if there exists a $W$-torsor $x$ on $S$ such that
\[
x\stackrel{W}{\wedge} W\omega_{\alpha}
\]
is a connected object in the Galois category of locally constant
constructible sheaves on $S$, that is, if at some geometric point $s$ of
$S$ the image of the monodromy representation
\[
\rho_{x,s}: \pi_1(S,s)\to \mathrm{Aut}((x\stackrel{W}{\wedge}
W\omega_{\alpha})_s)
\]
acts transitively on the fibre $(x\stackrel{W}{\wedge} W\omega_{\alpha})_s$. Every
such $W$-torsor $x$ is said to be elliptic for $(D,\alpha)$.}
\smallskip

\smallskip

One has the following result :
\smallskip

\smallskip

{\bf Theorem 2.3.} \emph{Let $\eta$ be the spectrum of a complete
discretely valued field of characteristic zero with finite residue
field. Let $G$ be a semi-simple algebraic group over $\eta$ with Dynkin
diagram $D$ and let $\rho_V: G\to \mathrm{GL}(V)$ be an absolutely irreducible
$\eta$-linear algebraic representation with finite kernel.}
\smallskip

\emph{Then there exists a maximal torus of $G$
acting irreducibly on $V$ if and only if $V$ is minuscule and
$(D,\alpha)$ is elliptic, $\alpha$ being the minuscule section
corresponding to $V$.}
\smallskip

\smallskip

For its proof, we may and do assume $G$ to be simply connected.
\smallskip

Observe that if $G$ admits a maximal torus $\mathfrak{T}$ which acts irreducibly on $V$,
then the weights of $V_{\overline{\eta}}$ relative to
$\mathfrak{T}_{\overline{\eta}}$ are permuted transitively by
$\pi_1(\eta,\overline{\eta})$. A priori, all these weights have the
same length, and so $V$ is minuscule (\cite{bourbaki}, Chapitre VIII,
\S 7, Proposition 6, p. 127).
\smallskip

In the following we suppose that $V$ is minuscule. Let $\alpha=\sum \alpha_i$ denote its corresponding minuscule section of $D$.
\smallskip

\smallskip

{\bf Lemma 2.4.} \emph{For each anisotropic maximal torus
$\mathfrak{T}$ of $G$, if $\mathfrak{T}^{\mathrm{ad}}$ denotes its image in
the adjoint group $G^{\mathrm{ad}}$, the canonical map
\begin{displaymath}
H^1(\eta,\mathfrak{T}^{\mathrm{ad}})\to H^1(\eta,G^{\mathrm{ad}})
\end{displaymath}
is surjective and $H^2(\eta,\mathfrak{T})=0$.}

\begin{proof} Notice that
$H^1(\eta,G)=0$, as $G$ is by assumption simply connected (Kneser). Let $Z$ be the center of $G$. The central extension
\[
1\to Z\to G\to G^{\mathrm{ad}}\to 1
\]
induces the cohomology sequence
\[
H^1(\eta,G)\to
H^1(\eta,G^{\mathrm{ad}})\stackrel{\partial}{\longrightarrow}
H^2(\eta,Z)\] from which it follows that
\[
\partial: H^1(\eta,G^{\mathrm{ad}})\to H^2(\eta,Z)
\]
is injective. To show that
\[
H^1(\eta,\mathfrak{T}^{\mathrm{ad}})\to H^1(\eta,G^{\mathrm{ad}})
\]
is surjective, it suffices to show that the composition
\[
\delta: H^1(\eta,\mathfrak{T}^{\mathrm{ad}})\to
H^1(\eta,G^{\mathrm{ad}})\stackrel{\partial}{\hookrightarrow}
H^2(\eta,Z)
\]
is surjective. The map
\[
\delta: H^1(\eta,\mathfrak{T}^{\mathrm{ad}})\to H^2(\eta,Z)
\]
is a coboundary map induced by the central extension
\[
1\to Z\to \mathfrak{T}\to \mathfrak{T}^{\mathrm{ad}}\to 1
\]
and the cohomology sequence
\[
H^1(\eta,\mathfrak{T}^{\mathrm{ad}})\stackrel{\delta}{\longrightarrow}
H^2(\eta,Z)\to H^2(\eta,\mathfrak{T})
\]
implies that
\[
\delta: H^1(\eta,\mathfrak{T}^{\mathrm{ad}})\to H^2(\eta,Z)
\]
is surjective if $H^2(\eta,\mathfrak{T})=0$.
\smallskip

So it remains to show that $H^2(\eta, \mathfrak{T})=0$. Since the Yoneda pairing
\[
\mathrm{Hom}_\eta(\mathfrak{T},\mathbf{G}_m)\times
H^2(\eta,\mathfrak{T})\to H^2(\eta,\mathbf{G}_m)=\mathrm{Br}(\eta)\
\widetilde{\to}\ \mathbf{Q}/\mathbf{Z}
\]
is non-degenerate (Nakayama--Tate), it suffices to show that
\[
\mathrm{Hom}_\eta(\mathfrak{T},\mathbf{G}_m)=0.
\]

But this latter is precisely the condition that $\mathfrak{T}$ is anisotropic.

\end{proof}
\smallskip

Let the quasi-\'{e}pingl\'{e} semi-simple simply connected
$\eta$-group scheme $Q$, the ``isomorphisme ext\'{e}rieur'' $u\in
\mathrm{Isom.ext}_{\eta}(Q,G)$ and the bitorsor
$\underline{\mathrm{Isom.int}}_{\eta}(Q,G)$ be as in (2.1).
\smallskip

Let $T\subset B$ be the canonical maximal torus and Borel subgroup
of $Q$, $N$ the normalizer of $T$ in $Q$, $W=N/T$, $C$ the center of
$Q$ and $T^{\mathrm{ad}}$ (resp. $N^{\mathrm{ad}}$) the image of $T$
(resp. $N$) in the adjoint group $Q^{\mathrm{ad}}$.
\smallskip

Let $E_{\alpha}=\otimes E_{\alpha_i}$ be the minuscule
representation of $Q$ of fundamental weight $\omega_{\alpha}$.
\smallskip

\smallskip

{\bf Lemma 2.5.} 1) \emph{The $Q^{\mathrm{ad}}(\eta)$-conjugacy
classes of maximal tori of $Q$ are in bijective correspondence with
the elements of $H^1(\eta,N)$.}
\smallskip

2) \emph{The map $H^1(\eta,N)\to H^1(\eta,W)$ is injective whose
image contains those isomorphism classes of $W$-torsors $x$ on $\eta$
such that $x\stackrel{W}{\wedge}T$ is anisotropic.}

\begin{proof} 1) The set $(Q/N)(\eta)$ classifies the maximal tori
of $Q$ because locally on $\eta$ for the \'{e}tale topology they are
all conjugate to $T$ by sections of $Q$.
\smallskip

The exact sequence of pointed sets
\[
Q^{\mathrm{ad}}(\eta)\to (Q/N)(\eta)\to H^1(\eta,N^{\mathrm{ad}})\to
H^1(\eta,Q^{\mathrm{ad}})
\]
shows that the $Q^{\mathrm{ad}}(\eta)$-orbits in $(Q/N)(\eta)$ are
in one-to-one correspondence with the elements of the kernel of the map
\[H^1(\eta,N^{\mathrm{ad}})\to H^1(\eta,Q^{\mathrm{ad}}).\]

Observe that in the cohomology sequence
\[
H^1(\eta,Q)\to
H^1(\eta,Q^{\mathrm{ad}})\stackrel{\partial}{\longrightarrow}
H^2(\eta,C)
\]
induced by the central extension
\[
1\to C\to Q\to Q^{\mathrm{ad}}\to 1,
\]
the map
\[
\partial: H^1(\eta,Q^{\mathrm{ad}})\to H^2(\eta,C)
\]
is injective since
\[
H^1(\eta,Q)=0,
\]
$Q$ being simply connected.
\smallskip

Hence, the kernel of the map
\[
H^1(\eta,N^{\mathrm{ad}})\to H^1(\eta,Q^{\mathrm{ad}})
\]
is equal to the kernel of the composition
\[
\delta: H^1(\eta,N^{\mathrm{ad}})\to
H^1(\eta,Q^{\mathrm{ad}})\stackrel{\partial}{\hookrightarrow}
H^2(\eta,C).
\]

This
\[
\delta: H^1(\eta,N^{\mathrm{ad}})\to H^2(\eta,C)
\]
is a coboundary map induced by the central extension
\[
1\to C\to N\to N^{\mathrm{ad}}\to 1.
\]

From the exact sequence
\[
H^1(\eta,C)\to H^1(\eta,N)\to
H^1(\eta,N^{\mathrm{ad}})\stackrel{\delta}{\longrightarrow}
H^2(\eta,C),
\]
one finds that $H^1(\eta, N)$ is mapped onto $\mathrm{Ker}(\delta)$ by
\[
H^1(\eta,N)\to H^1(\eta,N^{\mathrm{ad}}).
\]

To conclude that $H^1(\eta,N)$ is isomorphic to this image, it suffices
to show that the map
\[
H^1(\eta,C)\to H^1(\eta,N)
\]
is $0$ or, by the factorization
\[
H^1(\eta,C)\to H^1(\eta,T)\to H^1(\eta,N),
\]
that
\[
H^1(\eta,T)=0.
\]

This latter vanishing follows from the identity
\[
H^1(\eta,T)=H^1(D,\mathbf{G}_m)
\]
(SGA 3, \'{E}xpos\'{e} XXIV, Corollaire 3.14) and by Satz 90 :
\[
H^1(D,\mathbf{G}_m)=0,
\]
the Dynkin diagram $D$ being representable by a finite
\'{e}tale $\eta$-scheme.
\smallskip

\smallskip

2) That
\begin{displaymath}
H^1(\eta,N)\to H^1(\eta,W)
\end{displaymath}
is injective results from the cohomology sequence
\begin{displaymath}
H^1(\eta,T)\to H^1(\eta,N)\to H^1(\eta,W)
\end{displaymath}
and by $H^1(\eta,T)=0$.
\smallskip

The class of a $W$-torsor $x$ on $\eta$ lies in the image of the map
\begin{displaymath}
H^1(\eta,N)\to H^1(\eta,W)
\end{displaymath}
if and only if an obstruction
\begin{displaymath}
o(x)\in H^2(\eta,x\stackrel{W}{\wedge}T)
\end{displaymath}
vanishes.
\smallskip

When $x\stackrel{W}{\wedge}T$ is anisotropic, one has in fact $H^2(\eta,x\stackrel{W}{\wedge}T)=0$ (2.4).

\end{proof}
\smallskip

{\bf Lemma 2.6.} \emph{If a torus of $G$ acts irreducibly on
$V$, it is anisotropic.}

\begin{proof} A torus is anisotropic if and only if it has no
diagonalizable sub-torus other than $1$.
\smallskip

Recall that the kernel of the representation
\begin{displaymath}
\rho_V: G\to \mathrm{GL}(V)
\end{displaymath}
is finite. And $\mathrm{det}(\rho_V)=1$, as $G$ is semi-simple.
\smallskip

Suppose that a certain torus of $G$ acts irreducibly on $V$. If a $\mathbf{G}_m$ were in this torus, it would act on $V$ by a
character $z\mapsto z^n$ for some integer $n$ and thus on
$\mathrm{det}(V)$ by the character $z\mapsto z^{nd}$, where
$d=\mathrm{dim}(V)$. So $nd=0$, i.e., $n=0$ and $\mathbf{G}_m$ was
contained in $\mathrm{Ker}(\rho_V)$.

\end{proof}
\smallskip

{\bf Lemma 2.7.} \emph{The group $G$ has a maximal torus acting
irreducibly on $V$ if and only if the group $Q$ has a maximal torus
acting irreducibly on $E_{\alpha}$.}

\begin{proof} Suppose that a maximal torus $\mathfrak{T}$ of $G$ acts
irreducibly on $V$. By (2.6), $\mathfrak{T}$ is anisotropic. And by (2.4), the map
\[
H^1(\eta,\mathfrak{T}^{\mathrm{ad}})\to H^1(\eta,G^{\mathrm{ad}})
\]
is surjective. The $G^{\mathrm{ad}}$-torsor
\[
\underline{\mathrm{Isom.int}}_{\eta}(Q,G)
\]
is in particular the image of a $\mathfrak{T}^{\mathrm{ad}}$-torsor,
which means (SGA 3, \'{E}xpos\'{e} XXIV, Proposition 2.11) that
$\mathfrak{T}$ imbeds into $Q$ as a maximal torus and the scheme
\[
\mathfrak{I}=\underline{\mathrm{Isom.int}}_{\eta}(Q,G; \mathrm{Id\
on}\ \mathfrak{T})
\]
of ``isomorphismes int\'{e}rieurs'' from $Q$ to $G$ that induce the
identity automorphism on $\mathfrak{T}$ is nonempty.
\smallskip

Let $\overline{\eta}$ be a geometric point of $\eta$. The choice of
a section $\iota\in \mathfrak{I}(\overline{\eta})$ identifies the
sheaves of weights of $V$ and of $E_{\alpha}$ relative to
$\mathfrak{T}$. So $E_{\alpha}$ is isomorphic to $V$ as a
$\mathfrak{T}$-module. So $\mathfrak{T}$ acts irreducibly on $E_{\alpha}$.
\smallskip

The other direction is proven similarly.

\end{proof}
\smallskip

2.8. \emph{Proof of Theorem} 2.3.
\smallskip

\smallskip

By (2.7) it suffices to show that $(D,\alpha)$ is elliptic if and
only if $Q$ has some maximal torus acting irreducibly on $E_{\alpha}$.
\smallskip

Suppose first that $Q$ admits a maximal torus acting irreducibly on
$E_{\alpha}$.
\smallskip

This torus has then the form $z\stackrel{N}{\wedge}T$ for an $N$-torsor $z$ (2.5).
Relative to this torus the sheaf of weights of $E_{\alpha}$ is
\begin{displaymath}
z\stackrel{N}{\wedge}W\omega_{\alpha}\subset
z\stackrel{N}{\wedge}\underline{\mathrm{Hom}}_{\eta}(T,\mathbf{G}_m).
\end{displaymath}

The condition that $z\stackrel{N}{\wedge}T$ acts irreducibly on $E_{\alpha}$ is
equivalent to the condition that $z\stackrel{N}{\wedge}W\omega_{\alpha}$ is a connected object
in the Galois category of locally constant constructible sheaves on
$\eta$. So $z\stackrel{N}{\wedge}W$ is a $W$-torsor elliptic for $(D,\alpha)$.
\smallskip

Suppose next that $(D,\alpha)$ is elliptic and that $x$ is a
$W$-torsor elliptic for $(D, \alpha)$.
\smallskip

Let $\rho: Q\to \mathrm{GL}(E_{\alpha})$ denote the minuscule representation corresponding to $\alpha$ and let
$\rho_T$ be its restriction to $T$.
\smallskip

One has that $\mathrm{Ker}(\rho_T)$
is finite and that $\mathrm{det}(\rho_T)=1$. The torsor $x$ twists $\rho_T$ to a representation of $x\stackrel{W}{\wedge}T$,
\[
\rho_{x,T}: x\stackrel{W}{\wedge}T\to \mathrm{GL}(E_{\alpha}),
\]
which has $x\stackrel{W}{\wedge}W\omega_{\alpha}$ as its sheaf of weights.
In particular, $\rho_{x,T}$ is irreducible. Moreover, being a twist of $\rho_T$, $\rho_{x, T}$
has finite kernel and determinant $1$. As in (2.6), $x\stackrel{W}{\wedge}T$
is anisotropic. Thus it can be imbedded into $Q$ (2.5). So $x\stackrel{W}{\wedge}T$ is a sought-after maximal torus of $Q$ acting irreducibly on $E_{\alpha}$.

\begin{flushright}
$\square$
\end{flushright}

\smallskip

\smallskip

\smallskip


\section{Simple elliptic pairs}
\smallskip

\smallskip

Let $S$ be a scheme. Recall that an $S$-Dynkin diagram $D$ is also a $\pi_0(D)$-Dynkin diagram, where $\pi_0(D)$ is the finite \'{e}tale $S$-scheme, the quotient of $D$ by the $S$-equivalence relation generated by the $S$-binary relation $L$ (``liaisons'') (2.1). The fiber $D\times_{\pi_0(D)}z$ is a connected Dynkin diagram for every geometric point $z$ of $\pi_0(D)$.
\smallskip

Suppose that $S$ is connected. Let $(D, \alpha)$ be as in (2.2). Suppose that $\pi_0(D)=S$. Then in the notations of Bourbaki--Tits (\cite{bourbaki}, Chapitre VI, Planches I--IX, p. 250--275, and
\cite{tits}, p. 54--61), if $D$ is non-constant, $(D,\alpha)$ can
only be $({}^2A_n,\alpha_{\frac{n+1}{2}})$, $n$ odd $\geq 3$, or
$({}^2D_n,\alpha_1)$, $n\geq 5$, or $({}^2D_4,\alpha_i)$, $i=1,3,4$.
\smallskip

Let $s$ be a geometric point of $S$. We write down the condition
that $(D,\alpha)$ be elliptic.
\smallskip

\smallskip

{\bf Lemma 3.1.} 1) \emph{$(A_n,\alpha_r)$, $r\in [1,n]$, is
elliptic if and only if there is a monodromy representation in the symmetric group of $n+1$ letters
\[
\rho: \pi_1(S,s)\to \mathfrak{S}_{n+1}
\]
whose image permutes transitively the subsets of $\{1,\cdots,n+1\}$
of cardinality $r$.}
\smallskip

2) \emph{$(B_n,\alpha_n)$ is elliptic if and only if there is a
representation
\[
\rho: \pi_1(S,s)\to \mathrm{GL}_n(\mathbf{Z})
\]
whose image lies in the group generated by the diagonal matrices and
monomial matrices and acts transitively on the set
\[
\{\pm e_1\pm \cdots \pm e_n\},
\]
where $e_1,\cdots,e_n$ denote the standard basis of $\mathbf{Z}^n$.}
\smallskip

3) \emph{$(C_n,\alpha_1)$ is elliptic if and only if there is a
representation
\[
\rho: \pi_1(S,s)\to \mathrm{GL}_n(\mathbf{Z})
\]
whose image lies in the group generated by the diagonal matrices and
monomial matrices and acts transitively on the set
\[
\{e_1,\cdots,e_n,-e_1,\cdots,-e_n\},
\]
where $e_1,\cdots,e_n$ denote the standard basis of $\mathbf{Z}^n$.}
\smallskip

4) \emph{$(D_n,\alpha_1)$ is elliptic if and only if there is a
representation
\[
\rho: \pi_1(S,s)\to \mathrm{GL}_n(\mathbf{Z})
\]
whose image lies in the group generated by the diagonal matrices of
determinant $1$ and monomial matrices and acts transitively on the set
\[
\{e_1,\cdots,e_n,-e_1,\cdots,-e_n\},
\]
where $e_1,\cdots,e_n$ denote the standard basis of $\mathbf{Z}^n$.}
\smallskip

5) \emph{$(D_n,\alpha_{n-1})$ (resp. $(D_n,\alpha_n)$) is elliptic
if and only if there is a representation
\[
\rho: \pi_1(S,s)\to \mathrm{GL}_n(\mathbf{Z})
\]
whose image lies in the group generated by the diagonal matrices of
determinant $1$ and monomial matrices and permutes transitively the
vectors
\[
s_1e_1+\cdots+s_ne_n,
\]
where $s_i\in \{1,-1\}$, $s_1\cdots s_n=-1$ (resp. $s_1\cdots
s_n=1$) and $e_1,\cdots,e_n$ denote the standard basis of
$\mathbf{Z}^n$.}
\smallskip

6) \emph{$(E_6,\alpha_i)$, $i=1,6$, are elliptic if and only if there
is a representation
\[
\rho: \pi_1(S,s)\to \mathrm{O}(\mathbf{F}_2^6,q)
\]
whose image permutes transitively the nonzero $q$-singular vectors in
$\mathbf{F}_2^6$, where $q$ is the quadratic form such that
\[ q(e_i)=q(f_j)=1,\ q(e_i+e_j)=q(f_i+f_j)=0,\
q(e_i+f_j)=\delta_{ij},\] where $e_i, f_j$, $1\leq i,j\leq 3$, are a
basis of $\mathbf{F}_2^6$ and where $\delta_{ij}=1$, if $i=j$, and $\delta_{ij}=0$, if
$i\neq j$.}
\smallskip

7) \emph{$(E_7,\alpha_7)$ is elliptic if and only if there is a
representation
\[
\rho: \pi_1(S,s)\to \{1,-1\}\times \mathrm{Sp}_6(\mathbf{F}_2)
\]
whose image acts transitively on $\{1,-1\}\times
(\mathrm{Sp}_6(\mathbf{F}_2)/\mathrm{O}(q))$, $q$ being the
quadratic form on $\mathbf{F}_2^6$ such that
\[
q(e_i)=q(f_j)=1,\ q(e_i+e_j)=q(f_i+f_j)=0,\ q(e_i+f_j)=\delta_{ij},
\]
where $e_i, f_j$ are the standard symplectic base of
$\mathbf{F}_2^6$ and where $\delta_{ij}=1$, if $i=j$, and $\delta_{ij}=0$, if $i\neq j$.}
\smallskip

8) \emph{$({}^{2}A_n,\alpha_{\frac{n+1}{2}})$, $n$ odd $\geq 3$, is
elliptic if and only if there is a representation
\[
\rho=(\rho_1,\rho_2): \pi_1(S,s)\to \{1,-1\}\times
\mathfrak{S}_{n+1}
\]
whose image permutes transitively the subsets of $\{1,\cdots,n+1\}$
of cardinality $(n+1)/2$ and whose component $\rho_1$ is the index
of ${}^2A_n$. Here $-1: Y\mapsto \{1,\cdots,n+1\}\backslash Y$, for
any $Y\subset \{1,\cdots,n+1\}$ of cardinality $(n+1)/2$.}
\smallskip

9) \emph{$({}^{2}D_n,\alpha_1)$, $n\geq 5$, or
$({}^{2}D_n,\alpha_i)$, $n=4$, $i=1,3,4$, are elliptic if and only if there
is a representation
\[
\rho: \pi_1(S,s)\to \mathrm{GL}_n(\mathbf{Z})
\]
whose image
lies in the group $\mathfrak{W}$ generated by the diagonal matrices
and monomial matrices and acts transitively on the set $\{\pm e_1,\cdots,\pm e_n\}$ and which when composed with the projection $\mathfrak{W}\to\mathfrak{W}/\mathfrak{W}_1=\{1,-1\}$ induces the index of ${}^2D_n$ :
\[\rho_{{}^2D_n}:\pi_1(S,s)\stackrel{\rho}{\longrightarrow}\mathfrak{W}\to
\mathfrak{W}/\mathfrak{W}_1=\{1,-1\},\] where $\mathfrak{W}_1$ is the subgroup of $\mathfrak{W}$ generated
by the diagonal matrices of determinant $1$ and monomial matrices and where $e_1,\cdots,e_n$ denote the standard basis of
$\mathbf{Z}^n$.}

\begin{proof} Let $Q$ be a quasi-\'{e}pingl\'{e} semi-simple simply connected $S$-group scheme which has $D$ as its $S$-Dynkin diagram (2.1). Let $T$ be the canonical maximal torus of $Q$. Let $R$ (resp. $W$) be the root system (resp. Weyl group) of $Q$ relative to $T$. One has
the following canonical exact sequence of sheaves of $S$-groups for the \'{e}tale topology :
\[ 1\to
W\to\underline{\mathrm{Aut}}_S(R)\to\underline{\mathrm{Aut}}_S(D)\to
1.\]

This exact sequence induces the cohomology sequence :
\[ H^1(S,W)\to H^1(S,\underline{\mathrm{Aut}}_S(R))\to
H^1(S,\underline{\mathrm{Aut}}_S(D)),\] by which one concludes that
\smallskip

\emph{An $S$-form of $R$, $R_1$, is isomorphic to $x\stackrel{W}{\wedge}R$ for some $W$-torsor $x$ if and only if $R_1$ has its Dynkin diagram isomorphic to $D$.}

\smallskip

\smallskip

When a geometric point $s$ of the connected scheme $S$ is given, the following two conditions are equivalent :
\smallskip

--- \emph{$R_1$ has $D$ as its $S$-Dynkin diagram.}
\smallskip

--- \emph{the composition}
\[
\pi_1(S,s)\stackrel{\rho_{R_1,s}}{\longrightarrow}\mathrm{Aut}(R_s)\to\mathrm{Aut}(D_s)\]
\emph{is the index of $D$ at $s$, where $\rho_{R_1, s}$ denotes the monodromy representation
associated with $R_1$ at $s$.}
\smallskip

Let $x$ be a $W$-torsor and $R^x:=x\stackrel{W}{\wedge}R$. Observe that the monodromy $\mathrm{Im}(\rho_{R^x, s})$ at $s$ associated with every such form $R^x$ normalizes the weights $W_s\omega_{\alpha}$. The following two conditions are equivalent :
\smallskip

--- \emph{$x\stackrel{W}{\wedge}W\omega_{\alpha}$ is a connected object in the Galois category of locally constant constructible sheaves on $S$.}
\smallskip

--- \emph{the monodromy $\mathrm{Im}(\rho_{R^x, s})$ acts transitively on the weights $W_s\omega_{\alpha}$.}
\smallskip

\smallskip

In brief, $(D, \alpha)$ is elliptic if and only if
\smallskip

\emph{There is a representation
\[\rho: \pi_1(S, s)\to\mathrm{Aut}(R_s)\] which satisfies the following two properties }:
\smallskip

\noindent --- \emph{When composed with the projection $\mathrm{Aut}(R_s)\to\mathrm{Aut}(D_s)$ it induces the index of $D$ at $s$ :
\[\rho_{D}: \pi_1(S, s)\stackrel{\rho}{\longrightarrow}\mathrm{Aut}(R_s)\to\mathrm{Aut}(D_s).\]}
\noindent --- \emph{The image of $\rho$ acts transitively on $W_s\omega_{\alpha}$.}

\smallskip

\smallskip

If $D$ is constant, then $W$ and $R$ are constant and the class of a $W$-torsor ``is'' a $W$-conjugacy class of monodromy representations in $W$. This criterion simplifies then to
\smallskip

\emph{There is a representation}
\[\rho: \pi_1(S, s)\to W\] \emph{whose image acts transitively on the weights $W\omega_{\alpha}$.}

\smallskip

\smallskip

For type $(A_n, \alpha_r)$, this says that
\smallskip

\emph{There is a representation}
\[\rho: \pi_1(S, s)\to\mathfrak{S}_{n+1}\] \emph{whose image permutes transitively the subsets of $\{1,\cdots, n+1\}$ of cardinality $r$.}
\smallskip

Indeed, in this case,
\smallskip

--- the Weyl group ``is'' the symmetric group $\mathfrak{S}_{n+1}$ of $n+1$ letters.
\smallskip

--- the Weyl orbit $W\omega_r$ of the minuscule weight $\omega_r$ ``is'' the collection of subsets of $\{1,\cdots, n+1\}$ of cardinality $r$ equipped with its canonical permutation action by $\mathfrak{S}_{n+1}$.

\smallskip

One proceeds similarly for other types provided given a description of $\mathrm{Aut}(R)$, of the Weyl
group $W$, of the minuscule vertex $\alpha$ and of the weights
$W\omega_{\alpha}$.
\smallskip

These for $(B_n,\alpha_n)$,
$(C_n,\alpha_1)$, $(D_n,\alpha_i)$, $i=1,n-1,n$, $(E_6,\alpha_i)$,
$i=1,6$, $({}^2A_n,\alpha_{\frac{n+1}{2}})$, $({}^2D_n,\alpha_1)$,
$({}^2D_4,\alpha_i)$, $i=1,3,4$ follow from Bourbaki \cite{bourbaki},
Chapitre VI, Planches and Chapitre VI, $n^o$4, Exercice 2.
\smallskip

For $(E_7, \alpha_7)$, one can almost quote Bourbaki \cite{bourbaki}, Chapitre VI, $n^o$4, Exercices 3+2 :
\smallskip

Let
$Q(E_7)$ be the root lattice and $P(E_7)$ the weight lattice of a root system of type $E_7$. Then $2P(E_7)\subset
Q(E_7)$ and the quotient $E=Q(E_7)/2P(E_7)$ is a $6$-dimensional
$\mathbf{F}_2$-vector space on which the Killing form $(,)$ induces
a symplectic form. The Weyl group $W(E_7)$ acts on $E$ preserving $(,)$ and it maps onto
$\mathrm{Sp}(E)$ with kernel $\{1,-1\}$ of order $2$, \emph{loc.cit}. The central extension
\[
1\to \{1,-1\}\to W(E_7)\to \mathrm{Sp}(E)\to 1
\]
splits. Let $\{\alpha_1,\cdots, \alpha_7\}$ be a base of $E_7$ so that
$\{\alpha_1,\cdots,\alpha_6\}$ generates a root system of type $E_6$. Observe that the roots of this sub-system
\[
e_1=\alpha_1+\alpha_2+2\alpha_3+2\alpha_4+\alpha_5+\alpha_6,
\]
\[
e_2=\alpha_1+\alpha_2+\alpha_3+\alpha_4+\alpha_5,
\]
\[
e_3=\alpha_2+\alpha_4,
\]
\[
f_1=\alpha_1+\alpha_3+\alpha_4,
\]
\[
f_2=\alpha_4+\alpha_5+\alpha_6,
\]
\[
f_3=\alpha_3+\alpha_4+\alpha_5
\]
satisfy the orthogonality relations
\[
(e_i,e_j)=2\delta_{ij},\ (f_i,f_j)=2\delta_{ij},\
(e_i,f_j)=\delta_{ij}
\]
and that their images in $E$ form a symplectic base. In particular,
\[
F=Q(E_6)/2Q(E_6)\ \widetilde{\to}\ Q(E_7)/2P(E_7)=E
\]
is a bijection, where $Q(E_6)$ denotes the root lattice of $E_6$.
\smallskip

When
$F$ is equipped with the quadratic form $q=\frac{1}{2}(,)$, $W(E_6)$
is identified with $\mathrm{O}(q)$ (\emph{loc.cit.}). Hence,
\[W(E_7)\omega_7=W(E_7)/W(E_6)=\{1,-1\}\times
(\mathrm{Sp}(E)/\mathrm{O}(q)).\]

\end{proof}
\smallskip

\smallskip

It is evident that ellipticity is a nonempty condition only when the base scheme has a rather ``small'' fundamental group.
\smallskip

\smallskip

\smallskip

{\bf Theorem 3.2.} \emph{Let $S$ be the spectra of a complete
discrete valuation ring, $\eta$ (resp. $s$) its generic (resp. closed) point and $\overline{\eta}$ a geometric generic point. Suppose that $k(\eta)$ is of characteristic
zero and that $k(s)$ is finite of characteristic $\ell$.}
\smallskip

\emph{Then the elliptic minuscule pairs $(D,\alpha)$ over $\eta$ such that $D_{\overline{\eta}}$ is a connected Dynkin diagram are enumerated in the following list }:
\smallskip

\smallskip

\noindent \emph{$A$) $(A_n,\alpha_1)$, $(A_n,\alpha_n)$, $n\geq 1$, every prime $\ell$,}
\smallskip

\emph{$(A_{\ell^d-1},\alpha_2)$, $(A_{\ell^d-1},\alpha_{\ell^d-2})$,
$d$ an integer $\geq 1$, every prime $\ell$,}
\smallskip

\emph{$(A_{p-1},\alpha_2)$, $(A_{p-1},\alpha_{p-2})$, $p$ prime, $p\equiv 1$ mod $4$, $\mathrm{Card}(k(s))$ mod $p$ generates
$\mathbf{F}_p^{\times}$,}
\smallskip

\emph{$(A_{p-1}, \alpha_2)$, $(A_{p-1}, \alpha_{p-2})$, $p$ prime, $p\equiv 3$ mod $4$, $\mathrm{Card}(k(s))$ mod $p$ generates a subgroup of $\mathbf{F}_p^{\times}$ of index $\leq 2$,}
\smallskip

\emph{$(A_7,\alpha_3)$, $(A_7,\alpha_5)$, $\ell=2$,}
\smallskip

\emph{$(A_{31},\alpha_3)$, $(A_{31},\alpha_{29})$, $\ell=2$, $5\nmid
[s:\mathbf{F}_2]$};
\smallskip

\smallskip

\noindent \emph{${}^2A$) $({}^2A_3,\alpha_2)$, every prime $\ell$,}
\smallskip

\emph{$({}^2A_5,\alpha_3)$, $\ell=5$,}
\smallskip

\emph{$({}^2A_5, \alpha_3)$, ${}^2A_5$ ramified
over $S$, $\mathrm{Card}(k(s))$ mod $5$ generates
$\mathbf{F}_5^{\times}$};
\smallskip

\smallskip

\noindent \emph{$B$) $(B_3,\alpha_3)$, $(B_4,\alpha_4)$, every prime $\ell$,}
\smallskip

\emph{$(B_n,\alpha_n)$, $n\geq 5$, $\ell=2$};
\smallskip

\smallskip

\noindent \emph{$C$) $(C_n,\alpha_1)$, $n\geq 2$, every prime $\ell$};
\smallskip

\smallskip

\noindent \emph{$D$) $(D_n,\alpha_1)$, $n$ odd $\geq 5$, $\ell=2$,}
\smallskip

\emph{$(D_n, \alpha_1)$, $n$ even $\geq 4$, every prime $\ell$,}
\smallskip

\emph{$(D_5,\alpha_4)$, $(D_5,\alpha_5)$, every prime $\ell$,}
\smallskip

\emph{$(D_n,\alpha_{n-1})$, $(D_n,\alpha_n)$, $n\geq 6$, $\ell=2$};
\smallskip

\smallskip

\noindent \emph{${}^2D$) $({}^2D_n,\alpha_1)$, $n\geq 5$, every prime $\ell$};
\smallskip

\smallskip

\noindent \emph{$E_6$) $(E_6,\alpha_1)$, $(E_6,\alpha_6)$, $\ell=3$,}
\smallskip

\emph{$(E_6, \alpha_1)$, $(E_6, \alpha_6)$,
$\mathrm{Card}(k(s))\equiv\pm 2, \pm 4$ mod $9$};
\smallskip

\smallskip

\noindent \emph{$E_7$) $(E_7,\alpha_7)$, $\ell=2$.}
\smallskip

\smallskip

\smallskip

This list is justified in the remaining sections.
\smallskip

\smallskip

\smallskip


\section{Two lemmas}
\smallskip

\smallskip

Let $S$ be the spectra of a complete discrete valuation ring and $\eta$ (resp. $s$) its generic (resp. closed) point. Suppose that $k(\eta)$ is of characteristic zero and that $k(s)$ is finite of characteristic $\ell$. Let $\overline{\eta}$ (resp. $\overline{s}$) be the spectrum of an algebraic closure of $k(\eta)$ (resp. $k(s)$).
\smallskip

As $S$ is complete along $s$, the inclusion $s\hookrightarrow S$ induces a bijection
\[\pi_1(s, \overline{s})\ \widetilde{\to}\ \pi_1(S, \overline{s}).\]

The group $\pi_1(s, \overline{s})$ is isomorphic to $\widehat{\mathbf{Z}}$ with the Frobenius $F_s$ as its canonical generator. For each integer $N\geq 1$ there is thus up to isomorphisms a unique spectra $S_N$ of a discrete valuation ring such that $S_N$ is finite \'{e}tale Galois over $S$ with cyclic Galois group of order $N$.

\smallskip

Let $S_{(\overline{s})}$ be the strict localization of $S$ at $\overline{s}$ and $\eta^{hs}$ the generic point of $S_{(\overline{s})}$. The open immersion $\eta\hookrightarrow S$ induces a surjection
\[\pi_1(\eta, \overline{\eta})\to \pi_1(S, \overline{\eta})\simeq \pi_1(S, \overline{s}),\] whose kernel, the inertia subgroup of $\pi_1(\eta, \overline{\eta})$, is isomorphic to $\pi_1(\eta^{hs}, \overline{\eta})$. This inertia subgroup admits a canonical surjection
\[\pi_1(\eta^{hs}, \overline{\eta})\to \prod_{p\neq \ell} \mathbf{Z}_p(1),\] which corresponds by Galois theory to the subextension of $k(\overline{\eta})/k(\eta^{hs})$ obtained by joining to $k(\eta^{hs})$ all $N$-th roots of a uniformizer of $S_{(\overline{s})}$ for all integers $N$ prime to $\ell$. The kernel of this surjection, the wild inertia subgroup of $\pi_1(\eta, \overline{\eta})$, is a pro-$\ell$-group and normal in $\pi_1(\eta, \overline{\eta})$.
\smallskip

In particular, the group $\pi_1(\eta, \overline{\eta})$ is pro-solvable.
\smallskip

The quotient of $\pi_1(\eta, \overline{\eta})$ by its wild inertia subgroup is denoted by $\pi_1^t(\eta, \overline{\eta})$, which as a profinite group admits $2$ generators $F, T$ and $1$ single relation :
\[FTF^{-1}=T^q,\] where $q=\mathrm{Card}(k(s))$.

\smallskip

A monodromy representation
$\pi_1(\eta, \overline{\eta})\to \mathfrak{G}$ is said to be unramified (resp. tamely ramified) over $S$ if its kernel contains the inertia (resp. wild inertia) subgroup. A quotient $\mathfrak{G}$ of $\pi_1(\eta, \overline{\eta})$ is said to be unramified (resp. tamely ramified) over $S$ if the quotient homomorphism $\pi_1(\eta, \overline{\eta})\to \mathfrak{G}$ is.

\smallskip

We will apply the following two simple lemmas a few times.

\smallskip

\smallskip

{\bf Lemma 4.1.} \emph{Let $N$ be an integer $\geq 1$. Let
$\zeta\in\mathrm{GL}_N(\mathbf{F}_{\ell})$ be such that
\[ \zeta: e_1\mapsto e_2,\ e_2\mapsto e_3,\ \cdots,\ e_N\mapsto e_1,\] where
$e_1,\cdots,e_N$ denote the standard basis of $\mathbf{F}_{\ell}^N$.}
\smallskip

\emph{Then the semi-direct product
$\langle\zeta\rangle\mathbf{F}_{\ell}^N$ is a quotient of
$\pi_1(\eta,\overline{\eta})$. If $(\ell,N)=1$ and if $V$ is an
irreducible $\mathbf{F}_{\ell}$-linear representation of
$\langle\zeta\rangle$, then $\langle\zeta\rangle V$ is a quotient of
$\pi_1(\eta,\overline{\eta})$.}

\begin{proof} Let $\pi\in\Gamma(S, \mathcal{O}_S)$ be a uniformizer. Let $S'$ be the spectra of a discrete valuation ring such that $S'$ is finite \'{e}tale Galois over $S$ with cyclic Galois group of order $N$. Let $\eta'$ (resp. $s'$) be the generic (resp. closed) point of $S'$, $\zeta$ a generator of $\mathrm{Gal}(S'/S)$ and
let
$u'\in\Gamma(S', \mathcal{O}_{S'})^{\times}$ be a unit such that the images of
$u',\zeta(u'),\cdots, \zeta^{N-1}(u')$ in $k(s')$ form a normal base
over $k(s)$. Then
\[\eta'[x_1,\cdots,x_N]/(x_1^{\ell}-x_1-\zeta(u')\pi^{-1},\cdots,x_N^{\ell}-x_N-\zeta^N(u')\pi^{-1})\] is
connected and Galois over $\eta$ with Galois group
$\langle\zeta\rangle\mathbf{F}_{\ell}^N$. If $(\ell, N)=1$,
$\langle\zeta\rangle V$ is a quotient of
$\langle\zeta\rangle\mathbf{F}_{\ell}^N$ and hence is a quotient of
$\pi_1(\eta,\overline{\eta})$.

\end{proof}
\smallskip

{\bf Lemma 4.2.} \emph{Let $p$ be a prime number different from
$\ell$.}
\smallskip

1) \emph{If the underlying group of an $\mathbf{F}_p$-vector space
$V$ is a normal subgroup of a finite quotient $\mathfrak{G}$ of $\pi_1(\eta, \overline{\eta})$ such that $\mathfrak{G}$
acts irreducibly on $V$ by conjugation, then $\mathrm{dim}\
V=1$.}
\smallskip

2) \emph{There is a unique group of affine linear transformations of
$\mathbf{F}_p$ which contains all translations and which is a quotient
of $\pi_1(\eta,\overline{\eta})$ ramified over $S$. This group has
cardinality $pN$, where $N$ is the order of the element
$\mathrm{Card}(k(s))$ mod $p$ in $\mathbf{F}_p^{\times}$.}

\begin{proof} 1) Let $I$ (resp. $P$) be the image in $\mathfrak{G}$ of the
inertia (resp. wild inertia) subgroup of $\pi_1(\eta, \overline{\eta})$. Notice that $V\cap P=1$.
The intersection $V\cap I$ being normal in
$\mathfrak{G}$ is a sub-$\mathfrak{G}$-module of $V$. As $V$ is by assumption an irreducible $\mathfrak{G}$-module, one has $V\cap I=1$ or $V$.
\smallskip

If $V\cap I=1$, then $V$ is isomorphic to a subgroup of
$\mathfrak{G}/I$ and thus is cyclic.
\smallskip

If $V\cap I=V$, then
$V$ is isomorphic to a subgroup of $I/P$ and thus
is again cyclic.
\smallskip

2) Let $q=\mathrm{Card}(k(s))$. Let $t: x\mapsto x+1$, $\forall\ x\in\mathbf{F}_p$. For every
$a\in\mathbf{F}_p^{\times}$, let $l_a: x\mapsto ax$,
$\forall\ x\in\mathbf{F}_p$. The following relation holds :
\[l_atl_a^{-1}=t^a: x\mapsto x+a,\ \forall\ x\in\mathbf{F}_p.\]

In particular, writing $N$ for the order of $q$ mod $p$ as an element of $\mathbf{F}_p^{\times}$, the group generated by $\{l_q, t\}$ has order $pN$ and it is a quotient of $\pi_1(\eta, \overline{\eta})$ tamely ramified over $S$ :
\[\pi_1^t(\eta, \overline{\eta})\to \langle l_q, t\rangle, \ F\mapsto l_q,\ T\mapsto t.\]

Suppose that another representation of $\pi_1(\eta, \overline{\eta})$ in the group of affine linear transformations of $\mathbf{F}_p$ is ramified over $S$ and has $t$ in its image $\mathfrak{G}$. Let $I$ (resp. $P$) be the image in $\mathfrak{G}$ of the inertia (resp. wild inertia) subgroup of $\pi_1(\eta, \overline{\eta})$.
\smallskip

In the group of affine linear transformations of
$\mathbf{F}_p$, the subgroup of translations is its own centralizer and it intersects $P$ in $1$. So $P=1$.
So $I=I/P$ is cyclic and $\neq 1$. Either $I$ contains $t$ or it intersects the group of translations in $1$. In both cases, $t$ commutes with all elements of $I$. Hence $I$ is the group of all translations.
\smallskip

In brief, the quotient homomorphism $\pi_1(\eta, \overline{\eta})\to\mathfrak{G}$ factors through $\pi_1^t(\eta, \overline{\eta})=\langle F, T\rangle$ and it maps $T$ to a non-zero translation.
\smallskip

Let the image of $F$ (resp. $T$) in $\mathfrak{G}$ be $l_at^b$ (resp. $t^c$), where $a, c\in\mathbf{F}_p^{\times}, b\in\mathbf{F}_p$. The identity
\[(l_at^b)t^c(l_at^b)^{-1}=(t^c)^q\] says that $ac=qc$, namely, that $a=q$ mod $p$. So
$\mathfrak{G}=\langle l_qt^b, t^c\rangle=\langle l_q, t\rangle$.

\end{proof}
\smallskip

\smallskip

\smallskip


\section{Type $A$}
\smallskip

\smallskip

Let $(S,\eta,s)$, $\mathrm{char}(s)=\ell$, be as in $\S 4$.
\smallskip

\smallskip

{\bf Proposition 5.1.} \emph{For every integer $n\geq 1$,
$(A_n,\alpha_1)$ and $(A_n,\alpha_n)$ are elliptic over $\eta$.}

\begin{proof} The subgroup of $\mathfrak{S}_{n+1}$
generated by the cycle $(12\cdots n+1)$ acts transitively on $\{1,\cdots, n+1\}$ and
permutes transitively the subsets of $\{1,\cdots, n+1\}$ of
cardinality $n$. As $\langle(12\cdots
n+1)\rangle=\mathbf{Z}/(n+1)\mathbf{Z}$ is a quotient of
$\pi_1(\eta,\overline{\eta})$ (\S 4), both $(A_n,\alpha_1)$ and $(A_n, \alpha_n)$
are elliptic over $\eta$ (3.1), 1).

\end{proof}
\smallskip

{\bf Lemma 5.2.} \emph{Let $X$ be a finite set of cardinality $q\geq
4$. Let $r$ be an integer such that $2\leq r\leq q/2$. Suppose that the subsets of $X$ of cardinality $r$ are permuted transitively by a solvable subgroup $\mathfrak{G}$ of $\mathrm{Aut}(X)$. Then $r<4$. Moreover,}
\smallskip

1) \emph{If $r=2$, $\mathfrak{G}$ acts $2$-transitively on $X$ unless }:
\smallskip

--- \emph{$X=\mathbf{F}_q$, $q\equiv 3$ mod $4$ and, for some subfield $k$ of $\mathbf{F}_q$, $\mathfrak{G}$ consists of all transformations of the form }:
\[ x\mapsto a^2\varphi(x)+b,\ \forall\ x\in\mathbf{F}_q\] \emph{where $a\in\mathbf{F}_q^{\times}$, $b\in\mathbf{F}_q$, $\varphi\in\mathrm{Gal}(\mathbf{F}_q/k)$.}
\smallskip

2) \emph{If $r=3$, then $X=\mathbf{F}_{32}$ or $\mathbf{F}_8$. When
$X=\mathbf{F}_{32}$, $\mathfrak{G}$ consists of all affine
semi-linear transformations of $X$. When $X=\mathbf{F}_8$,
$\mathfrak{G}$ consists of either all affine semi-linear
transformations or only of the affine linear transformations of $X$.}

\begin{proof} That $r<4$ as well as 2) is extracted from
\cite{livingstone}, p. 402--403.
\smallskip

Suppose that $r=2$ and that $\mathfrak{G}$ does not act $2$-transitively on $X$. By
\emph{loc.cit.}, then $X=\mathbf{F}_{p^d}$, $p$ prime $\equiv 3$ mod $4$,
$d$ is odd and $\mathfrak{G}=\mathfrak{L}\mathfrak{T}$, where
$\mathfrak{L}\leq\mathrm{GL}_d(\mathbf{F}_p)$ has odd order and where $\mathfrak{T}$ is the group of all translations of $X$. Observe that
$-1$ then normalizes $\mathfrak{G}$ and that
$\{1,-1\}\mathfrak{G}$ acts $2$-transitively on $X$, where $-1: x\mapsto -x$, $\forall\ x\in X$. Now 1) follows by the classification of $2$-transitive solvable permutation
groups.

\end{proof}
\smallskip

{\bf Corollary 5.3.} \emph{If $4\leq r\leq (n+1)/2$, then
$(A_n,\alpha_r)$ and $(A_n,\alpha_{n+1-r})$ are not elliptic over
$\eta$. The pairs $(A_n,\alpha_3)$ and $(A_n,\alpha_{n-2})$ are
elliptic over $\eta$ only if $n=7$ or $31$. The pairs
$(A_n,\alpha_2)$ and $(A_n,\alpha_{n-1})$ are elliptic over $\eta$ only
if $n=p^d-1$, $p$ prime, $d\geq 1$.}

\begin{proof} This is immediate from (5.2)+(3.1), 1). Recall that the group $\pi_1(\eta, \overline{\eta})$ is pro-solvable (\S 4).

\end{proof}
\smallskip

{\bf Proposition 5.4.} \emph{Let $p$ be a prime number and $d$ an integer $\geq 1$. The pairs $(A_{p^d-1},\alpha_2)$ and $(A_{p^d-1},\alpha_{p^d-2})$ are elliptic
over $\eta$ if $p=\ell$ and only if $p=\ell$ when $d\geq 2$.}

\begin{proof} If a solvable subgroup of $\mathfrak{S}_{p^d}$
permutes transitively the $2$-point subsets of $\mathbf{F}_p^d=V$, then it is of the form $\mathfrak{G}=\mathfrak{L}\mathfrak{T}$, where $\mathfrak{L}$ is a certain subgroup of $\mathrm{GL}(V)$ acting irreducibly on $V$ and where $\mathfrak{T}$ is the group of all translations of $V$ (5.2), 1).
\smallskip

If $p\neq\ell$ and if $d\geq 2$, $\pi_1(\eta,\overline{\eta})$ has no such quotient as $\mathfrak{G}$ (4.2), 1) and hence $(A_{p^d-1},\alpha_2)$ and $(A_{p^d-1},\alpha_{p^d-2})$ are not elliptic over $\eta$ (3.1), 1).
\smallskip

Suppose that $p=\ell$. On $\mathbf{F}_{\ell^d}$ the group $\mathfrak{G}$ of all affine linear transformations acts $2$-transitively. And by (4.1) $\mathfrak{G}$ is a
quotient of $\pi_1(\eta,\overline{\eta})$. So
$(A_{\ell^d-1},\alpha_2)$ and $(A_{\ell^d-1},\alpha_{\ell^d-2})$ are elliptic over $\eta$ (3.1), 1).

\end{proof}
\smallskip

{\bf Proposition 5.5.} \emph{Let $p$ be an odd prime different from
$\ell$.}
\smallskip

--- \emph{Case $p\equiv 1$ mod $4$} : \emph{Then $(A_{p-1},\alpha_2)$ and
$(A_{p-1},\alpha_{p-2})$ are elliptic over $\eta$ if and only if
$\mathrm{Card}(k(s))$ mod $p$ generates $\mathbf{F}_p^{\times}$.}
\smallskip

--- \emph{Case $p\equiv 3$ mod $4$} : \emph{Then $(A_{p-1},\alpha_2)$ and $(A_{p-1},\alpha_{p-2})$ are elliptic over $\eta$ if and only if $\mathrm{Card}(k(s))$ mod $p$
generates a subgroup of $\mathbf{F}_p^{\times}$ of index $\leq 2$.}

\begin{proof} By (3.1), 1) the pairs $(A_{p-1},\alpha_2)$ and $(A_{p-1},\alpha_{p-2})$ are
elliptic over $\eta$ if and only if there is a representation
$\pi_1(\eta,\overline{\eta})\to\mathfrak{S}_p$ whose image $\mathfrak{G}$ permutes transitively
the $2$-point subsets of $\mathbf{F}_p$.
\smallskip

By (5.2), 1) and by the classification of $2$-transitive solvable
permutation groups of degree $p$, such $\mathfrak{G}$ can only be
\smallskip

--- (Case $p\equiv 1$ mod $4$) the group of all affine linear transformations of $\mathbf{F}_p$.
\smallskip

--- (Case $p\equiv 3$ mod $4$) either the group of all affine linear transformations of $\mathbf{F}_p$ or the subgroup consisting of all transformations of the form $x\mapsto a^2x+b$, $\forall\ x\in\mathbf{F}_p$, where $a\in\mathbf{F}_p^{\times}$, $b\in\mathbf{F}_p$.

\smallskip

Now by (4.2), 2) the lemma follows.

\end{proof}
\smallskip

{\bf Proposition 5.6.} \emph{The pairs $(A_7,\alpha_3)$ and $(A_7,\alpha_5)$ are elliptic over $\eta$ if and only if $k(s)$ is of characteristic $2$.}

\begin{proof} In (5.2), 2) either of the two solvable subgroups of $\mathfrak{S}_8$ that permute transitively the $3$-point subsets of
$\mathbf{F}_8$ contains all translations of $\mathbf{F}_8$.
So $(A_7,\alpha_3)$ and $(A_7,\alpha_5)$ are elliptic over $\eta$ only if $k(s)$ is of characteristic $\ell=2$ (3.1), 1)+(4.2), 1).
\smallskip

If $\ell=2$, the group of all affine linear transformations of
$\mathbf{F}_8$ is a quotient of $\pi_1(\eta,\overline{\eta})$
(4.1) and hence $(A_7,\alpha_3)$ and $(A_7,\alpha_5)$ are elliptic over $\eta$ (3.1), 1)+(5.2), 2).

\end{proof}
\smallskip

{\bf Proposition 5.7.} \emph{The pairs $(A_{31},\alpha_3)$ and $(A_{31},\alpha_{29})$ are elliptic over $\eta$ if and only if $\ell=2$, $5\nmid[s:\mathbf{F}_2]$.}

\begin{proof} The pairs $(A_{31},\alpha_3)$ and $(A_{31},\alpha_{29})$
are elliptic over $\eta$ if and only if $\pi_1(\eta, \overline{\eta})$ has as quotient the group $\mathfrak{G}$ of
all affine semi-linear transformations of $\mathbf{F}_{32}$ (3.1), 1), (5.2), 2).
\smallskip

By (4.2), 1) $\mathfrak{G}$ is a quotient of
$\pi_1(\eta,\overline{\eta})$ only if $k(s)$ is of characteristic $\ell=2$.
\smallskip

Suppose that $\ell=2$.
\smallskip

\emph{Suppose that $\pi_1(\eta, \overline{\eta})$ has $\mathfrak{G}$ as a quotient. Then $5\nmid [s:\mathbf{F}_2]$.}

\smallskip

Let $I$ (resp. $P$) be the image in $\mathfrak{G}$ of the inertia (resp. wild inertia) subgroup of $\pi_1(\eta, \overline{\eta})$. It is immediate that $P$ (resp. $I$) must consist of all translations (resp. all affine linear transformations) of $\mathbf{F}_{32}$. The subgroup of $\mathfrak{G}$ generated by the Frobenius $F: x\mapsto x^2$ and the scalar multiplications $l_a: x\mapsto ax$ is isomorphic to $\mathfrak{G}/P$. By (4.2), 2) one concludes that the element
$\mathrm{Card}(k(s))$ mod $31$ must be of order $5$ in $\mathbf{F}_{31}^{\times}$. That is, $5\nmid [s:\mathbf{F}_2]$, since
$2$ mod $31$ is of order $5$ in $\mathbf{F}_{31}^{\times}$.

\smallskip

\emph{Suppose that $5\nmid [s:\mathbf{F}_2]$. Then $\mathfrak{G}$ is a quotient of $\pi_1(\eta, \overline{\eta})$.}

\smallskip

Let $S'$ be the spectra of a discrete valuation ring such that $S'$ is finite \'{e}tale Galois over $S$ with cyclic Galois group of order $5$ (\S 4). Let $\eta'$ (resp. $s'$) be the generic (resp. closed) point of $S'$, $\zeta\in\mathrm{Gal}(S'/S)$ a generator, $\pi\in\Gamma(S, \mathcal{O}_S)$ a uniformizer and let $u'\in\Gamma(S', \mathcal{O}_{S'})^{\times}$ be a unit such that the images of $u', \zeta(u'),\cdots, \zeta^4(u')$ in $k(s')$ form a normal base over $k(s)$. Then
\[\eta'[z,x_1,\cdots,x_5]/(z^{31}-\pi, x_1^2-1-z\zeta(u'), \cdots, x_5^2-1-z\zeta^5(u'))\] is connected and Galois over $\eta$ with Galois group $\mathfrak{G}$.

\end{proof}

\smallskip

\smallskip

\smallskip


\section{Type ${}^2A$}
\smallskip

\smallskip

\smallskip

{{\bf Proposition 6.1.} \emph{Let $X$ be a finite
set of even cardinality $2d$. Let $\mathfrak{G}$ be a solvable subgroup of $\mathrm{Aut}(X)$ which permutes the subsets of $X$ of cardinality $d$ in $2$ orbits.}
\smallskip

\emph{The following list enumerates such $(X,\mathfrak{G})$ up to equivalence }:
\smallskip

1) \emph{$X=\{o,1\}$, $\mathfrak{G}=1$.}
\smallskip

2) \emph{$X=\{o,1,2,3\}$, $\mathfrak{G}$ fixes $o$ and on
$\{1,2,3\}$ it is either $\mathfrak{S}_3$ or $\mathfrak{A}_3$.}
\smallskip

3) \emph{$X=\{o\}\cup \mathbf{F}_5$, $\mathfrak{G}$ fixes $o$ and on $\mathbf{F}_5$ it is
the group of all affine linear transformations.}
\smallskip

4) \emph{$X=\mathbf{Z}/4\mathbf{Z}$, $\mathfrak{G}$ consists of either
all transformations
\[
x\mapsto ax+b,\ \forall\ x\in\mathbf{Z}/4\mathbf{Z}
\] where
$a\in(\mathbf{Z}/4\mathbf{Z})^{\times}$, $b\in\mathbf{Z}/4\mathbf{Z}$ or only of the translations
\[x\mapsto x+b,\ \forall\ x\in\mathbf{Z}/4\mathbf{Z}\] where $b\in\mathbf{Z}/4\mathbf{Z}$.}
\smallskip

5) \emph{$X=\{1,\cdots,6\}$, either $\mathfrak{G}$ is the normalizer $\mathfrak{N}$
in $\mathrm{Aut}(X)$ of a partition $X=\{a,b,c\}\cup\{a',b',c'\}$ or it is the subgroup of $\mathfrak{N}$ generated by $\mathfrak{Alt}(\{a,b,c\})
\mathfrak{Alt}(\{a',b',c'\})$ and one of the following subgroups }:
\smallskip

--- \emph{$\langle(aa')(bb')(cc')\rangle$}
\smallskip

--- \emph{$\langle (aa'bb')(cc')\rangle$}
\smallskip

--- \emph{$\langle(aa')(bb')(cc'), (ab)(a'b')\rangle$}
\smallskip

6) \emph{$X=\{1,\cdots,6\}$, either $\mathfrak{G}$ is the normalizer
$\mathfrak{N}$ in $\mathrm{Aut}(X)$ of a partition
$X=\{a,a'\}\cup\{b,b'\}\cup\{c,c'\}$ or it is the subgroup of
$\mathfrak{N}$ generated by $\{(aa'), (bb'), (cc'),
(abc)(a'b'c')\}$.}
\smallskip

7) \emph{$X=\mathbf{F}_8$, $\mathfrak{G}$ consists of either all
affine semi-linear transformations
\[
x\mapsto ax^{2^c}+b, \ \forall\ x\in\mathbf{F}_8
\] where
$a\in\mathbf{F}_8^{\times}$, $b\in\mathbf{F}_8$,
$c\in\mathbf{Z}/3\mathbf{Z}$ or only of the affine linear
transformations
\[
x\mapsto ax+b, \ \forall\ x\in\mathbf{F}_8
\] where
$a\in\mathbf{F}_8^{\times}$, $b\in\mathbf{F}_8$.}
\smallskip

\smallskip

\smallskip

The proof is divided into several parts : (6.2), (6.4), (6.5), (6.6).
\smallskip

\smallskip

{\bf Lemma 6.2.} \emph{With the notations and assumptions of $(6.1)$, suppose furthermore that $\mathfrak{G}$ does not act transitively on $X$.}
\smallskip

\emph{The following list enumerates all such $(X, \mathfrak{G})$ up to equivalence }:
\smallskip

1) \emph{$X=\{o, 1\}$, $\mathfrak{G}=1$.}
\smallskip

2) \emph{$X=\{o, 1, 2, 3\}$, $\mathfrak{G}$ fixes $o$ and on $\{1, 2, 3\}$ it is either $\mathfrak{S}_3$ or $\mathfrak{A}_3$.}
\smallskip

3) \emph{$X=\{o\}\cup\mathbf{F}_5$, $\mathfrak{G}$ fixes $o$ and on $\mathbf{F}_5$ it is the group of affine linear transformations.}

\begin{proof} Choose $o\in X$ such that $O=\mathfrak{G}o$ has cardinality $\leq d=\mathrm{Card}(X)/2$. Such a point exists since by assumption $\mathfrak{G}$ does not act transitively on $X$.
\smallskip

Choose a subset $Y$ (resp. $Z$) of $X$ with $d$ elements such
that $Y$ (resp. $Z$) contains (resp. is disjoint with) $O$.
One has $gY\supset O$ and $gZ\cap O=\emptyset$, $\forall\ g\in\mathfrak{G}$. So $\mathfrak{G}Y$ and $\mathfrak{G}Z$ are these two $\mathfrak{G}$-orbits in the collection of $d$-point subsets of $X$.
\smallskip

Choose a point $z\in Z$. The set
$\{o\}\cup Z\backslash \{z\}$
has $d$ elements and it intersects $O$ in $\{o\}$. So $O=\{o\}$.
\smallskip

Now $X\backslash \{o\}$ has $2d-1$ elements and its subsets of
cardinality $d$ form a single $\mathfrak{G}$-orbit $\mathfrak{G}Z$. The
following lemma applies.

\end{proof}

\smallskip

{\bf Lemma 6.3.} \emph{Let $X$ be a finite set
of odd cardinality $2d-1$. Let $\mathfrak{G}$ be a solvable subgroup of $\mathrm{Aut}(X)$ which permutes transitively the subsets of $X$ of cardinality $d$.}
\smallskip

\emph{The following list enumerates such $(X, \mathfrak{G})$ up to equivalence }:
\smallskip

1) \emph{$X=1$, $\mathfrak{G}=1$.}
\smallskip

2) \emph{$X=\{1,2,3\}$, $\mathfrak{G}=\mathfrak{S}_3$ or
$\mathfrak{A}_3$.}
\smallskip

3) \emph{$X=\mathbf{F}_5$, $\mathfrak{G}$ consists of all affine
linear transformations
\[x\mapsto ax+b,\ \forall\ x\in\mathbf{F}_5\] where
$a\in\mathbf{F}_5^{\times}$, $b\in\mathbf{F}_5$.}

\begin{proof} If $d=1$, then $X=1$ and $\mathfrak{G}=1$, hence 1).
\smallskip

Suppose $d>1$. Notice that
\smallskip

\smallskip

--- \emph{The group $\mathfrak{G}$ acts transitively on $X$ }:
\smallskip

\smallskip

Otherwise, some
$\mathfrak{G}$-orbit in $X$, say $O$, has $<d$ elements. Choose
a subset $Y$ of $X$ with $d$ elements so that $Y$ contains $O$. For all $g\in\mathfrak{G}$, $O\subset gY$. Namely, $O$ is contained in every subset of $X$ of
cardinality $d$. The complement of $O$ in $X$ has $>(2d-1)-d=d-1$ elements. Hence $X\backslash O$ contains at least one set of cardinality $d$. A contradiction.
\smallskip

Fix a point $o\in X$. Then
\smallskip

\smallskip

--- \emph{The stabilizer $\mathfrak{G}_o$ of $o$ in $\mathfrak{G}$ is a maximal subgroup of $\mathfrak{G}$} :
\smallskip

\smallskip

Assume $\mathfrak{G}_o<\mathfrak{H}<\mathfrak{G}$ for a group
$\mathfrak{H}$. Then
$1<(\mathfrak{G}:\mathfrak{H}),(\mathfrak{H}:\mathfrak{G}_o)<d$,
because
\[
(\mathfrak{G}:\mathfrak{H})(\mathfrak{H}:\mathfrak{G}_o)=(\mathfrak{G}:\mathfrak{G}_o)=\mathrm{Card}(\mathfrak{G}.o)=\mathrm{Card}(X)=2d-1.
\]
As $\mathfrak{H}.o\simeq \mathfrak{H}/\mathfrak{G}_o$, $X\backslash
(\mathfrak{H}.o)$ has cardinality $>(2d-1)-d=d-1$. Pick a set $Y\subset
X\backslash (\mathfrak{H}.o)$ with $d$ elements. Then $gY\cap
g\mathfrak{H}.o=\emptyset$, $\forall\ g\in \mathfrak{G}$. Therefore, each
subset of $X$ of cardinality $d$ is disjoint with at least one translate
$g\mathfrak{H}.o$ of $\mathfrak{H}.o$. Let $\mathfrak{R}$ be a set of representatives
for $\mathfrak{G}/\mathfrak{H}$, which has cardinality
$(\mathfrak{G}:\mathfrak{H})<d$. So
$\mathfrak{R}.o$ is contained in some set of cardinality $d$ in $X$, say $Z$. But $Z$ intersects every $g\mathfrak{H}.o$, $\forall\ g\in \mathfrak{G}$. A contradiction.
\smallskip

\smallskip

--- \emph{The group $\mathfrak{G}_o$ contains no normal subgroups of
$\mathfrak{G}$ other than $1$} :
\smallskip

\smallskip

Let $\mathfrak{N}$ be a subgroup of $\mathfrak{G}_o$ such that $\mathfrak{N}$ is normal in $\mathfrak{G}$. Then
$\mathfrak{N}g.o=g\mathfrak{N}.o=g.o$, $\forall\ g\in \mathfrak{G}$. That is, $\mathfrak{N}$ fixes every element of $\mathfrak{G}.o=X$. So $\mathfrak{N}=1$.
\smallskip

\smallskip

Let $\mathfrak{U}$ be the last term $>1$ in the derived series of
$\mathfrak{G}$. Then
$[\mathfrak{U},\mathfrak{U}]=1$, as $\mathfrak{G}$ is solvable. So $\mathfrak{U}$ is abelian on which $\mathfrak{G}$ acts by conjugation. Let $V\subset \mathfrak{U}$ be a simple
sub-$\mathfrak{G}$-module\ ; it is an $\mathbf{F}_p$-vector
space for some prime number $p$. Let $f=\mathrm{dim}\ V$.
\smallskip

\smallskip

--- \emph{One has $V\mathfrak{G}_o=\mathfrak{G}$ and $V\cap\mathfrak{G}_o=1$ }:
\smallskip

\smallskip

The maximal subgroup $\mathfrak{G}_o$ does not contain $V$, as $V$ is normal in $\mathfrak{G}$. So $V\mathfrak{G}_o$ contains $\mathfrak{G}_o$
properly and so $V\mathfrak{G}_o=\mathfrak{G}$. The intersection $V\cap \mathfrak{G}_o$ is
normalized by $\mathfrak{G}_o$ and by $V$, $V$ being abelian, and thus
by $V\mathfrak{G}_o=\mathfrak{G}$. So $V\cap \mathfrak{G}_o$ is
a sub-$\mathfrak{G}$-module of $V$ distinct from $V$. So $V\cap\mathfrak{G}_o=1$.
\smallskip

\smallskip

--- \emph{The map $V\to X$, $v\mapsto v.o$, is a bijection} :
\smallskip

\smallskip

It is surjective
because $X=\mathfrak{G}.o=V\mathfrak{G}_o.o=V.o$. It is injective
because if $v.o=v'.o$, then $v^{-1}v'\in V\cap \mathfrak{G}_o=1$ and $v=v'$.
\smallskip

Now $p^f=\mathrm{Card}(V)=\mathrm{Card}(X)=2d-1$. So $p>2$.
\smallskip

\smallskip

--- \emph{The representation $\mathfrak{G}_o\to \mathrm{GL}(V)$,
$g\mapsto \mathrm{Int}(g)$, is faithful }:
\smallskip

\smallskip

Let $g\in \mathfrak{G}_o$ be such that
$\mathrm{Int}(g)=1$ on $V$. Then
$gv.o=gvg^{-1}.o=\mathrm{Int}(g)(v).o=v.o$, $\forall\ v\in V$. So
$g$ fixes each point of $V.o=X$.
\smallskip

\smallskip

Pick a prime $p'$ such that $d<p'<2d$ (Bertrand's postulate).
\smallskip

\smallskip

--- \emph{Then $p'=p$ }:
\smallskip

\smallskip

Suppose $p'\neq p$. By its choice, $p'$ divides
$N:={2d-1\choose d}$. Notice that $X$ has $N$ subsets of cardinality $d$. These $N$ subsets are permuted transitively by $\mathfrak{G}$. And $\mathfrak{G}=\mathfrak{G}_oV$ imbeds into $\mathrm{GL}(V)V$ by the faithful representation $\mathrm{Int}:\mathfrak{G}_o\hookrightarrow\mathrm{GL}(V)$. So $p'$ divides the order of $\mathrm{GL}(V)V$. So $p'$ divides $p^i-1$ for some $i\in\{1,\cdots, f\}$, as $p'\neq p$. But this is absurd. For, $p'$ is odd, $p^i-1$ is even and $p^i-1\leq p^f-1=2d-2<2p'-2$.
\smallskip

\smallskip

--- \emph{Then $f=1$} :
\smallskip

\smallskip

For, $p^f=2d-1<2p'-1=2p-1$.
\smallskip

\smallskip

--- \emph{One has $d\leq 3$} :
\smallskip

\smallskip

This is immediate from the division :
\[{2d-1\choose d}\ |\ \mathrm{Card}(\mathrm{GL}(V)V)=p(p-1)=(2d-1)(2d-2).\]

--- \emph{Case $d=2$. Then $\mathfrak{G}=\mathfrak{S}_3$ or $\mathfrak{A}_3$ on $X=\{1, 2, 3\}$ }:
\smallskip

\smallskip

The set $X$ has $2d-1=3$ elements. The transitivity of the $\mathfrak{G}$-action on the $2$-point subsets of $X$ is equivalent to the transitivity of the $\mathfrak{G}$-action on $X$. So $\mathfrak{G}$ is either $\mathfrak{S}_3$ or $\mathfrak{A}_3$.
\smallskip

\smallskip

--- \emph{Case $d=3$. Then $\mathfrak{G}$ consists of all affine linear transformations of $\mathbf{F}_5=X$} :

\smallskip

\smallskip

The set $X$ as well as $V$ has $2d-1=5$ elements. So $V=\mathbf{F}_5$. And $\mathrm{GL}(V)V$ is the group of all affine linear transformations of $\mathbf{F}_5$ which acts $2$-transitively on $\mathbf{F}_5$. Indeed, if $a, b$ are two distinct points of $\mathbf{F}_5$, the affine linear transformation $x\mapsto (a-b)x+b$ maps $0$ to $b$ and maps $1$ to $a$. In particular, $\mathrm{GL}(V)V$ permutes transitively the $2$-point subsets, or what amounts to the same, the $3$-point subsets, of $\mathbf{F}_5$.
\smallskip

The unique index $2$ subgroup $H$ of $\mathrm{GL}(V)V$ consists of all transformations of the form :
\[x\mapsto a^2x+b,\ \forall\ x\in\mathbf{F}_5\] where $a\in\mathbf{F}_5^{\times}$, $b\in\mathbf{F}_5$. The $2$-point subsets $\{u, v\}$ of $\mathbf{F}_5$ are divided into $2$ $H$-orbits according to whether or not $u-v$ is a square in $\mathbf{F}_5^{\times}$. Notice that $-1$ is a square in $\mathbf{F}_5^{\times}$. The assertion evidently follows.

\end{proof}

\smallskip

{\bf Lemma 6.4.} \emph{With the notations and assumptions of $(6.1)$, let $o\in X$ be a point and $\mathfrak{G}_o$ its stabilizer in $\mathfrak{G}$. Suppose furthermore that $\mathfrak{G}$ acts transitively on $X$ and that the following condition holds }:
\smallskip

--- \emph{There is a subgroup $\mathfrak{H}$ of even index in $\mathfrak{G}$ such that $\mathfrak{H}$ contains $\mathfrak{G}_o$ properly.}
\smallskip

\emph{Then the following list enumerates such $(X, \mathfrak{G})$ up to equivalence }:

\smallskip

1) \emph{$X=\mathbf{Z}/4\mathbf{Z}$, $\mathfrak{G}$ consists of either
all transformations
\[
x\mapsto ax+b,\ \forall\ x\in\mathbf{Z}/4\mathbf{Z}
\] where
$a\in(\mathbf{Z}/4\mathbf{Z})^{\times}$, $b\in\mathbf{Z}/4\mathbf{Z}$ or only of the translations
\[x\mapsto x+b,\ \forall\ x\in\mathbf{Z}/4\mathbf{Z}\] where $b\in\mathbf{Z}/4\mathbf{Z}$.}
\smallskip

2) \emph{$X=\{1,\cdots,6\}$, $\mathfrak{G}$ is either the normalizer $\mathfrak{N}$
in $\mathrm{Aut}(X)$ of a partition $X=\{a,b,c\}\cup\{a',b',c'\}$ or it is the subgroup of $\mathfrak{N}$ generated by $\mathfrak{Alt}(\{a,b,c\})
\mathfrak{Alt}(\{a',b',c'\})$ and one of the following subgroups :}
\smallskip

--- \emph{$\langle(aa')(bb')(cc')\rangle$}
\smallskip

--- \emph{$\langle(aa'bb')(cc')\rangle$}
\smallskip

--- \emph{$\langle(aa')(bb')(cc'), (ab)(a'b')\rangle$}

\begin{proof} Let $(\mathfrak{G}:\mathfrak{H})=2r$ and let $\mathfrak{R}=\{g_1,\cdots, g_{2r}\}$ be a set of representatives for $\mathfrak{G}/\mathfrak{H}$. Notice that
\[
d=\frac{\mathrm{Card}(X)}{2}=\frac{(\mathfrak{G}:\mathfrak{H})}{2}(\mathfrak{H}:\mathfrak{G}_o)=r\
\mathrm{Card}(\mathfrak{H}.o)\geq 2r.
\]

In particular, if $\mathfrak{I}$ is a subset of $\mathfrak{R}$ of cardinality $r$, then
\[Z=\mathfrak{I}\mathfrak{H}.o\] has $d$ elements.
\smallskip

As $\mathrm{Card}(\mathfrak{R}.o)\leq
\mathrm{Card}(\mathfrak{R})=2r\leq d$, there is some set $Z'$ in $X$ with $d$ elements which contains $\mathfrak{R}.o$. By its choice $Z'$ intersects every $g\mathfrak{H}.o$,
$\forall\ g\in\mathfrak{G}$.
\smallskip

Since $\mathfrak{G}$ permutes the $d$-point subsets of $X$ in $2$ orbits, each of these sets satisfies one or the other of the following conditions :
\smallskip

\smallskip

i) \emph{It is equal to $\mathfrak{I}\mathfrak{H}.o$ for a subset $\mathfrak{I}$ of $\mathfrak{R}$, where $\mathfrak{I}$ has $r$ elements.}
\smallskip

ii) \emph{It intersects every translate $g\mathfrak{H}.o$, $\forall\ g\in\mathfrak{G}$.}

\smallskip

\smallskip

--- \emph{Then $r=1$} :

\smallskip

\smallskip

Assume $r>1$. Then the set
\[E:=\{g_1,\cdots,g_r\}\mathfrak{H}.o\cup \{g_{r+1}.o\}\backslash
\{g_1.o\}\] has $d$ elements and is disjoint with
$g_{2r}\mathfrak{H}.o$. But $E$ is not of the form $\mathfrak{I}\mathfrak{H}.o$
for any subset $\mathfrak{I}$ of $\mathfrak{R}$.
\smallskip

So $\mathfrak{R}=\{g_1,g_2\}$, $\mathrm{Card}(\mathfrak{H}.o)=d$,
$X=\mathfrak{H}.o\cup \tau\mathfrak{H}.o$, where $\tau:=g_1^{-1}g_2$, and
the $d$-point subsets of $X$ distinct from
$\mathfrak{H}.o$ and $\tau\mathfrak{H}.o$ are permuted transitively by
$\mathfrak{G}$.
\smallskip

\smallskip

--- \emph{Then $d\leq 3$ }:
\smallskip

\smallskip

Suppose $d>3$. Choose a point $o'\in \mathfrak{H}.o\backslash
\{o\}$. Both sets
\[
Y=\{o\}\cup \tau\mathfrak{H}.o\backslash\{\tau.o\} \ ,\
Y'=\{o,o'\}\cup\tau\mathfrak{H}.o\backslash\{\tau.o,\tau.o'\}
\]
have $d$ elements. Both are distinct from $\mathfrak{H}.o$ and
$\tau\mathfrak{H}.o$. But $Y\neq gY'$, $\forall\ g\in\mathfrak{G}$.
For, $Y\cap \mathfrak{H}.o$ consists of $1$ element, while $gY'\cap
\mathfrak{H}.o=g(Y'\cap g^{-1}\mathfrak{H}.o)$ consists of either $2$ or
$d-2$ elements, $\forall\ g\in\mathfrak{G}$.
\smallskip

\smallskip

--- \emph{Case $d=2$ }:
\smallskip

\smallskip

The set $X$ has $2d=4$ elements. Both $\mathfrak{H}.o$ and $\tau\mathfrak{H}.o$ have $2$ elements. As
$\mathfrak{H}$ is a subgroup of
$\mathrm{Aut}(\mathfrak{H}.o)\times\mathrm{Aut}(\tau\mathfrak{H}.o)$, it has $2$ or $4$ elements.
\smallskip

Suppose first that $\mathfrak{H}$ has $2$ elements. Then $\mathfrak{G}_o=1$, $|\mathfrak{G}|=4$ and $\mathfrak{G}$ acts simply transitively on $X$.
\smallskip

Notice that the translation action on itself of $\mathbf{Z}/4\mathbf{Z}$ permutes the $2$-point subsets $\{u, v\}$ of $\mathbf{Z}/4\mathbf{Z}$ in $2$ orbits according to whether or not $u-v$ belongs to the subgroup $\mathfrak{H}=2\mathbf{Z}/4\mathbf{Z}$. And, the translation action on itself of $\mathbf{Z}/2\mathbf{Z}\times\mathbf{Z}/2\mathbf{Z}$ permutes its $2$-point subsets $\{u, v\}$ in $3$ orbits according to which subgroup $u-v$ generates.
\smallskip

Suppose next that $\mathfrak{H}$ has $4$ elements. Then $\mathfrak{G}$ is a $2$-Sylow subgroup of $\mathrm{Aut}(X)$. It is isomorphic to the group of all transformations :
\[x\mapsto ax+b,\ \forall\ x\in\mathbf{Z}/4\mathbf{Z}\] where $a\in(\mathbf{Z}/4\mathbf{Z})^{\times}$, $b\in\mathbf{Z}/4\mathbf{Z}$. The group $\mathfrak{G}$ permutes the $2$-point subsets $\{u, v\}$ of $\mathbf{Z}/4\mathbf{Z}$ in $2$ orbits according to whether or not $u-v$ lies in $2\mathbf{Z}/4\mathbf{Z}$.

\smallskip

\smallskip

--- \emph{Case $d=3$ }:

\smallskip

\smallskip

The set $X$ has $2d=6$ elements. Both $\mathfrak{H}.o$ and $\tau\mathfrak{H}.o$ have $3$ elements. Let $\mathfrak{N}$ denote the normalizer in $\mathrm{Aut}(X)$ of the partition
\[X=\mathfrak{H}.o\cup \tau\mathfrak{H}.o.\] The group $\mathfrak{N}$ has $72$ elements. Besides $\mathfrak{H}.o$ and $\tau\mathfrak{H}.o$, there are $18$ subsets in $X$ of cardinality $3$. These $18$ sets are permuted transitively by $\mathfrak{G}$. So $\mathfrak{G}$ is of index $1$, $2$ or $4$ in $\mathfrak{N}$.

\smallskip

We write $\mathfrak{H}.o=\{1, 2, 3\}$ and
$\tau\mathfrak{H}.o=\{4, 5, 6\}$.

\smallskip

Let $\mathfrak{P}:=\mathfrak{Alt}(\{1, 2, 3\})\times\mathfrak{Alt}(\{4, 5, 6\})$. It is the unique $3$-Sylow subgroup of $\mathfrak{N}$ and of $\mathfrak{G}$. Let $\mathfrak{Q}$ be a $2$-Sylow subgroup of $\mathfrak{G}$. Thus $\mathfrak{G}=\mathfrak{P}\mathfrak{Q}$ and $\mathfrak{Q}$ is of order $2$, $4$ or $8$.
\smallskip

\smallskip

i) \emph{Case $\mathrm{Card}(\mathfrak{Q})=2$ }:
\smallskip

\smallskip

Let $\mathfrak{Q}=\{1,\alpha\}$, where $\alpha$ transforms $\{1,2,3\}$ to $\{4,5,6\}$. If say
\[\alpha: 1\mapsto 4,\ 2\mapsto 5,\ 3\mapsto 6\] then $\alpha=(14)(25)(36)$.
\smallskip

\smallskip

ii) \emph{Case $\mathrm{Card}(\mathfrak{Q})=4$, $\mathfrak{Q}$ cyclic }:
\smallskip

\smallskip

Let $\alpha$ be a generator of $\mathfrak{Q}$ which then transforms $\{1, 2, 3\}$ to $\{4, 5, 6\}$. So $\alpha^2$ normalizes $\{1, 2, 3\}$ and, being of order $2$, $\alpha^2$ fixes a point, say $3$, in $\{1,2,3\}$. If say
\[\alpha: 1\mapsto 4,\ 2\mapsto 5,\ 3\mapsto 6\] then $\alpha=(1425)(36)$.

\smallskip

\smallskip

iii) \emph{Case $\mathrm{Card}(\mathfrak{Q})=4$, $\mathfrak{Q}$ non cyclic }:
\smallskip

\smallskip

Let $\mathfrak{Q}=\{1, \alpha, \beta, \gamma\}$. Suppose that $\alpha$ and $\beta$ (resp. $\gamma$) transform $\{1, 2, 3\}$ to $\{4, 5, 6\}$ (resp. normalizes $\{1, 2, 3\}$). Then $\gamma$ fixes a point, say $3$, in $\{1, 2, 3\}$. If say
\[\alpha: 1\mapsto 4,\ 2\mapsto 5,\ 3\mapsto 6\] then $\alpha=(14)(25)(36)$, $\gamma=(12)(45)$, $\beta=(15)(24)(36)$.

\smallskip

\smallskip

iv) \emph{Case $\mathrm{Card}(\mathfrak{Q})=8$ }:

\smallskip

\smallskip

Then $\mathfrak{G}=\mathfrak{N}$.
\smallskip

\smallskip

It remains to verify that in all these cases $\mathfrak{G}$ permutes transitively the $3$-point subsets $Y$ of $X$ other than $\{1, 2, 3\}$ and $\{4, 5, 6\}$.
\smallskip

Given such a subset $Y$ of $X$, notice that there is an element $p\in\mathfrak{P}$ such that $pY$ is either $\{1,2,6\}$ or $\{3, 4, 5\}$. Then, in the notations of i)--iii), $\alpha$ transforms $\{1,2,6\}$ to $\{3,4,5\}$.

\end{proof}
\smallskip

{\bf Lemma 6.5.} \emph{With the notations and assumptions of $(6.1)$, let $o\in X$ be a point and $\mathfrak{G}_o$ its stabilizer in $\mathfrak{G}$. Suppose furthermore that $\mathfrak{G}_o$ acts transitively on $X$ and that the following condition holds }:
\smallskip

--- \emph{There is a subgroup $\mathfrak{H}$ of odd index $>1$ in $\mathfrak{G}$ such that $\mathfrak{H}$ contains $\mathfrak{G}_o$ properly.}
\smallskip

\emph{Then $X=\{1,\cdots,6\}$,
$\mathfrak{G}$ is either the normalizer $\mathfrak{N}$ in
$\mathrm{Aut}(X)$ of a partition
$X=\{a,a'\}\cup\{b,b'\}\cup\{c,c'\}$ or it is the subgroup of
$\mathfrak{N}$ generated by $\{(aa'), (bb'), (cc'),
(abc)(a'b'c')\}$.}

\begin{proof} Let $r$ be an integer $\geq 1$ such that $(\mathfrak{G}:\mathfrak{H})=2r+1$. Let
$\mathfrak{R}=\{g_1,\cdots,g_{2r+1}\}$ be a set of
representatives for $\mathfrak{G}/\mathfrak{H}$.
\smallskip

The identity
\[
d=\frac{\mathrm{Card}(X)}{2}=\frac{(\mathfrak{G}:\mathfrak{H})}{2}(\mathfrak{H}:\mathfrak{G}_o)=(r+\frac{1}{2})\
\mathrm{Card}(\mathfrak{H}.o),
\]
implies in particular that $\mathfrak{H}.o$ has even, say $2f$, elements.
\smallskip

Choose a subset
$B\subset g_{r+1}\mathfrak{H}.o\backslash \{g_{r+1}.o\}$ of
cardinality $f$. Then
\[Y=\{g_1,\cdots,g_r\}\mathfrak{H}.o\cup B\] has
$d$ elements.
\smallskip

As $
\mathrm{Card}(\mathfrak{R})\leq d$, $\mathfrak{R}.o$ is contained in some set $Y'$
with $d$ elements. This $Y'$ intersects every $g\mathfrak{H}.o$,
$\forall\ g\in\mathfrak{G}$.
\smallskip

By assumption $\mathfrak{G}$ has $2$ orbits in the collection of $d$-point subsets of $X$. Each of these sets satisfies thus one or the other of the following two conditions :
\smallskip

--- \emph{It is equal to $\mathfrak{I}\mathfrak{H}.o\cup B'$ for some subset $\mathfrak{I}$ of $\mathfrak{R}$ of cardinality $r$ and some subset $B'$ of $z\mathfrak{H}.o$ of cardinality $f$, where $z$ is an element of $\mathfrak{R}\backslash\mathfrak{I}$.}
\smallskip

--- \emph{It intersects every $g\mathfrak{H}.o$, $\forall\ g\in\mathfrak{G}$.}

\smallskip

Notice that every such set $\mathfrak{I}\mathfrak{H}.o\cup B'$
intersects precisely $r+1$ members among
\[g_1\mathfrak{H}.o,\ \cdots,\ g_{2r+1}\mathfrak{H}.o.\]

--- \emph{Then $f=1$ }:
\smallskip

\smallskip

Assume $f>1$. Then the set
\[
E:=\{g_1,\cdots,g_{r-1}\}\mathfrak{H}.o\cup(g_r\mathfrak{H}.o\backslash
\{g_r.o\})\cup (B\cup\{g_{r+1}.o\})
\]
has $d$ elements and is disjoint with $g_{2r+1}\mathfrak{H}.o$. But $E$ is
not of the form $\mathfrak{I}\mathfrak{H}.o\cup B'$ for any
$\mathfrak{I}\subset\mathfrak{R}$ of cardinality $r$, any $z\in\mathfrak{R}\backslash\mathfrak{I}$ and any $B'\subset
z\mathfrak{H}.o$ of cardinality $f$.
\smallskip

\smallskip

It follows that $B$ consists of $f=1$ element and that $d=2r+1$.
\smallskip

\smallskip

--- \emph{Then $r=1$ }:
\smallskip

\smallskip

Suppose $r>1$. Then the set
\[
F:=\{g_1,\cdots,g_{r-1}\}\mathfrak{H}.o\cup
(g_r\mathfrak{H}.o\backslash\{g_r.o\})\cup
\{g_{r+1}.o\}\cup\{g_{2r+1}.o\}
\]
has $d$ elements and is disjoint with $g_{r+2}\mathfrak{H}.o$. But $F$
intersects $r+2$, rather than $r+1$, members among
\[g_1\mathfrak{H}.o,\ \cdots,\ g_{2r+1}\mathfrak{H}.o.\]

Hence, $d=2r+1=3$, $\mathfrak{R}=\{g_1,g_2,g_3\}$, the set
$\mathfrak{H}.o$ has $2$ elements and the set $X$ has $6$ elements.
\smallskip

In $X$ there are $20$ subsets of cardinality $3$. Among these, $8$ members intersect all three cosets $g\mathfrak{H}$, $\forall\ g\in\{g_1, g_2, g_3\}$. So $|\mathfrak{G}|$ is divisible by $8$ and by $20-8=12$. That is, $|\mathfrak{G}|$ is a multiple of $24$.
\smallskip

So $\mathfrak{G}$ is either the normalizer $\mathfrak{N}$ in $\mathrm{Aut}(X)$ of the partition
\[X=g_1\mathfrak{H}.o\cup g_2\mathfrak{H}.o\cup g_3\mathfrak{H}.o=\{a,a'\}\cup\{b,b'\}\cup\{c,c'\}\] or it is the index $2$ subgroup $\mathfrak{M}$ of $\mathfrak{N}$ generated by
\[(aa'),\ (bb'),\ (cc'),\ (abc)(a'b'c').\]

\smallskip

It remains to verify that $\mathfrak{M}$ as well as $\mathfrak{N}$ has $2$ orbits in the collection of $3$-point subsets of $X$ :
\smallskip

Let $Z$ be a subset of $X$ of cardinality $3$. Then
\smallskip

i) either $Z$ intersects all three : $\{a, a'\}$, $\{b, b'\}$, $\{c, c'\}$
\smallskip

ii) or $Z$ is disjoint with exactly one among $\{a, a'\}$, $\{b, b'\}$, $\{c, c'\}$.
\smallskip

In the first case, there is an element $g\in\langle (aa'), (bb'), (cc')\rangle$ such that $gZ=\{a, b,c\}$. In the latter, there is an element $g\in\langle(abc)(a'b'c')\rangle$ such that $gZ$ is either $\{a, a', b\}$ or $\{a, a', b'\}$. Then note that the cycle $(bb')$ transforms $\{a, a', b\}$ to $\{a, a', b'\}$.

\end{proof}
\smallskip

{\bf Lemma 6.6.} \emph{With the notations and assumptions of $(6.1)$, let $o\in X$ be a point and $\mathfrak{G}_o$ its stabilizer
in $\mathfrak{G}$. Suppose furthermore that $\mathfrak{G}$ acts transitively on $X$ and that
$\mathfrak{G}_o$ is a maximal subgroup of $\mathfrak{G}$.}
\smallskip

\emph{Then
$X=\mathbf{F}_8$, $\mathfrak{G}$ consists of either all affine
semi-linear transformations
\[
x\mapsto ax^{2^c}+b, \ \forall\ x\in\mathbf{F}_8
\] where
$a\in\mathbf{F}_8^{\times}$, $b\in\mathbf{F}_8$,
$c\in\mathbf{Z}/3\mathbf{Z}$ or only of the affine linear
transformations
\[
x\mapsto ax+b, \ \forall\ x\in\mathbf{F}_8
\] where
$a\in\mathbf{F}_8^{\times}$, $b\in\mathbf{F}_8$.}

\begin{proof} As in (6.3) one argues that there exists a normal subgroup $V$ of $\mathfrak{G}$ which has the following properties :
\smallskip

--- \emph{$\mathfrak{G}=V\mathfrak{G}_o$, $V\cap\mathfrak{G}_o=1$.}
\smallskip

--- \emph{$V$ acts simply transitively on $X$.}
\smallskip

--- \emph{$V$ is an $\mathbf{F}_p$-vector space for some prime $p$, and $\mathfrak{G}_o$ acts faithfully and irreducibly on $V$.}
\smallskip

\smallskip

We identify $V$ with $X$ by the bijection $v\mapsto v.o$.
\smallskip

Let $f=\mathrm{dim}\ V$. Then
$p^f=\mathrm{Card}(V)=\mathrm{Card}(X)=2d$. So
$p=2$ and $d=2^{f-1}$. Clearly, $f>1$.
\smallskip

\smallskip

--- \emph{Then $f>2$ }:
\smallskip

\smallskip

Suppose $f=2$. As $\mathfrak{G}_o$ acts irreducibly on $V$, it
cannot be a $2$-group. So $|\mathfrak{G}_o|$ is divisible by $3$. So
$\mathfrak{G}$ is $\mathrm{Aut}(X)=\mathfrak{S}_4$ or $\mathfrak{A}_4$. But both permute transitively the $2$-point subsets of $X$ rather than have $2$ orbits.
\smallskip

\smallskip

One has now $d=2^{f-1}\geq 4$.
\smallskip

Notice that every hyperplane of $V$ has $2^{f-1}=d$ elements. If $H_1,H_2$ are two distinct hyperplanes, the intersection $H_1\cap H_2$ has dimension
$f-2$ and cardinality $2^{f-2}=d/2$. And $H_2\backslash
H_1$ has $d/2$ elements. Given every $g\in\mathfrak{G}$, either $gH$ or $V\backslash
gH$ is a hyperplane. Hence $gH\backslash H$ has
$0,d$ or $d/2$ elements.
\smallskip

Fix a point $v\in V\backslash H$. The set
\[Y=\{v\}\cup H\backslash \{0\}\]
has $d$ elements. As $Y\backslash H$ consists only of one point, neither $Y$ nor its complement is a hyperplane.
\smallskip

Therefore, $\mathfrak{G}H$ and $\mathfrak{G}Y$ are these $2$ orbits of $\mathfrak{G}$ in the collection of $d$-point subsets of $X$.
\smallskip

\smallskip

--- \emph{Then $f=3$ }:
\smallskip

\smallskip

Assume $f>3$. Choose a point $u\in H\backslash\{0\}$. The set
\[
Z=\{v,u+v\}\cup H\backslash \{0,u\}
\]
has $d$ elements. But $Z$ is not a member of $\mathfrak{G}H$ or $\mathfrak{G}Y$. For, if $g$ is an element of $\mathfrak{G}$, then
\smallskip

--- the set $gH\backslash H$ has $0,d$ or $d/2$ elements,
\smallskip

--- the set  $gY\backslash H=g(Y\backslash g^{-1}H)$ has $1, d-1, d/2, (d/2)+1$ or $(d/2)-1$ elements,
\smallskip

--- while the set $Z\backslash H$ has $2$ elements.
\smallskip

Note that $2\notin \{0, 1, d, d-1, d/2, (d/2)+1, (d/2)-1\}$, as $d\geq 8$.

\smallskip

\smallskip

So $f=3$, and $\mathbf{P}(V)=\mathbf{P}^2$ is a projective plane over $\mathbf{F}_2$ which has $7$
$\mathbf{F}_2$-rational points. That is, $V$ has $7$ hyperplanes. These hyperplanes are permuted transitively by $\mathfrak{G}$. So $7$ divides
$|\mathfrak{G}|$ and $|\mathfrak{G}_o|$.
\smallskip

If one identifies $V$ with the underlying group of a finite field $\mathbf{F}_8$, a $7$-Sylow
subgroup of $\mathfrak{G}_o$ consists of all scalar multiplications
\[l_a: x\mapsto ax,\ \forall\ x\in\mathbf{F}_8\] where
$a\in\mathbf{F}_8^{\times}$.
\smallskip

\smallskip

--- \emph{The normalizer $\mathfrak{N}$ of the group $\{l_a,\ a\in\mathbf{F}_8^{\times}\}$ in $\mathrm{GL}(V)$ consists of all transformations of the form }:
\[x\mapsto aF^c(x),\ \forall\ x\in\mathbf{F}_8^{\times}\] \emph{where $a\in\mathbf{F}_8^{\times}$, $c\in\mathbf{Z}/3\mathbf{Z}$ and $F: x\mapsto x^2$, $\forall\ x\in\mathbf{F}_8$, is the Frobenius.}

\smallskip

\smallskip

Suppose that an element $g\in\mathrm{GL}(V)$ normalizes $\{l_a\}$. The characteristic polynomial of $g l_a g^{-1}$ on $V$ factors as :
\[
\mathrm{det}(T-g l_a g^{-1}, V)=\mathrm{det}(T-l_a,
V)=(T-a)(T-a^2)(T-a^4).\]

There is thus an element $c\in\mathbf{Z}/3\mathbf{Z}$
such that
\[g l_a g^{-1}=l_{F^c(a)}=F^c l_aF^{-c}.\] So $F^{-c}g$
commutes with all elements of the cyclic group $\{l_a\}$. So $F^{-c}g$ belongs to $\{l_a\}$. That is to say, $g$ is of the form
\[x\mapsto aF^c(x),\ \forall\ x\in\mathbf{F}_8\] for some
$a\in\mathbf{F}_8^{\times}$ and $c\in\mathbf{Z}/3\mathbf{Z}$. In particular, $\mathfrak{N}$ has $21$ elements.
\smallskip

\smallskip

--- \emph{The group $\mathfrak{G}_o$ is of odd order }:
\smallskip

\smallskip

As
$\mathfrak{G}_o$ is solvable, it has a Hall subgroup $\mathfrak{H}$
which is generated by
$\{l_a\}$ and a $2$-Sylow subgroup $\mathfrak{Q}$ of $\mathfrak{G}_o$.
Assume that $\mathfrak{H}$ is not of odd order. Thus $\mathfrak{H}$ is not a subgroup
of $\mathfrak{N}$. That is, $\{l_a\}$ is not normal in
$\mathfrak{H}$. So $\mathfrak{Q}$ is not of order $2$ or $4$. As $\mathrm{GL}(V)$ is of order $2^3.3.7$, $\mathfrak{Q}$ is of order $8$ and thus is normal in $\mathfrak{H}$. As $\mathfrak{Q}$ is $2$-Sylow in
$\mathrm{GL}(V)$, the center $\mathfrak{Z}$ of $\mathfrak{Q}$ is of order $2$, which is
normalized by $\{l_a\}$ and thus is centralized by $\{l_a\}$ and
thus is contained in $\mathfrak{N}$. This is absurd.
\smallskip

\smallskip

Now $|\mathfrak{G}_o|=7$ or $21$. In particular, $\{l_a\}$ is
normal in $\mathfrak{G}_o$. Hence, $\mathfrak{G}_o$ is contained in $\mathfrak{N}$. So $\mathfrak{G}_o$ is either $\mathfrak{N}$ or $\{l_a\}$. So $\mathfrak{G}$ is either $V\mathfrak{N}$, the group of all affine semi-linear transformations
\[
x\mapsto ax^{2^c}+b, \ \forall\ x\in\mathbf{F}_8
\] where
$a\in\mathbf{F}_8^{\times}$, $b\in\mathbf{F}_8$,
$c\in\mathbf{Z}/3\mathbf{Z}$ or it is $V\{l_a\}$ which
consists of all affine linear transformations
\[
x\mapsto ax+b,\ \forall\ x\in \mathbf{F}_8
\] where
$a\in\mathbf{F}_8^{\times}$, $b\in\mathbf{F}_8$.
\smallskip

\smallskip

It remains to verify that $V\{l_a\}$ as well as $V\mathfrak{N}$ has $2$ orbits in the collection of $4$-point subsets of $\mathbf{F}_8$ :
\smallskip

There are $70$ these subsets. Among them, the $7$ hyperplanes
and their complements, $14$ in number, form $1$ orbit under $V\{l_a\}$. For, if $E$ is one such set, then by a translation if necessary, one can transform $E$ to a hyperplane $H$. Now all $7$ hyperplanes are of the form $aH$, for $a\in\mathbf{F}_8^{\times}$.
\smallskip

Let $Y$ be a member in the rest $70-14=56$ sets of cardinality $4$. Let
$\mathfrak{S}$ denotes the normalizer of $Y$ in $V\{l_a\}$. Thus
$\mathfrak{S}$ is also a subgroup of $\mathrm{Aut}(Y)\times\mathrm{Aut}(X\backslash Y)=\mathfrak{S}_4\times\mathfrak{S}_4$, whose order is not a multiple of $7$.
So $\mathfrak{S}$ consists only of translations. So $\mathfrak{S}=1$ by the choice of $Y$. Therefore, these $56$ subsets are permuted simply transitively by $V\{l_a\}$.

\end{proof}
\smallskip

{\bf Lemma 6.7.} \emph{Let $X$ a finite nonempty set of even cardinality
$2d$. Let $\mathfrak{G}$ be a solvable subgroup of
$\mathrm{Aut}(X)$ which permutes transitively the subsets of $X$ of
cardinality $d$.}
\smallskip

\emph{Then up to equivalence $(X, \mathfrak{G})$ is one of the following }:
\smallskip

1) \emph{$X=\{1,2\}$, $\mathfrak{G}=\mathfrak{S}_2$.}
\smallskip

2) \emph{$X=\{1,2,3,4\}$, $\mathfrak{G}=\mathfrak{S}_4$ or
$\mathfrak{A}_4$.}

\begin{proof} Let $o\in X$ be a point and $\mathfrak{G}_o$ its
stabilizer in $\mathfrak{G}$. As in (6.3) there exists a normal subgroup $V$ of $\mathfrak{G}$ which has the following properties :
\smallskip

--- \emph{$\mathfrak{G}=V\mathfrak{G}_o$, $V\cap\mathfrak{G}_o=1$.}
\smallskip

--- \emph{$V$ acts simply transitively on $X$.}
\smallskip

--- \emph{$V$ is an $\mathbf{F}_p$-vector space for some prime $p$, and $\mathfrak{G}_o$ acts faithfully and irreducibly on $V$.}
\smallskip

We identify $V$ with $X$ by the bijection $v\mapsto v.o$. Let $f=\mathrm{dim}\ V$.
\smallskip

Then $p^f=\mathrm{Card}(V)=\mathrm{Card}(X)=2d$. So $p=2$ and $d=2^{f-1}$.
\smallskip

--- \emph{Then $f\leq 2$ }:
\smallskip

\smallskip

Fix a hyperplane $H$ in $V$ and let $v$ be a vector in the complement of $H$. By assumption each subset of $V$ of cardinality $d$ is a transform of $H$ by an element of $\mathfrak{G}$. These subsets are thus hyperplanes or complements of hyperplanes. But if $f>2$, the set
\[Y=\{v\}\cup H\backslash\{0\}\] is neither a hyperplane nor the complement of a hyperplane.
\smallskip

\smallskip

--- \emph{Case $f=1$ }:
\smallskip

\smallskip

The set $X$ has $2^f=2$ elements. And $\mathfrak{G}$ permutes the subsets of $X$ of cardinality $d=1$ transitively. Hence, $\mathfrak{G}=\mathrm{Aut}(X)$.

\smallskip

\smallskip

--- \emph{Case $f=2$ }:
\smallskip

\smallskip

The set $X$ as well as $V$ has $2^f=4$ elements. On $V$ the group $\mathfrak{G}_o$ acts irreducibly. So $\mathfrak{G}_o$ cannot be a $2$-group. So $\mathfrak{G}=V\mathfrak{G}_o$ is either $\mathrm{Aut}(X)$ or $\mathfrak{Alt}(X)$. Both do permute transitively the $2$-point subsets of $X$.

\end{proof}
\smallskip

{\bf Proposition 6.8.} \emph{Let $X$ a finite
set of even cardinality $2d\geq 4$. Let
$\mathfrak{G}$ be a solvable subgroup of $\{1,-1\}\times\mathrm{Aut}(X)$ which
permutes transitively the subsets of $X$ of cardinality $d$, where
$-1$ transforms every subset $Y$ of $X$ of
cardinality $d$ to $X\backslash Y$. Suppose furthermore that $\mathfrak{G}$ is not a subgroup of $\mathrm{Aut}(X)$.}
\smallskip

\emph{The following list enumerates such $(X, \mathfrak{G})$ up to equivalence }:
\smallskip

1) \emph{$X=\{1,2,3,4\}$, $\mathfrak{G}=\{1,-1\}\times\mathfrak{S}_4$.}
\smallskip

2) \emph{$X=\{1,2,3,4\}$, $\mathfrak{G}=\{1,-1\}\times\mathfrak{A}_4$.}
\smallskip

3) \emph{$X=\{1,2,3,4\}$, $\mathfrak{G}$ consists of $\mathfrak{A}_4$ and of all elements of the form $-1.\alpha$, where $\alpha$ is an odd permutation of $X$.}
\smallskip

4) \emph{$X=\{o,1,2,3\}$, $\mathfrak{G}=\{1,-1\}\times\mathrm{Aut}(\{1,2,3\})$.}
\smallskip

5) \emph{$X=\{o,1,2,3\}$, $\mathfrak{G}=\{1,-1\}\times\mathfrak{Alt}(\{1,2,3\})$.}
\smallskip

6) \emph{$X=\{o,1,2,3\}$, $\mathfrak{G}$ consists of $\mathfrak{Alt}(\{1,2,3\})$ and of all elements of the form $-1.\alpha$ where $\alpha$ is an odd permutation of $\{1,2,3\}$.}
\smallskip

7) \emph{$X=\{o\}\cup\mathbf{F}_5$, $\mathfrak{G}=\{1,-1\}\times\mathfrak{H}$ where $\mathfrak{H}$ is the group of affine linear transformations of $\mathbf{F}_5$.}

\begin{proof} Let $\mathfrak{H}=\mathfrak{G}\cap
\mathrm{Aut}(X)$, which is of index $2$ in $\mathfrak{G}$, as by assumption $\mathfrak{G}$ is not contained in $\mathrm{Aut}(X)$. It follows that
\smallskip

\smallskip

\emph{The collection of $d$-point subsets of $X$ are permuted by $\mathfrak{H}$ either transitively or in $2$ orbits of the same cardinality.}
\smallskip

\smallskip

--- \emph{Case where $\mathfrak{H}$ permutes transitively }:
\smallskip

\smallskip

By (6.7) the set $X$ has $4$ elements and $\mathfrak{H}=\mathrm{Aut}(X)$ or $\mathfrak{Alt}(X)$.
\smallskip

i) If $\mathfrak{H}=\mathrm{Aut}(X)$, then $\mathfrak{G}=\{1,-1\}\times\mathrm{Aut}(X)$.
\smallskip

ii) If $\mathfrak{H}=\mathfrak{Alt}(X)$, then either $\mathfrak{G}$ is $\{1,-1\}\times\mathfrak{Alt}(X)$ or it consists of $\mathfrak{Alt}(X)$ and of all elements of the form $-1.\alpha$, where $\alpha$ is an odd permutation of $X$.
\smallskip

\smallskip

--- \emph{Case where $\mathfrak{H}$ permutes with $2$ orbits of the same cardinality }:
\smallskip

\smallskip

By the proof of (6.1) precisely the following two occur :
\smallskip

--- \emph{$X=\{o,1,2,3\}$, $\mathfrak{H}$ fixes $o$ and on $\{1,2,3\}$ it is $\mathfrak{S}_3$ or $\mathfrak{A}_3$.}
\smallskip

--- \emph{$X=\{o\}\cup\mathbf{F}_5$, $\mathfrak{H}$ fixes $o$ and on $\mathbf{F}_5$ it is the
group of affine linear transformations.}
\smallskip

\smallskip

Let $\mathfrak{N}$ denote the normalizer of $\mathfrak{H}$ in
$\{1,-1\}\times\mathrm{Aut}(X)$.
\smallskip

Suppose first that $X=\{o, 1,2,3\}$.
\smallskip

Then $\mathfrak{N}=\{1,-1\}\times\mathrm{Aut}(\{1,2,3\})$ for both groups $\mathrm{Aut}(\{1,2,3\})$ and $\mathfrak{Alt}(\{1,2,3\})$. For, in $\mathrm{Aut}(X)$, the subgroup $\mathrm{Aut}(\{1,2,3\})$ is maximal and not normal.
\smallskip

Suppose next that $X=\{o\}\cup\mathbf{F}_5$.
\smallskip

Then $\mathfrak{N}=\{1,-1\}\times\mathfrak{H}$. Indeed, if $g\in\mathfrak{N}\cap\mathrm{Aut}(X)$, then $\mathfrak{H}g.o=g\mathfrak{H}.o=g.o$. So $g.o=o$ and $g$ normalizes the subset $\mathbf{F}_5$. As $\mathfrak{H}$ acts $2$-transitively on $\mathbf{F}_5$, there is an element $h\in\mathfrak{H}$ such that $hg$ fixes at least $2$ points of $\mathbf{F}_5$. In particular, $hg$ is of order $1$, $2$ or $3$. In $\mathfrak{H}$ the subgroup $\mathfrak{T}$ consisting of all translations is the unique $5$-Sylow subgroup. So $\mathfrak{T}$ is normalized and thus is centralized by $hg$. It follows that $hg$ fixes all points of $\mathbf{F}_5$. So $hg=1$ and $g=h^{-1}\in\mathfrak{H}$.
\smallskip

The pair $(X, \mathfrak{G})$ appears hence in the following list :
\smallskip

iii) $X=\{o,1,2,3\}$, $\mathfrak{G}=\{1,-1\}\times\mathrm{Aut}(\{1,2,3\})$.
\smallskip

iv) $X=\{o,1,2,3\}$, $\mathfrak{G}=\{1,-1\}\times\mathfrak{Alt}(\{1,2,3\})$.
\smallskip

v) $X=\{o,1,2,3\}$, $\mathfrak{G}$ consists of $\mathfrak{Alt}(\{1,2,3\})$ and of all elements of the form $-1.\alpha$, where $\alpha$ is an odd permutation of $\{1,2,3\}$.
\smallskip

vi) $X=\{o\}\cup\mathbf{F}_5$, $\mathfrak{G}=\{1,-1\}\times\mathfrak{H}$ where $\mathfrak{H}$ is the group of affine linear transformations of $\mathbf{F}_5$.

\smallskip

One inspects in each of the cases iii)--vi) that $\mathfrak{G}$ permutes the $d$-point subsets of $X$ transitively.

\end{proof}
\smallskip

{\bf Proposition 6.9.} \emph{Let $(S, \eta, s)$ be as in $\S 4$.
Then every $({}^2A_3,\alpha_2)$ over $\eta$ is elliptic. If $n>5$, then
$({}^2A_n,\alpha_{\frac{n+1}{2}})$ is not elliptic over $\eta$.}

\begin{proof} By (3.1), 8) and (6.8), $({}^2A_n,\alpha_{\frac{n+1}{2}})$ is not elliptic over $\eta$ if $n>5$. Suppose $n=3$ and suppose given a $({}^2A_3,\alpha_2)$ over $\eta$. Let
\[\rho_1: \pi_1(\eta, \overline{\eta})\to\{1,-1\}\] denote the index of ${}^2A_3$. Let
\[\rho_2: \pi_1(\eta, \overline{\eta})\to\pi_1(S, \overline{\eta})\to\mathfrak{Alt}(\{1,2,3\})\] be a surjective homomorphism (\S 4). Then
\[\rho=(\rho_1, \rho_2): \pi_1(\eta, \overline{\eta})\to\{1,-1\}\times\mathfrak{Alt}(\{1,2,3\})\] is surjective. Thus, by (6.8), 5)+(3.1), 8), it follows that
$({}^2A_3,\alpha_2)$ is elliptic over $\eta$.

\end{proof}
\smallskip

{\bf Proposition 6.10.} \emph{Let $(S, \eta, s)$, $\mathrm{char}(s)=\ell$, be as in $\S 4$. If $\ell=5$, then every $({}^2A_5,\alpha_3)$ over
$\eta$ is elliptic. When $(\ell,5)=1$, a $({}^2A_5,\alpha_3)$
over $\eta$ is elliptic if and only if ${}^2A_5$ is ramified over
$S$ and $\mathrm{Card}(k(s))$ mod $5$ generates
$\mathbf{F}_5^{\times}$.}

\begin{proof} Suppose given a $({}^2A_5, \alpha_3)$ over $\eta$. Let
\[\rho_1: \pi_1(\eta, \overline{\eta})\to \{1,-1\}\] denote its index.
By (3.1), 8)+(6.8), 2) this $({}^2A_5,\alpha_3)$
is elliptic if and only if there is a surjective homomorphism :
\[ \rho=(\rho_1,\rho_2): \pi_1(\eta,\overline{\eta})\to
\{1,-1\}\times\mathfrak{H}=:\mathfrak{G}\] where $\mathfrak{H}$ consists
of all affine linear transformations of $\mathbf{F}_5$.
\smallskip

Notice that $\mathfrak{H}$ is a quotient of
$\pi_1(\eta,\overline{\eta})$ if and only if one of the following two holds :
\smallskip

--- $\ell=5$ (4.1).
\smallskip

--- $(\ell, 5)=1$ and $\mathrm{Card}(k(s))$ mod $5$ generates $\mathbf{F}_5^{\times}$ (4.2), 2).
\smallskip

Note also that
\smallskip

\smallskip

\emph{If $(\ell, 5)=1$ and if $\rho_1$ is unramified over $S$, then $({}^2A_5, \alpha_3)$ is not elliptic over $\eta$.}
\smallskip

\smallskip

Otherwise, the image of the inertia subgroup of $\pi_1(\eta, \overline{\eta})$ in $\mathfrak{G}$ would be the subgroup $\mathfrak{T}$ of all translations of $\mathbf{F}_5$. But $\mathfrak{G}/\mathfrak{T}$ is not cyclic.
\smallskip

\smallskip

--- \emph{Case $\ell=5$, $\rho_1$ unramified over $S$ }:
\smallskip

\smallskip

Let $\pi\in\Gamma(S, \mathcal{O}_S)$ be a uniformizer. Then
\[\eta[z, x]/(z^4-\pi, x^5-x-z^{-1})\] is connected, totally ramified over $S$ and Galois over $\eta$ with Galois group $\mathfrak{H}$. If its corresponding monodromy representation is
\[\rho_2: \pi_1(\eta, \overline{\eta})\to\mathfrak{H},\] then
\[\rho=(\rho_1, \rho_2): \pi_1(\eta, \overline{\eta})\to \mathfrak{G}\] is surjective.

\smallskip

\smallskip

--- \emph{Case $\ell=5$, $\rho_1$ ramified over $S$ }:
\smallskip

\smallskip

Let $\pi\in\Gamma(S, \mathcal{O}_S)$ be a uniformizer. Let $S'$ be the spectra of a discrete valuation ring such that $S'$ is finite \'{e}tale Galois over $S$ with cyclic Galois group of order $4$, $\eta'$ (resp. $s'$) the generic (resp. closed) point of $S'$, $\zeta\in\mathrm{Gal}(S'/S)$ a generator and let $u'\in\Gamma(S', \mathcal{O}_{S'})^{\times}$ a unit such that the images of $u',\cdots, \zeta^3(u')$ in $k(s')$ form a normal base over $k(s)$. Then
\[\eta'[x_1,\cdots, x_4]/(x_1^5-x_1-\zeta(u')\pi^{-1},\cdots, x_4^5-x_4-\zeta^4(u')\pi^{-1})\] is connected and Galois over $\eta$ with Galois group $\mathfrak{H}$. If its corresponding monodromy representation is
\[\rho_2: \pi_1(\eta, \overline{\eta})\to\mathfrak{H},\] then
\[\rho=(\rho_1, \rho_2): \pi_1(\eta, \overline{\eta})\to\mathfrak{G}\] is surjective.
\smallskip

\smallskip

--- \emph{Case $(\ell, 5)=1$, $\mathrm{Card}(k(s))$ mod $5$ generates $\mathbf{F}_5^{\times}$, $\rho_1$ ramified over $S$ }:
\smallskip

\smallskip

By (4.2), 2) $\mathfrak{H}$ is realizable as a tame quotient of $\pi_1(\eta, \overline{\eta})$, say
\[\rho_2: \pi_1(\eta, \overline{\eta})\to\mathfrak{H}.\] Now
\[\rho=(\rho_1, \rho_2): \pi_1(\eta, \overline{\eta})\to\mathfrak{G}\] is surjective.

\end{proof}
\smallskip

\smallskip


\section{Type $B$}
\smallskip

\smallskip

\smallskip

Let $S,\eta,s$, $\mathrm{char}(s)=\ell$, be as in $\S 4$.
\smallskip

Let $n$ be an integer $\geq 3$. Let $e_1,\cdots,e_n$ be the standard basis
of $\mathbf{Z}^n$. We denote the group of diagonal (resp. monomial) matrices in $\mathrm{GL}_n(\mathbf{Z})$ by $\mathfrak{D}$ (resp. $\mathfrak{M}$). Let $\mathfrak{W}=\mathfrak{D}\mathfrak{M}$.

\smallskip

\smallskip

\smallskip

{\bf Proposition 7.1.} \emph{Suppose that $k(s)$ is of characteristic $\ell=2$. Then
$(B_n,\alpha_n)$ is elliptic over $\eta$.}

\begin{proof} Let $\zeta\in\mathrm{GL}_n(\mathbf{Z})$ be such that
\[\zeta: e_1\mapsto e_2,\ e_2\mapsto e_3,\ \cdots, \ e_n\mapsto e_1.\]
Let $\mathfrak{G}$ be the subgroup of $\mathfrak{W}$ generated by $\zeta$ and all diagonal matrices. The group $\mathfrak{G}$ permutes the vectors
\[\pm e_1\pm\cdots \pm e_n\] transitively. Moreover, $\mathfrak{G}$ is a quotient of
$\pi_1(\eta,\overline{\eta})$ (4.1). So $(B_n,\alpha_n)$ is
elliptic over $\eta$ (3.1), 2).

\end{proof}
\smallskip

{\bf Proposition 7.2.} \emph{The pair $(B_3,\alpha_3)$ is elliptic
over $\eta$.}

\begin{proof} The following two elements of $\mathrm{GL}_3(\mathbf{Z})$
\[a: e_1\mapsto e_1,\ e_2\mapsto e_3,\ e_3\mapsto -e_2\]
\[b: e_1\mapsto -e_1,\ e_2\mapsto e_2,\ e_3\mapsto e_3\]
satisfy the relations
\[a^4=b^2=1,\  ab=ba.\]

The group $\mathfrak{G}$ they generate is isomorphic to $\mathbf{Z}/4\mathbf{Z}\times\mathbf{Z}/2\mathbf{Z}$. One verifies that $\mathfrak{G}$
permutes the vectors
\[\pm e_1\pm e_2\pm e_3\] simply transitively. Moreover, $\mathfrak{G}$ is a quotient of $\pi_1(\eta, \overline{\eta})$. Indeed, let
\[\rho_1: \pi_1(\eta, \overline{\eta})\to\pi_1(S, \overline{\eta})\to\mathbf{Z}/4\mathbf{Z}\] be a surjective homomorphism (\S 4) and let
\[\rho_2: \pi_1(\eta, \overline{\eta})\to\mathbf{Z}/2\mathbf{Z}\] be the monodromy representation corresponding to the quadratic extension
\[k(\eta)[x]/(x^2-\pi)\] of $k(\eta)$, where $\pi\in\Gamma(S,\mathcal{O}_S)$ is a uniformizer. Then
\[\rho=(\rho_1, \rho_2): \pi_1(\eta, \overline{\eta})\to\mathbf{Z}/4\mathbf{Z}\times\mathbf{Z}/2\mathbf{Z}\] is surjective. By (3.1), 2) $(B_3,\alpha_3)$ is thus elliptic over $\eta$.

\end{proof}
\smallskip

{\bf Proposition 7.3.} \emph{The pair
$(B_4,\alpha_4)$ is elliptic over $\eta$.}

\begin{proof} By (7.1) one can suppose $\ell>2$.
\smallskip

The following elements of $\mathrm{GL}_4(\mathbf{Z})$
\[a: e_1\mapsto e_2,\ e_2\mapsto -e_1,\ e_3\mapsto e_3,\ e_4\mapsto e_4\]
\[b: e_1\mapsto e_1,\ e_2\mapsto e_2,\ e_3\mapsto e_4,\ e_4\mapsto -e_3\]
\[c: e_1\mapsto e_2,\ e_2\mapsto e_3,\ e_3\mapsto e_4,\ e_4\mapsto -e_1\]
\[d: e_1\mapsto e_3,\ e_2\mapsto -e_4,\ e_3\mapsto -e_1,\ e_4\mapsto e_2\]
satisfy the relations
\[ a^4=b^4=1,\ ab=ba.\]
\[ c^8=d^4=1,\ cdc^{-1}=d^{-1}.\]

The group $\mathfrak{G}_1$ generated by $\{a,b\}$ is isomorphic to
$\mathbf{Z}/4\mathbf{Z}\times\mathbf{Z}/4\mathbf{Z}$. The group $\mathfrak{G}_2$ generated by $\{
c,d\}$ is quaternion of order $16$. Both permute simply transitively
the vectors
\[\pm e_1\pm e_2\pm e_3\pm e_4.\]

If $\mathrm{Card}(k(s))\equiv 1$ mod $4$ (resp.
$\mathrm{Card}(k(s))\equiv -1$ mod $4$), then $\mathfrak{G}_1$ (resp.
$\mathfrak{G}_2$) is a quotient of
$\pi_1^t(\eta,\overline{\eta})$ (\S 4).
So $(B_4,\alpha_4)$ is elliptic over $\eta$ (3.1), 2).

\end{proof}
\smallskip

{\bf Proposition 7.4.} \emph{Suppose $\ell>2$, $n>4$. Then
$(B_n,\alpha_n)$ is not elliptic over $\eta$.}

\begin{proof} By (3.1), 2) $(B_n, \alpha_n)$ is elliptic if and only if there is
a representation
\[\rho: \pi_1(\eta,\overline{\eta})\to\mathfrak{W}\] whose image permutes transitively
the vectors
\[\pm e_1\pm\cdots\pm e_n.\]

Suppose that such a representation exists. Let $\mathfrak{G}$ be its image and $I$ (resp. $P$) the image of the inertia (resp. wild inertia) subgroup. Let
\[X=\{\pm e_1\pm\cdots\pm e_n\}.\]

Observe that if an element $g$ of $\mathfrak{W}$ fixes all points of $X$, then $g=1$.
\smallskip

\smallskip

--- \emph{One has $P=1$ }:
\smallskip

\smallskip

For, $P$ being normal in $\mathfrak{G}$, the $P$-orbits in $X$
all have the same cardinality, say $r$, which divides both $2^n$ and $\mathrm{Card}(P)$. So $r=1$. So $P=1$.
\smallskip

\smallskip

Thus $I=I/P$ is cyclic.

\smallskip

\smallskip

--- \emph{The cyclic group $I$ is a $2$-group }:
\smallskip

\smallskip

The maximal odd order subgroup of $I$, $I'$, is normal in $\mathfrak{G}$. Thus the $I'$-orbits in $X$ all have the same cardinality, say $r'$, which divides both $2^n$ and $\mathrm{Card}(I')$. That is, $r'=1$. So $I'=1$.
\smallskip

\smallskip

Now, as $\mathfrak{G}/I$ is cyclic, $\mathfrak{G}$ has a unique $2$-Sylow subgroup.
\smallskip

\smallskip

--- \emph{The unique $2$-Sylow subgroup $\mathfrak{H}$ of $\mathfrak{G}$ acts transitively on $X$ }:
\smallskip

\smallskip

As $\mathfrak{H}$ is normal in $\mathfrak{G}$, the $\mathfrak{H}$-orbits in $X$ all have the same cardinality. These orbits, which form a set of cardinality dividing $2^n$, are permuted transitively by the quotient $\mathfrak{G}/\mathfrak{H}$. So there is only one orbit.
\smallskip

\smallskip

To $\mathfrak{H}$ there corresponds a subextension $k(\eta')/k(\eta)$ of $k(\overline{\eta})/k(\eta)$ so that $\mathfrak{H}$ is the image of the composition
\[\pi_1(\eta', \overline{\eta})\to\pi_1(\eta, \overline{\eta})\stackrel{\rho}{\longrightarrow}\mathfrak{G}.\]

Replacing $\eta$ by $\eta'$ if necessary, one can assume that $\mathfrak{G}=\mathfrak{H}$ is a $2$-group.
\smallskip

Consider the exact sequence
\[1\to I\cap\mathfrak{D}\to\mathfrak{G}\cap\mathfrak{D}\to \mathfrak{G}/I.\]

Notice that
\smallskip

--- the group $I\cap\mathfrak{D}$ has $\leq 2$ elements. For, $I\cap\mathfrak{D}$ is both cyclic and an elementary $2$-group.
\smallskip

--- the image $\mathfrak{Q}$ of $\mathfrak{G}\cap\mathfrak{D}$ in $\mathfrak{G}/I$ has $\leq 2$ elements. For, being a subgroup of $\mathfrak{G}/I$, $\mathfrak{Q}$ is cyclic. And being a quotient of $\mathfrak{G}\cap\mathfrak{D}$, $\mathfrak{Q}$ is an elementary $2$-group.

\smallskip

Therefore, $\mathfrak{G}\cap\mathfrak{D}$ is of order $1, 2$ or $4$.

\smallskip

The quotient $\mathfrak{G}/(\mathfrak{G}\cap\mathfrak{D})$ is isomorphic to a group of monomial matrices. Thus this $2$-group is of order $2^e$ for an integer $e\leq \mathrm{ord}_2(n!)$. One has $\mathrm{ord}_2(n!)\leq n-1$, where the equality holds if and only if $n$ is a power of $2$.

\smallskip

As $\mathfrak{G}$ acts transitively on $X$, one of the following three holds :
\smallskip

1) $\mathfrak{G}\cap\mathfrak{D}$ has $4$ elements, $\mathfrak{G}/(\mathfrak{G}\cap\mathfrak{D})$ has $2^{n-2}$ elements and is $2$-Sylow in $\mathfrak{M}$, $n$ is not a power of $2$.
\smallskip

2) $\mathfrak{G}\cap\mathfrak{D}$ has $2$ elements, $\mathfrak{G}/(\mathfrak{G}\cap\mathfrak{D})$ has $2^{n-1}$ elements and is $2$-Sylow in $\mathfrak{M}$, $n$ is a power of $2$.
\smallskip

3) $\mathfrak{G}\cap\mathfrak{D}$ has $4$ elements, $\mathfrak{G}/(\mathfrak{G}\cap\mathfrak{D})$ has $2^{n-2}$ or $2^{n-1}$ elements and is of index $\leq 2$ in a $2$-Sylow subgroup of $\mathfrak{M}$, $n$ is a power of $2$.
\smallskip

Next, from the exact sequence
\[1\to I/(I\cap\mathfrak{D})\to \mathfrak{G}/(\mathfrak{G}\cap\mathfrak{D})\to\mathfrak{G}/I(\mathfrak{G}\cap\mathfrak{D})\to 1\]
one deduces that $\mathfrak{G}/(\mathfrak{G}\cap\mathfrak{D})$ does not have elementary $2$-subgroups of $2$-rank $\geq 3$. So
\smallskip

--- $n\leq 5$ in case 1)
\smallskip

--- $n\leq 4$ in case 2)
\smallskip

--- $n\leq 4$ in case 3)
\smallskip

It remains to consider the case $n=5$ :
\smallskip

By 1) then $|\mathfrak{G}\cap\mathfrak{D}|=4$, $|I\cap\mathfrak{D}|=2$, $|\mathfrak{G}|=32$.
\smallskip

The group $\mathfrak{G}/(\mathfrak{G}\cap\mathfrak{D})$ has $8$ elements and is $2$-Sylow in $\mathfrak{M}=\mathfrak{S}_5$ and is an extension of the cyclic group $\mathfrak{G}/I(\mathfrak{G}\cap\mathfrak{D})$ by the cyclic group $I/(I\cap\mathfrak{D})$. So $I/(I\cap\mathfrak{D})$ has $4$ elements. So $I$ has $8$ elements and $\mathfrak{G}/I$ has $4$ elements.
\smallskip

Let $t$ be a generator of $I$. Choose an element $f$ of $\mathfrak{G}$ such that its image in $\mathfrak{G}/I$ is a generator. Then
$ftf^{-1}=t^q$ for an odd integer $q$. So $f^2$ commutes with $t$, for
$f^2tf^{-2}=t^{q^2}=t$, as $q^2\equiv 1$ mod $8$.
\smallskip

Observe that $\mathfrak{G}$ normalizes the set
\[Y=\{e_1,\cdots, e_5, -e_1,\cdots, -e_5\}\] in which one and only one $I$-orbit $O=-O$ has $8$ elements. Let $O'=Y\backslash O$, which consists of two eigenvectors of $t$.
\smallskip

Now $f$ normalizes $O$ as well as $O'$. On $O'$, $f^2$ acts as the identity. On $O$, $f^2$ acts as $t^r$ for an even integer $r$ since it commutes with $t$. So $f^2=t^r\in I$. So $\mathfrak{G}/I$ has $\leq 2$ elements. A contradiction.

\end{proof}

\smallskip

\smallskip


\section{Type $C$}
\smallskip

\smallskip

Let $(S,\eta,s)$ be as in \S 4.
\smallskip

\smallskip

\smallskip

{\bf Proposition 8.1.} \emph{For every integer $n\geq 1$,
$(C_n,\alpha_1)$ is elliptic over $\eta$.}

\begin{proof} Let $\zeta, \tau\in\mathrm{GL}_n(\mathbf{Z})$ be such that
\[\zeta: e_1\mapsto e_2, \ e_2\mapsto e_3, \cdots, e_n\mapsto e_1\]
\[\tau: e_1\mapsto -e_1,\ e_i\mapsto e_i,\ \forall\ i>1\]
where $e_1,\cdots, e_n$ denote the standard basis of $\mathbf{Z}^n$.
\smallskip

The cyclic group $\langle\tau\zeta\rangle$ generated by $\tau\zeta$ permutes the vectors
\[e_1,\cdots, e_n, -e_1,\cdots, -e_n\] simply transitively. And, $\langle\tau\zeta\rangle$ is a quotient of $\pi_1(\eta, \overline{\eta})$ (\S 4). So $(C_n,\alpha_1)$ is
elliptic over $\eta$ (3.1), 3).

\end{proof}
\smallskip

\smallskip

\smallskip


\section{Type $D$}
\smallskip

\smallskip

Let $S,\eta,s$, $\mathrm{char}(s)=\ell$, be as in \S 4.
\smallskip

Let $n$ be an integer $\geq 4$. Let $e_1,\cdots,e_n$ be the standard basis
of $\mathbf{Z}^n$. We denote the group of diagonal (resp. monomial) matrices in $\mathrm{GL}_n(\mathbf{Z})$ by $\mathfrak{D}$ (resp. $\mathfrak{M}$). Let $\mathfrak{D}_1$ be the subgroup of $\mathfrak{D}$ consisting of all diagonal matrices of determinant $1$.
\smallskip

Let $\mathfrak{W}=\mathfrak{D}\mathfrak{M}$ and $\mathfrak{W}_1=\mathfrak{D}_1\mathfrak{M}$.

\smallskip

\smallskip

{\bf Proposition 9.1.} \emph{Suppose that $n$ is even. Then $(D_n, \alpha_1)$ is elliptic over $\eta$.}

\begin{proof} As $n$ is even, the diagonal matrix $-1\in\mathrm{GL}_n(\mathbf{Z})$ has determinant $1$. Let $\zeta\in\mathrm{GL}_n(\mathbf{Z})$ be such that
\[\zeta: e_1\mapsto e_2, \cdots, e_n\mapsto e_1.\]

The subgroup $\mathfrak{G}$ of $\mathfrak{W}_1$ generated by $\{\zeta, -1\}$ permutes the vectors
\[e_1,\cdots, e_n, -e_1,\cdots, -e_n\] simply transitively. And $\mathfrak{G}$, which
is isomorphic to $\mathbf{Z}/n\mathbf{Z}\times\mathbf{Z}/2\mathbf{Z}$, is a quotient of $\pi_1(\eta, \overline{\eta})$. Indeed, let
\[\rho_1: \pi_1(\eta, \overline{\eta})\to\pi_1(S, \overline{\eta})\to\mathbf{Z}/n\mathbf{Z}\] be a surjective homomorphism (\S 4) and let
\[\rho_2: \pi_1(\eta, \overline{\eta})\to\mathbf{Z}/2\mathbf{Z}\] be the monodromy representation corresponding to the quadratic extension
\[k(\eta)[x]/(x^2-\pi)\] of $k(\eta)$, where $\pi\in\Gamma(S,\mathcal{O}_S)$ is a uniformizer. Then
\[\rho=(\rho_1, \rho_2): \pi_1(\eta, \overline{\eta})\to\mathbf{Z}/n\mathbf{Z}\times\mathbf{Z}/2\mathbf{Z}\] is surjective. So $(D_n, \alpha_1)$ is elliptic over $\eta$ (3.1), 4).

\end{proof}

\smallskip

{\bf Proposition 9.2.} \emph{Suppose $\ell=2$. Then $(D_n, \alpha_1)$ is elliptic over $\eta$.}

\begin{proof} Let $\zeta\in\mathrm{GL}_n(\mathbf{Z})$ be such that
\[\zeta: e_1\mapsto e_2,\cdots, e_n\mapsto e_1.\] The subgroup $\mathfrak{G}$ of $\mathfrak{W}_1$ generated by $\zeta$ and all diagonal matrices of determinant $1$ permutes transitively the vectors
\[e_1,\cdots, e_n, -e_1,\cdots, -e_n.\] Moreover, $\mathfrak{G}$ is a quotient of $\pi_1(\eta, \overline{\eta})$ (4.1). So $(D_n, \alpha_1)$ is elliptic over $\eta$ (3.1), 4).

\end{proof}

\smallskip

{\bf Lemma 9.3.} \emph{Every odd order subgroup of $\mathfrak{W}_1$ is conjugate to a subgroup of $\mathfrak{M}$.}

\begin{proof} Let $\mathfrak{H}$ be an odd order subgroup of $\mathfrak{W}_1$. Consider the split exact sequence
\[1\to\mathfrak{D}_1\to\mathfrak{W}_1\to \mathfrak{M}\to 1\] and let $\mathfrak{H}'$ be the image of $\mathfrak{H}$ in $\mathfrak{M}$. Then $\mathfrak{D}_1\mathfrak{H}=\mathfrak{D}_1\mathfrak{H}'=:\mathfrak{G}$. As $H^1(\mathfrak{H}', \mathfrak{D}_1)=0$, every two splittings of the exact sequence
\[1\to\mathfrak{D}_1\to\mathfrak{G}\to\mathfrak{H}'\to 1,\] especially the ones corresponding to $\mathfrak{H}$ and to $\mathfrak{H'}$, are conjugate to each other by an element of $\mathfrak{G}$.

\end{proof}

\smallskip

{\bf Proposition 9.4.} \emph{Suppose $\ell>2$ and that $n$ is odd $\geq 3$. Then $(D_n, \alpha_1)$ is not elliptic over $\eta$.}

\begin{proof} This conclusion holds when $n=3$ by (5.3)+(5.4), as $(D_3, \alpha_1)=(A_3, \alpha_2)$.
\smallskip

Suppose $n\geq 5$.
\smallskip

By (3.1), 4) $(D_n, \alpha_1)$ is elliptic if and only if there is a representation
\[\rho: \pi_1(\eta, \overline{\eta})\to\mathfrak{W}_1\] whose image acts transitively on the set
\[X=\{e_1,\cdots, e_n, -e_1,\cdots, -e_n\}.\]

Suppose that $n$ is the smallest odd integer $\geq 5$ for a representation $\rho$ as such exists. Let $\mathfrak{G}$ be its image and $I$ (resp. $P$) the image of the inertia (resp. wild inertia) subgroup of $\pi_1(\eta, \overline{\eta})$. By (9.3) one may assume $P$ to be a subgroup of monomial matrices. Thus, $P$ normalizes
\[X_{+}=\{e_1,\cdots, e_n\}.\]

Notice that, $P$ being normal in $\mathfrak{G}$, the $P$-orbits in $X$ all have the same cardinality, say $r$, which divides both $\mathrm{Card}(P)$ and $n$. Let $d=n/r$. Let $E_1,\cdots, E_d$ be the $P$-orbits in $X_{+}$; the other $P$-orbits are $-E_1,\cdots, -E_d$. These $P$-orbits are permuted transitively by $\mathfrak{G}$.
\smallskip

Let $g=\delta p$ be an element of $\mathfrak{G}$, where $\delta\in\mathfrak{D}_1$, $p\in\mathfrak{M}$. Let $E$ be a $P$-orbit in $X_{+}$.
Suppose that
$g(E)=\chi E'$, where $\chi\in\{1,-1\}$, $E'\subset X_{+}$. Then
$p(E)=\chi \delta(E')$. Namely, $p(E)=E'$, $\delta|E'=\chi$.
\smallskip

It follows that when $E_1,\cdots, E_d$ is considered as a base of a free $\mathbf{Z}$-module $\mathbf{Z}^d$ the permutation action of $\mathfrak{G}$ on
\[\{E_1,\cdots, E_d, -E_1,\cdots, -E_d\}\] induces a representation of $\mathfrak{G}$ in $\mathrm{GL}_d(\mathbf{Z})$ whose image lies in the group generated by the diagonal matrices of determinant $1$ and monomial matrices.
\smallskip

Thus, in view of the choice of $n$, one has $d=n$. So $r=1$, $P=1$.
\smallskip

So $I=I/P$ is cyclic. The maximal odd order subgroup of $I$, which is normal in $\mathfrak{G}$, is $1$ by the same argument as for $P$. That is, $I$ is a cyclic $2$-group.
\smallskip

As $I$ is normal in $\mathfrak{G}$ and is commutative, the $I$-orbits in $X$ all have the same cardinality $|I|$, which divides $2n$. So $|I|=1$ or $2$. So $\mathfrak{G}$ is commutative of order $2n$. The unique index $2$ subgroup of $\mathfrak{G}$ is again $1$ by the same argument as for $P$. So $n=1$. A contradiction.

\end{proof}

\smallskip

{\bf Proposition 9.5.} \emph{The pairs $(D_n, \alpha_{n-1})$ and $(D_n, \alpha_n)$ are elliptic over $\eta$ if $n=4$ or $5$.}

\begin{proof} By comparing (3.1), 5) with (3.1), 2), it is evident that $(D_n, \alpha_{n-1})$ and $(D_n, \alpha_n)$ are elliptic if $(B_{n-1}, \alpha_{n-1})$ is elliptic. One now applies (7.2)+(7.3).

\end{proof}

\smallskip

{\bf Proposition 9.6.} \emph{Suppose $\ell=2$. Then $(D_n,\alpha_{n-1})$ and
$(D_n,\alpha_n)$ are elliptic over $\eta$.}

\begin{proof} Let $\zeta\in\mathrm{GL}_n(\mathbf{Z})$ be such that
\[\zeta: e_1\mapsto e_2,\cdots, e_n\mapsto e_1.\] Let $\mathfrak{G}$ be the subgroup of $\mathfrak{W}_1$ generated by $\zeta$ and $\mathfrak{D}_1$. The group $\mathfrak{G}$ acts transitively on
\[ \{s_1e_1+\cdots +s_ne_n,\ s_i\in\{1,-1\},\ s_1\cdots s_n=-1\}\]
and on
\[ \{s_1e_1+\cdots +s_ne_n,\ s_i\in\{1,-1\},\ s_1\cdots s_n=1\}.\]

Moreover, $\mathfrak{G}$ is
a quotient of $\pi_1(\eta,\overline{\eta})$. Indeed, let $S'$ be the spectra of a discrete valuation ring such that $S'$ is finite \'{e}tale Galois over $S$ with cyclic Galois group of order $n$ (\S 4), $\eta'$ (resp. $s'$) the generic (resp. closed) point of $S'$, $\pi\in\Gamma(S, \mathcal{O}_S)$ a uniformizer and let $u'\in\Gamma(S',\mathcal{O}_{S'})^{\times}$ be a unit such that the images of
$u',\zeta(u'),\cdots,\zeta^{n-1}(u')$ in $k(s')$ form a normal basis
over $k(s)$. Put $b':=1+u'\pi$. Then
\[\eta'[z_1,\cdots,z_n]/(z_1^2-\frac{\zeta(b')}{b'},\cdots,
z_n^2-\frac{\zeta^n(b')}{\zeta^{n-1}(b')}, 1-z_1\cdots z_n)\] is connected and
Galois over $\eta$ with Galois group $\mathfrak{G}$.
So $(D_n,\alpha_{n-1})$ and $(D_n,\alpha_n)$ are elliptic
(3.1), 4), 5).

\end{proof}

\smallskip

{\bf Proposition 9.7.} \emph{Suppose $\ell>2$, $n>5$. Then $(D_n, \alpha_{n-1})$ and $(D_n, \alpha_n)$ are not elliptic over $\eta$.}

\begin{proof} It suffices to consider $(D_n,\alpha_n)$ only. The same argument applies for $(D_n,\alpha_{n-1})$.
By (3.1), 5) $(D_n, \alpha_n)$ is elliptic if and only if there is a representation
\[\rho: \pi_1(\eta, \overline{\eta})\to \mathfrak{W}_1\] whose image acts transitively on the set
\[X=\{s_1e_1+\cdots+s_ne_n, s_i\in\{1, -1\}, s_1\cdots s_n=1\}.\]

Suppose that such a representation $\rho$ exists. Let $\mathfrak{G}$ be its image and $I$ (resp. $P$) the image in $\mathfrak{G}$ of the inertia (resp. wild inertia) subgroup of $\pi_1(\eta, \overline{\eta})$.

\smallskip

As in (7.4), extending if necessary $k(\eta)$ to a finite extension $k(\eta')$ which is unramified over $S$ and of odd degree over $k(\eta)$, one can assume that $\mathfrak{G}$ is a $2$-group. In particular, $P=1$ and $I$ is cyclic.

\smallskip

Consider the exact sequence
\[1\to I\cap\mathfrak{D}_1\to\mathfrak{G}\cap\mathfrak{D}_1\to \mathfrak{G}/I.\]

As both $I\cap\mathfrak{D}_1$ and $\mathfrak{G}/I$ are cyclic, the elementary $2$-group $\mathfrak{G}\cap\mathfrak{D}_1$ is of order $1, 2$ or $4$.
The quotient $\mathfrak{G}/(\mathfrak{G}\cap\mathfrak{D}_1)$, which is isomorphic to a group of monomial matrices, is of order $2^e$ for an integer $e\leq\mathrm{ord}_2(n!)$. Notice that $\mathrm{ord}_2(n!)\leq n-1$, where the equality holds if and only if $n$ is a power of $2$.
\smallskip

As $\mathfrak{G}$ acts transitively on $X$, one of the following five holds :
\smallskip

1) $|\mathfrak{G}\cap\mathfrak{D}_1|=2$ or $4$, $\mathfrak{G}/(\mathfrak{G}\cap\mathfrak{D}_1)$ has $2^{n-2}$ elements and is $2$-Sylow in $\mathfrak{M}$, $n$ is not a power of $2$.
\smallskip

2) $|\mathfrak{G}\cap\mathfrak{D}_1|=4$, $\mathfrak{G}/(\mathfrak{G}\cap\mathfrak{D}_1)$ has $2^{n-3}$ elements and is of index $\leq 2$ in a $2$-Sylow subgroup of $\mathfrak{M}$, $n$ is not a power of $2$.
\smallskip

3) $|\mathfrak{G}\cap\mathfrak{D}_1|=1,2$ or $4$, $\mathfrak{G}/(\mathfrak{G}\cap\mathfrak{D}_1)$ has $2^{n-1}$ elements and is $2$-Sylow in $\mathfrak{M}$, $n$ is a power of $2$.
\smallskip

4) $|\mathfrak{G}\cap\mathfrak{D}_1|=2$ or $4$, $\mathfrak{G}/(\mathfrak{G}\cap\mathfrak{D}_1)$ has $2^{n-2}$ elements and is of index $\leq 2$ in a $2$-Sylow subgroup of $\mathfrak{M}$, $n$ is a power of $2$.
\smallskip

5) $|\mathfrak{G}\cap\mathfrak{D}_1|=4$, $\mathfrak{G}/(\mathfrak{G}\cap\mathfrak{D}_1)$ has $2^{n-3}$ elements and is of index $1, 2$ or $4$ in a $2$-Sylow subgroup of $\mathfrak{M}$, $n$ is a power of $2$.
\smallskip

Next, the exact sequence
\[1\to I/(I\cap\mathfrak{D}_1)\to \mathfrak{G}/(\mathfrak{G}\cap\mathfrak{D}_1)\to\mathfrak{G}/I(\mathfrak{G}\cap\mathfrak{D}_1)\to 1\] implies that $\mathfrak{G}/(\mathfrak{G}\cap\mathfrak{D}_1)$ does not contain elementary $2$-groups of $2$-rank $\geq 3$. So
\smallskip

--- $n\leq 5$ in case 1)
\smallskip

--- $n\leq 6$ in case 2)
\smallskip

--- $n\leq 4$ in case 3)
\smallskip

--- $n\leq 4$ in case 4)
\smallskip

--- $n\leq 8$ in case 5)
\smallskip

\smallskip

Now as $P=1$, $\mathfrak{G}$ is a quotient of $\pi_1^t(\eta, \overline{\eta})=\langle F, T\rangle$.
Let $t\in I$ (resp. $f$) be the image of $T$ (resp. $F$). Then $ftf^{-1}=t^q$ where $q=\mathrm{Card}(k(s))$ is a power of $\ell$.
\smallskip

Notice that $\mathfrak{W}_1$ does not have elements of order $16$ for $n\leq 8$. So $t^8=f^8=1$ and $|\mathfrak{G}|$ divides $64$. This rules out the possibility $n=8$, as $64<2^7$.
\smallskip

It remains to consider the case $n=6$.
\smallskip

By 2), $|\mathfrak{G}\cap\mathfrak{D}_1|=4$, $|\mathfrak{G}/(\mathfrak{G}\cap\mathfrak{D}_1)|=8$, $|\mathfrak{G}|=32$. Note that $f^2$ commutes with $t$. For, $f^2tf^{-2}=t^{q^2}=t$, as $q^2\equiv 1$ mod $8$.

\smallskip

Let
\[Y=\{e_1,\cdots, e_6, -e_1,\cdots, -e_6\}\] which is normalized by $\mathfrak{G}$.
\smallskip

--- \emph{Then $|I|\neq 4$ }:
\smallskip

\smallskip

Assume $|I|=4$. Then $f$ is of order $8$. Either $f$ commutes with $t$ or $ftf^{-1}=t^{-1}$. As $I$ is cyclic, one at least $I$-orbit in $Y$ has $4$ elements.
\smallskip

\smallskip

i) \emph{Case where exactly $1$ $I$-orbit in $Y$ has $4$ elements }:
\smallskip

This $I$-orbit, say $O$, is normalized by $f$, and $f$ acts simply transitively on $Y\backslash O=:O'$. On $O'$, as $t$ and $t^{-1}$ coincide, $t$ commutes with $f$ and thus acts as $f^4$ or $1$. If say
$O'=\{e_1,\cdots, e_4, -e_1,\cdots, -e_4\}$, then
$\{\pm e_1\pm\cdots\pm e_4\}$ is not acted transitively by $\langle f,t\rangle=\mathfrak{G}$ and thus
\[X=\{s_1e_1+\cdots+s_6e_6, s_i\in\{1,-1\}, s_1\cdots s_6=1\}\] is not acted transitively by $\mathfrak{G}$ either.
\smallskip

ii) \emph{Case where exactly $2$ $I$-orbits in $Y$ have $4$ elements }:
\smallskip

These two, say $O_1$, $O_2$, are exchanged by $f$. So $f^2$ normalizes and acts as $t$ or $t^{-1}$ on each, since $f^2$ commutes with $t$. Both $f^4$ and $t^2$ act as the identity on $Y\backslash (O_1\cup O_2)$. So $f^4=t^2$. But then $|\mathfrak{G}|$ divides $16$.
\smallskip

iii) \emph{Case where exactly $3$ $I$-orbits in $Y$ have $4$ elements }:
\smallskip

This contradicts the assumption that $t\in\mathfrak{W}_1$.
\smallskip

So $|I|=8$. Let $O$ be the unique $I$-orbit of cardinality $8$ in $Y$, say
\[O=\{e_1,\cdots, e_4, -e_1,\cdots, -e_4\}.\] Let $O':=Y\backslash O$. Then $f$ normalizes $O$ as well as $O'$. As $\mathfrak{G}$ acts transitively on $X$, it acts transitively on
$\{\pm e_1\pm\cdots\pm e_4\}$. That is, for each choice of $s_1,\cdots, s_4\in\{1,-1\}$, there are integers $i, j$ such that
\[s_1e_1+\cdots+s_4e_4=f^it^j(1+t+t^2+t^3)e_1.\]

In particular, there are $i, j\in\mathbf{Z}$ such that
\[(1-t+t^2+t^3)e_1=f^it^j(1+t+t^2+t^3)e_1.\]

Write $f(e_1)=t^{\mu}e_1$ for an integer $\mu$. Then $f^2=t^{(q+1)\mu}$ on $O$.
\smallskip

\smallskip

--- \emph{One has $q\not\equiv \pm 1$ mod $8$ }:
\smallskip

\smallskip

For, if $q\equiv 1$ mod $8$, then
\[f(1+t+t^2+t^3)e_1=(1+t+t^2+t^3)f(e_1)=t^{\mu}(1+t+t^2+t^3)e_1.\]
If $q\equiv -1$ mod $8$, then
\[f(1+t+t^2+t^3)e_1=(1+t^{-1}+t^{-2}+t^{-3})f(e_1)=t^{\mu-3}(1+t+t^2+t^3)e_1.\]

\smallskip

\smallskip

--- \emph{The group $I$ acts transitively on $O'$ }:
\smallskip

\smallskip

Assume that $O'$ consists of at least $2$ $I$-orbits. Choose $x'_1, x'_2\in O'$ such that
$(1+t+t^2+t^3)e_1+x'_1+x'_2\in X$. Now, $t$ is not $1$ or $-1$ on $O'$, since
\[t(1+t+t^2+t^3)e_1+t x'_1+t x'_2\in X.\]

One may assume $t x'_1=x'_1$, $t x'_2=-x'_2$. Then $f$ normalizes $\{x'_1, -x'_1\}$ as well as $\{x'_2, -x'_2\}$. So $f^2$ is the identity on $O'$. So $f^2=t^{(q+1)\mu}$. But then $|\mathfrak{G}|$ divides $16$.
\smallskip

As now $I$ acts transitively on $O'$, there exists $x'\in O'$ such that
$(1+t+t^2+t^3)e_1+(1+t)x'\in X$. One has
$f^2(x')\neq t^{(q+1)\mu}x'$. For otherwise $f^2=t^{(q+1)\mu}$ and $|\mathfrak{G}|$ divides $16$.
\smallskip

\smallskip

--- \emph{Then $q\not\equiv 3$ mod $8$ }:
\smallskip

\smallskip

Assume $q\equiv 3$ mod $8$. As $f^2(x')\neq t^{(q+1)\mu}x'=x'$, $f$ is of order $4$ on $O'$. So $f=t$ or $t^{-1}$ on $O'$. This contradicts the equation $ftf^{-1}=t^q=t^3$.
\smallskip

\smallskip

--- \emph{Then $q\not\equiv 5$ mod $8$ }:
\smallskip

\smallskip

Assume $q\equiv 5$ mod $8$. Then $f$ commutes with $t$ on $O'$ and so $f=t^{\nu}$ on $O'$ for an integer $\nu$. The condition $f^2(x')\neq t^{(q+1)\mu}x'$ says that $\nu-\mu$ is an odd integer. But $\nu-\mu$ should also be an even integer. For, the condition that $f$ normalizes $X$ implies that
\[t^{-\mu}f((1+t+t^2+t^3)e_1+(1+t)x')\in X,\] which is
\[(1-t+t^2-t^3)e_1+(1+t)t^{\nu-\mu}x'\in X.\]

\end{proof}
\smallskip

\smallskip

\smallskip


\section{Type ${}^2D$}
\smallskip

\smallskip

\smallskip

Let $(S,\eta,s)$, $\mathrm{char}(s)=\ell$, be as in $\S 4$.
\smallskip

Suppose given a $({}^2D_n, \alpha_1)$ over $\eta$, where $n$ is an integer $\geq 4$. Let
\[\rho_{{}^2D_n}: \pi_1(\eta, \overline{\eta})\to\{1,-1\}\] be the index of ${}^2D_n$. One says that ${}^2D_n$ is unramified (resp. ramified) over $S$ if its index is unramified (resp. ramified) over $S$ (\S 4).
\smallskip

Write $n=2^gr$, for an integer $g\geq 0$ and an odd integer $r\geq 1$.
\smallskip

Let $\mathbf{Z}^n$ be identified with
$\mathbf{Z}^{2^g}\otimes_{\mathbf{Z}}\mathbf{Z}^r$ in such a way that the standard basis
$e_1,\cdots,e_n$ of $\mathbf{Z}^n$ is identified with $e'_1\otimes e''_1,\cdots, e'_{2^g}\otimes
e''_r$, where $e'_1,\cdots, e'_{2^g}$ (resp. $e''_1,\cdots,e''_r$)
denote the standard basis of $\mathbf{Z}^{2^g}$ (resp. $\mathbf{Z}^r$).
\smallskip

We denote the group of diagonal (resp. monomial) matrices of $\mathrm{GL}_n(\mathbf{Z})$ by $\mathfrak{D}$ (resp. $\mathfrak{M}$). Let $\mathfrak{D}_1$ be the subgroup of $\mathfrak{D}$ consisting of all diagonal matrices of determinant $1$.
\smallskip

Let $\mathfrak{W}=\mathfrak{D}\mathfrak{M}$ and $\mathfrak{W}_1=\mathfrak{D}_1\mathfrak{M}$.

\smallskip

\smallskip

{\bf Proposition 10.1.} \emph{Suppose that ${}^2D_n$ is unramified over $S$. Then $({}^2D_n,\alpha_1)$ is elliptic.}

\begin{proof} Let $\zeta, \tau\in\mathfrak{W}$ be such that
\[\zeta: e_1\mapsto e_2,\cdots, e_n\mapsto e_1.\]
\[\tau: e_1\mapsto -e_1,\ e_i\mapsto e_i,\ \forall\ i>1.\]
The cyclic group $\mathfrak{G}$ generated by $\tau\zeta$ acts simply transitively on
\[\{e_1,\cdots,e_n,-e_1,\cdots,-e_n\}.\] Choose a surjective homomorphism (\S 4) :
\[\rho: \pi_1(\eta, \overline{\eta})\to\pi_1(S, \overline{\eta})\to\mathbf{Z}/2n\mathbf{Z}=\mathfrak{G}.\] The composition
\[\pi_1(\eta, \overline{\eta})\stackrel{\rho}{\longrightarrow}\mathfrak{G}\hookrightarrow\mathfrak{W}\to\mathfrak{W}/\mathfrak{W}_1=\{1,-1\}\] is the index of ${}^2D_n$, as ${}^2D_n$ is unramified over $S$. So $({}^2D_n,\alpha_1)$ is elliptic (3.1), 9).

\end{proof}
\smallskip

{\bf Proposition 10.2.} \emph{Suppose $\ell>2$ and that ${}^2D_n$ is ramified over $S$. Then $({}^2D_n, \alpha_1)$ is elliptic.}

\begin{proof} As $\ell>2$, the index of ${}^2D_n$ is tamely ramified over $S$ :
\[\rho_{{}^2D_n}:\pi_1(\eta, \overline{\eta})\to\pi_1^t(\eta, \overline{\eta})=\langle F, T\rangle\stackrel{\overline{\rho}}{\longrightarrow}\{1,-1\},\ \overline{\rho}:T\mapsto -1.\]

Let $q=\mathrm{Card}(k(s))$. Let $\tau', \sigma'\in\mathrm{GL}_{2^g}(\mathbf{Z})$ be such that
\[ \tau': e'_1\mapsto e'_2,\ \cdots,\ e'_{2^g-1}\mapsto e'_{2^g},\
e'_{2^g}\mapsto -e'_1,\]
\[\sigma'\tau'={\tau'}^{q^r}\sigma',\ \ \sigma': e'_1\mapsto e'_1.\]

Let $\tau\in\mathrm{GL}_n(\mathbf{Z})$ be such that
\[ \tau: e'_i\otimes e''_j\mapsto {\tau'}^{q^{j-1}}(e'_i)\otimes
e''_j,\ \forall\ j=1,\cdots,r,\ \forall\ i=1,\cdots,2^g.\]

And let $\sigma\in\mathrm{GL}_n(\mathbf{Z})$ be such that
\[ \sigma: e'_i\otimes e''_1\mapsto e'_i\otimes e''_2,\ \cdots,\
e'_i\otimes e''_{r-1}\mapsto e'_i\otimes e''_r,\ e'_i\otimes
e''_r\mapsto \sigma'(e'_i)\otimes e''_1\] $\forall\ i=1,\cdots,2^g$.
Then $\tau$ is of order $2^{g+1}$, $\sigma^r=\sigma'\otimes 1$ and
$\sigma\tau\sigma^{-1}=\tau^q$.
\smallskip

The subgroup $\mathfrak{G}$ of $\mathfrak{W}$ generated by $\{\sigma,\tau\}$ acts transitively on
\[\{e_1,\cdots,e_n,-e_1,\cdots,-e_n\}.\] Let
\[\rho: \pi_1^t(\eta, \overline{\eta})\to\mathfrak{G},\ T\mapsto \tau\] which maps $F$ to :
\smallskip

--- $\sigma$, if $\sigma\in\mathfrak{W}_1$, $\overline{\rho}: F\mapsto 1$.
\smallskip

--- $\sigma\tau$, if $\sigma\in\mathfrak{W}_1$, $\overline{\rho}: F\mapsto -1$.
\smallskip

--- $\sigma\tau$, if $\sigma\notin\mathfrak{W}_1$, $\overline{\rho}: F\mapsto 1$.
\smallskip

\smallskip

Then the composition
\[\pi_1(\eta, \overline{\eta})\stackrel{\rho}{\longrightarrow}\mathfrak{G}\hookrightarrow\mathfrak{W}\to\mathfrak{W}/\mathfrak{W}_1=\{1,-1\}\] is the index of ${}^2D_n$. So $({}^2D_n,\alpha_1)$ is elliptic (3.1), 9).

\end{proof}
\smallskip

Let $d$ (resp. $f$) be an integer $\geq 1$ (resp. $>1$). Let pro-$2$-groups $F_1, F_2, F_3, F_4$ be defined by generators and relations as :
\[F_1=\langle x_1,\cdots,x_{d+2}|\ x_1^{2^f}[x_1,x_2][x_3,x_4]\cdots
[x_{d+1},x_{d+2}]=1,\ d\ \mathrm{even}\rangle,
\]
\[F_2=\langle x_1,\cdots,x_{d+2}|\
x_1^2x_2^4[x_2,x_3]\cdots[x_{d+1},x_{d+2}]=1,\ d\
\mathrm{odd}\rangle,\]
\[F_3=\langle x_1,\cdots,x_{d+2}|\
x_1^{2+2^f}[x_1,x_2][x_3,x_4]\cdots[x_{d+1},x_{d+2}]=1,\ d\
\mathrm{even}\rangle,\]
\[F_4=\langle x_1,\cdots,x_{d+2}|\
x_1^2[x_1,x_2]x_3^{2^f}[x_3,x_4]\cdots[x_{d+1},x_{d+2}]=1,\ d\
\mathrm{even}\rangle,\] where
\[x,y\mapsto [x,y]=x^{-1}y^{-1}xy\] denotes the
commutator.
\smallskip

When $\ell=2$, $\pi_1(\eta,\overline{\eta})$
has one of the groups $F_1,F_2, F_3, F_4$ as the maximal pro-$2$-quotient, for $d=[\eta:\mathbf{Q}_2]$ and for a certain integer $f$ (\cite{serre}, p. 107--108).
\smallskip

\smallskip

{\bf Proposition 10.3.} \emph{Suppose $\ell=2$. Then
$({}^2D_n,\alpha_1)$ is elliptic.}

\begin{proof} Let $a',b'\in\mathrm{GL}_{2^g}(\mathbf{Z})$ be such that
\[ a': e'_1\mapsto -e'_1,\ e'_i\mapsto e'_i,\ \forall\ i>1,\]
\[ b': e'_1\mapsto e'_2,\ \cdots,\ e'_{2^g-1}\mapsto e'_{2^g},\
e'_{2^g}\mapsto e'_1.\]

Let $c''\in\mathrm{GL}_r(\mathbf{Z})$ be such that
\[ c'': e''_1\mapsto e''_2,\ \cdots,\ e''_{r-1}\mapsto e''_r,\
e''_r\mapsto e''_1.\]

Let $a=a'\otimes 1, b=b'\otimes 1, c=1\otimes
c''\in\mathfrak{W}$.
\smallskip

Notice that $a$
(resp. $b$, resp. $c$) has image $-1$ (resp. $1$, resp. $1$) by
\[\mathfrak{W}\to\mathfrak{W}/\mathfrak{W}_1=\{1,-1\}.\]

The group $\langle ab\rangle\times \langle c\rangle$ acts simply
transitive on \[\{e_1,\cdots,e_n,-e_1,\cdots,-e_n\}.\]

By (3.1), 9) it suffices to show that either $\langle
ab\rangle\times \langle c\rangle$ or $\langle a,b\rangle\times
\langle c\rangle$ is realizable as a quotient of $\pi_1(\eta,\overline{\eta})$
lifting the index of ${}^2D_n$. Now, this index factors through the maximal pro-$2$-quotient $F$ of $\pi_1(\eta, \overline{\eta})$ :
\[\pi_1(\eta,\overline{\eta})\to F\stackrel{\chi}{\longrightarrow}\{1,-1\}.\] And, the odd order cyclic subgroup $\langle c\rangle$ of $\mathfrak{W}_1$ is realizable as an unramified quotient
of $\pi_1(\eta,\overline{\eta})$. So it
suffices to show that every surjective homomorphism
\[\chi: F\to\{1,-1\}\] is a composition of the form
\[F\stackrel{\rho}{\longrightarrow}\langle a,b\rangle\hookrightarrow\mathfrak{W}\to\mathfrak{W}/\mathfrak{W}_1=\{1,-1\}\] for some representation
\[\rho: F\to \langle a, b\rangle\] whose image is $\langle
ab\rangle$ or $\langle a,b\rangle$.
\smallskip

Given the explicit structure of $F$ as above, the verification is straightforward. Consider for example the case where
\[F=\langle x,y,z| x^2y^4[y,z]=1\rangle\] and where $g\geq 2$. According to
the values of $\chi$ on $(x,y,z)$, one defines $\rho: F\to \langle
a,b\rangle$ as follows\ :
\smallskip

1) $(-1,1,1)$. Let $\rho: (x,y,z)\mapsto (a,1,b)$.
\smallskip

2) $(1,-1,1)$. Let $\rho: (x,y,z)\mapsto ((ab)^{-2}, ab, 1)$.
\smallskip

3) $(1,1,-1)$. Let $\rho: (x,y,z)\mapsto (1,1,ab)$.
\smallskip

4) $(-1,1,-1)$. Let $\rho: (x,y,z)\mapsto (a,1,ab)$.
\smallskip

5) $(1,-1,-1)$. Let $\rho: (x,y,z)\mapsto ((ab)^{-2}, ab, ab)$.
\smallskip

6) $(-1,-1,1)$. Let $\rho: (x,y,z)\mapsto (a,ab,ab^2ab^{-2})$, if
$g=2$, and let $\rho: (x,y,z)\mapsto (ab^2,ab^{-1}, ab^3ab^{-3})$, if
$g>2$.
\smallskip

7) $(-1,-1,-1)$. Let $\rho: (x,y,z)\mapsto (ab^2, ab, ab^{-1})$, if
$g=2$, and let $\rho: (x,y,z)\mapsto (b^{-1}ab^2aba, ab^{-1}, ab)$, if
$g>2$.

\end{proof}
\smallskip

\smallskip

\smallskip


\section{Type $E_6$}
\smallskip

\smallskip

\smallskip

Let $E$ be a $6$-dimensional $\mathbf{F}_2$-vector space. Let
$e_i,f_j$, $1\leq i,j\leq 3$, be a basis of $E$ and let $q$ be the quadratic form on $E$ such that
\[ q(e_i)=q(f_j)=1,\ q(e_i+e_j)=q(f_i+f_j)=0,\
q(e_i+f_j)=\delta_{ij},\] where $\delta_{ij}=1$, if $i=j$, and $\delta_{ij}=0$,
if $i\neq j$, $\forall\ i,j\in\{1,2,3\}$.
\smallskip

Let
\[X= \{v\in E\backslash\{0\},\ q(v)=0\}\] be the $q$-singular vectors of $E\backslash\{0\}$.
\smallskip

Let $V_i=\mathbf{F}_2e_i+\mathbf{F}_2f_i$,
$i=1,2,3$. The elements of $X$ are of the form
$v_i+v_j$, where $v_i\in V_i$, $v_j\in V_j$, $1\leq i, j\leq 3$, $i\neq j$, $v_i,v_j\neq 0$. The set $X$ consists of $27$ vectors which are permuted transitively by the orthogonal group $\mathrm{O}(q)$. The group $\mathrm{O}(q)$ has $2^7.3^4.5$ elements.
\smallskip

Observe that an element of $\mathrm{GL}(E)$ belongs to $\mathrm{O}(q)$ if and only if it normalizes $X$.
\smallskip


\smallskip

Note that, for each $i\in\{1,2,3\}$, one has
$\mathrm{O}(q|V_i)=\mathrm{GL}(V_i)$, because $q(e_i)=q(f_i)=q(e_i+f_i)=1$.
The subgroup $\mathrm{GL}(V_1)\times\mathrm{GL}(V_2)\times\mathrm{GL}(V_3)$ of $\mathrm{O}(q)$ consists of all elements $g$ such that $g(V_1)=V_1$, $g(V_2)=V_2$, $g(V_3)=V_3$.
\smallskip

Let $\mathfrak{N}$ be the subgroup of $\mathrm{O}(q)$ consisting of all elements $g$ such that $g(V_i)\in\{V_1, V_2, V_3\}$, $\forall\ i=1,2,3$. One has a split exact sequence
\[1\to\prod^3_1\mathrm{GL}(V_i)\to \mathfrak{N}\to\mathrm{Aut}(\{V_1, V_2, V_3\})\to 1.\]

The unique $3$-Sylow subgroup $\mathfrak{M}$ of $\prod_i\mathrm{GL}(V_i)$ is the unique abelian subgroup of order $27$ of $\mathfrak{N}$.
\smallskip

\smallskip

{\bf Lemma 11.1.} \emph{Suppose that a solvable subgroup $\mathfrak{G}$ of $\mathrm{O}(q)$ acts transitively on $X$. Then $5\nmid |\mathfrak{G}|$.}

\begin{proof} Let $\mathfrak{H}$ be a Hall subgroup of $\mathfrak{G}$ which is a product of a $3$-Sylow subgroup and a $5$-Sylow subgroup of $\mathfrak{G}$. Then $\mathfrak{H}$ also acts transitively on $X$.
\smallskip

\smallskip

--- \emph{Case $3^4$ does not divide $|\mathfrak{H}|$ }:
\smallskip

\smallskip

Then $\mathfrak{H}$ has a unique $5$-Sylow subgroup, say $\mathfrak{Q}$. The $\mathfrak{Q}$-orbits in $X$ all have the same cardinality, say $r$, which divides both $5$ and $\mathrm{Card}(X)=27$. So $r=1$. So $\mathfrak{Q}=1$. For, if an element $g\in\mathrm{GL}(E)$ restricts to the identity on $X$, then $g=1$.
\smallskip

\smallskip

--- \emph{Case $3^4$ divides $|\mathfrak{H}|$ }:
\smallskip

\smallskip

Then $\mathfrak{H}$ has a unique $3$-Sylow subgroup $\mathfrak{P}$ which one may, by conjugating $\mathfrak{H}$ in $\mathrm{O}(q)$, assume to be in $\mathfrak{N}$. In particular, $\mathfrak{P}$ contains $\mathfrak{M}$. Let $\mathfrak{Q}$ be a $5$-Sylow subgroup of $\mathfrak{H}$. Then $\mathfrak{Q}$ normalizes and thus centralizes $\mathfrak{M}$. Notice that $\mathfrak{M}$ has $3$-orbits in $X$ each of which consists of $3$ elements. The group $\mathfrak{Q}$ normalizes each of these $3$-orbits and so it fixes every point of $X$. So $\mathfrak{Q}=1$.

\end{proof}

\smallskip

{\bf Lemma 11.2.} \emph{Suppose that a solvable subgroup $\mathfrak{G}$ of $\mathrm{O}(q)$ acts transitively on $X$. Let $\mathfrak{A}$ be an abelian normal subgroup of $\mathfrak{G}$. Then $\mathfrak{A}$ is a $3$-group.}

\begin{proof} Let $\mathfrak{Q}$ be the unique $2$-Sylow subgroup of $\mathfrak{A}$. The group $\mathfrak{Q}$ is normal in $\mathfrak{G}$. So the $\mathfrak{Q}$-orbits in $X$ all have the same cardinality, say $r$, which divides both $|\mathfrak{Q}|$ and $27=\mathrm{Card}(X)$. So $r=1$. So $\mathfrak{Q}=1$. By (11.1), the lemma follows.

\end{proof}





\smallskip

{\bf Proposition 11.3.} \emph{Let $\mathfrak{H}$ be a cyclic subgroup of order $9$ of $\mathrm{O}(q)$. Let $h$ be a generator of $\mathfrak{H}$.}
\smallskip

\emph{Then the commutant of $\mathfrak{H}$ on $E$ is a field of cardinality $64$. And $\mathfrak{H}$ has $3$ orbits on $X$ each of which consists of $9$ points. Let $x, y\in X$ be in distinct $\mathfrak{H}$-orbits. Then there exists a unique element $g\in\mathrm{O}(q)$ of order $6$ which satisfies the following properties :
\[ghg^{-1}=h^2,\ g(x)=y.\] The group $\mathfrak{G}$ generated by $\{h,g\}$ acts transitively on $X$ and $\mathfrak{G}$ is of order $54$.}

\begin{proof} Notice that $[\mathbf{Q}(\mu_9):\mathbf{Q}]=6$ and that $2$ is inert in $\mathbf{Q}(\mu_9)$. So $\mathfrak{H}$ acts irreducibly on $E$ and the commutant $C$ of $\mathfrak{H}$ on $E$ is a field of cardinality $64$. In particular, neither $h$ nor $h^3$ fixes a nonzero vector in $E$, as $E$ is a $1$-dimensional $C$-vector space. Every $\mathfrak{H}$-orbit in $X$ consists of $9$ points.
\smallskip

Let $F: C\to C$, $c\mapsto c^2$, be the Frobenius automorphism of $C$. An element $g\in\mathrm{GL}(E)$ satisfying $ghg^{-1}=h^2$ is simply an $F$-linear automorphism of the $1$-dimensional $C$-vector space $E$. Let $g$ be an $F$-linear automorphism of $E$ and let $x\in E$ be a nonzero vector. For every integer $n$, $g^n$ is $F^n$-linear. So $g$ is of order a multiple of $6$. Write $g(x)=c.x$ for an element $c\in C^{\times}$. As $g^6$ is $C$-linear and as
\[g^6(x)=F^5(c)\cdots F(c)c.x=c^{63}.x=x,\] $g$ is of order $6$.
\smallskip

Note finally that $g$ lies in $\mathrm{O}(q)$ if and only if it normalizes $X$. From here, the claimed existence and uniqueness of $g$ as well as the last assertion immediately follow.

\end{proof}

\smallskip

{\bf Proposition 11.4.} \emph{Suppose that a solvable subgroup $\mathfrak{G}$ of $\mathrm{O}(q)$ acts transitively on $X$. Suppose furthermore that $\mathfrak{G}$ has a cyclic normal subgroup $\mathfrak{H}$ of generator $h$ of order $9$. Then $|\mathfrak{G}|=27$ or $54$.}
\smallskip

--- \emph{Case $|\mathfrak{G}|=27$. Then $\mathfrak{G}$ is generated by $\{h, g\}$ where the element $g$ is of order $3$ and satisfies $ghg^{-1}=h^4$.}
\smallskip

--- \emph{Case $|\mathfrak{G}|=54$. Then $\mathfrak{G}$ is generated by $\{h, g\}$ where the element $g$ is of order $6$ and satisfies $ghg^{-1}=h^2$.}

\begin{proof} Let $C$ be the commutant of $\mathfrak{H}$ on $E$. By (11.3), $C$ is a field of cardinality $64$ and $E$ is a $1$-dimensional $C$-vector space. The centralizer of $\mathfrak{H}$ in $\mathfrak{G}$ is the intersection $\mathfrak{G}\cap C^{\times}$, that is, $\mathfrak{H}$. Now, the exact sequence
\[1\to\mathfrak{H}\to\mathfrak{G}\stackrel{g\mapsto\mathrm{Int}(g)|\mathfrak{H}}{\longrightarrow}\mathrm{Aut}(\mathfrak{H})\] shows that $\mathfrak{G}$ is of order $27$ or $54$, as $\mathrm{Aut}(\mathfrak{H})$ is cyclic of order $6$.
Choose $g\in\mathfrak{G}$ such that $\mathrm{Int}(g)|\mathfrak{H}$ generates $\mathrm{Int}(\mathfrak{G})|\mathfrak{H}$. The automorphism $\mathrm{Int}(g)|\mathfrak{H}$ extends to an automorphism of the field $C$.
\smallskip

--- \emph{Case where $|\mathfrak{G}|=27$ }:
\smallskip

Replacing if necessary $g$ by its inverse, one may assume $ghg^{-1}=h^4$. As $\mathrm{GL}(E)$ has no element of order $27$, one has $g^3=h^{3n}$ for some integer $n$. The group $\mathfrak{G}$ is generated by $\{h, gh^{-n}\}$ and
\[(gh^{-n})^3=(gh^{-n}g^{-1})(g^2h^{-n}g^{-2})(g^3h^{-n}g^{-3})g^3=h^{-4n}h^{-16n}h^{-64n}g^3=1.\]

--- \emph{Case where $|\mathfrak{G}|=54$ }:
\smallskip

Replacing if necessary $g$ by its inverse, one may assume $ghg^{-1}=h^2$. As in (11.3), $g$ is of order $6$. And $\mathfrak{G}$ is generated by $\{h, g\}$.

\end{proof}

\smallskip

\smallskip

{\bf Proposition 11.5.} \emph{Let $(S, \eta, s)$, $\mathrm{char}(k(s))=\ell$, be as in $\S 4$. Suppose that $\ell=3$. Then $(E_6, \alpha_1)$ and $(E_6, \alpha_6)$ are elliptic over $\eta$.}

\begin{proof} The orthogonal group $\mathrm{O}(q)$ has a subgroup $\mathfrak{G}$ of order $27$ which acts transitively on $X$ and which is generated by $2$ elements $h,g$, where $h$ (resp. $g$) is of order $9$ (resp. $3$) and $ghg^{-1}=h^4$.
\smallskip

--- \emph{Case where $\mu_3(k(\eta))=1$ }:
\smallskip

In this case, the maximal pro-$3$-quotient of $\pi_1(\eta, \overline{\eta})$ is free of rank $\geq 2$ as a pro-$3$-group. In particular, $\mathfrak{G}$ is realizable as a quotient of $\pi_1(\eta, \overline{\eta})$. So $(E_6, \alpha_1)$ and $(E_6, \alpha_6)$ are elliptic over $\eta$ (3.1), 6).
\smallskip

--- \emph{Case where $\mu_3(k(\eta))=\mu_3(k(\overline{\eta}))$ }:
\smallskip

The maximal pro-$3$-quotient of $\pi_1(\eta, \overline{\eta})$ has then the presentation :
\[F=\langle x_1,\cdots,x_{d+2}\ |\
x_1^q[x_1,x_2][x_3,x_4]\cdots[x_{d+1},x_{d+2}]=1\rangle\] where
$d=[k(\eta):\mathbf{Q}_3]$, where $q$ is the maximal power of $3$ such that
$\mu_q(k(\eta))=\mu_q(k(\overline{\eta}))$ and where $(x,y)\mapsto[x,y]=x^{-1}y^{-1}xy$ is
the commutator. The homomorphism $\chi: F\to\mathfrak{G}$ such that
\[\chi: x_1\mapsto 1,\ x_2\mapsto h,\ x_3\mapsto g,\ x_i\mapsto 1, \forall\ i>3\] is surjective. So again $(E_6,\alpha_1)$ and $(E_6,\alpha_6)$ are elliptic over $\eta$ (3.1), 6).

\end{proof}

\smallskip

{\bf Proposition 11.6.} \emph{Let $(S,\eta,s)$, $\mathrm{char}(k(s))=\ell$, be as in $\S 4$. Suppose $\ell\neq 3$. Then $(E_6, \alpha_1)$ and $(E_6, \alpha_6)$ are elliptic over $\eta$ if and only if $\mathrm{Card}(k(s))\equiv \pm 2, \pm 4$ mod $9$.}

\begin{proof} By (3.1), 6), the pairs $(E_6,\alpha_1)$ and $(E_6,\alpha_6)$ are elliptic over $\eta$ if and only if there is a representation
\[\rho: \pi_1(\eta, \overline{\eta})\to\mathrm{O}(q)\] whose image acts transitively on $X$. Suppose that such a representation exists. Let $\mathfrak{G}$ be its image. Let $I$ (resp. $P$) be the image in $\mathfrak{G}$ of the inertia (resp. wild inertia) subgroup of $\pi_1(\eta,\overline{\eta})$.
\smallskip

As $P$ is normal in $\mathfrak{G}$, the $P$-orbits in $X$ all have the same cardinality, say $r$, which divides both $27$ and $|P|$. That is, $r=1$ and $P=1$.
\smallskip

So $I=I/P$ is cyclic of order a power of $3$ (11.2) and so $\mathfrak{G}$ has a unique $3$-Sylow subgroup, say $\mathfrak{H}$. As $\mathfrak{H}/I$ is cyclic, the group $I$ is cyclic of order $9$ and $\mathfrak{H}/I$ is of order $3$. The quotient $\mathfrak{G}/\mathfrak{H}$ is a cyclic $2$-group (11.1). So $|\mathfrak{G}|=27$ or $54$ (11.4). Write $\rho$ as a composition
\[\pi_1(\eta,\overline{\eta})\to\pi_1^t(\eta,\overline{\eta})\stackrel{\chi}{\longrightarrow}\mathfrak{G},\] where $\pi_1^t(\eta,\overline{\eta})=\langle F,T\rangle$ (\S 4). The image $t$ of $T$ in $\mathfrak{G}$ generates $I$. Let $f$ be the image of $F$ in $\mathfrak{G}$.
\smallskip

--- \emph{Case where $|\mathfrak{G}|=54$ }:
\smallskip

Then $f$ is of order $6$. One has $ftf^{-1}=t^2$ or $t^{32}=t^{-4}$.
\smallskip

--- \emph{Case where $|\mathfrak{G}|=27$ }:
\smallskip

Then $f$ is of order $3$ or $9$. One has $ftf^{-1}=t^4$ or $t^{16}=t^{-2}$.
\smallskip

Such groups do exist in $\mathrm{O}(q)$ (11.3)

\end{proof}

\smallskip

\smallskip

\smallskip


\section{Type $E_7$}
\smallskip

\smallskip

\smallskip

Let $E$ be a $6$-dimensional $\mathbf{F}_2$-vector
space equipped with a symplectic form $(,)$. Let $e_i,f_j$, $1\leq i,j\leq 3$, be a sympletic base of $E$. Let $q$ be the quadratic form on $E$ satisfying
\[q(e_i)=q(f_j)=1,\ q(e_i+e_j)=q(f_i+f_j)=0,\
q(e_i+f_j)=\delta_{ij}\] where $\delta_{ij}=1$, if $i=j$, and $\delta_{ij}=0$,
if $i\neq j$, $\forall\ i,j\in\{1,2,3\}$.
\smallskip

Observe that the orthogonal group $\mathrm{O}(q)$ is a subgroup of the symplectic group
$\mathrm{Sp}(E)$.
\smallskip

The group $\mathrm{Sp}(E)$ is of order $2^9.3^4.5.7$, the subgroup $\mathrm{O}(q)$ is of order $2^7.3^4.5$ and the homogenous space
\[X=\mathrm{Sp}(E)/\mathrm{O}(q)\] consists of $28$ elements.
\smallskip

We shall determine up to conjugation all solvable subgroups $\mathfrak{G}$ of
$\mathrm{Sp}(E)$ that act transitively on $X$.
\smallskip

Each such group $\mathfrak{G}$ contains a $7$-Sylow subgroup of $\mathrm{Sp}(E)$.
By conjugation in $\mathrm{Sp}(E)$, one may suppose that $\mathfrak{G}$ contains
$\zeta\in\mathrm{Sp}(E)$, where
\[ \zeta: \left\{ \begin{array}{ll}
          e_1\mapsto e_2,\ e_2\mapsto e_3,\ e_3\mapsto e_1+e_2 & \\
          f_1\mapsto f_1+f_2,\ f_2\mapsto f_3,\ f_3\mapsto f_1
          \end{array} \right. \]

Let
$V=\mathbf{F}_2e_1+\mathbf{F}_2e_2+\mathbf{F}_2e_3$,
$V^{\vee}=\mathbf{F}_2f_1+\mathbf{F}_2f_2+\mathbf{F}_2f_3$. Then
\[ \mathrm{det}(T-\zeta,V)=T^3+T+1\ ,\
\mathrm{det}(T-\zeta,V^{\vee})=T^3+T^2+1\] and
\[ \mathrm{det}(T-\zeta,E)=(T^3+T+1)(T^3+T^2+1)=(T^7-1)/(T-1).\]

As $\zeta$-modules, $V,V^{\vee}$ are irreducible mutually
non-isomorphic. The subspaces $0,V,V^{\vee},E$ are the only
sub-$\zeta$-modules of $E$.
\smallskip

The commutant $\mathrm{End}_{\zeta}(E)$ is equal to
$\mathbf{F}_2[\zeta|V]\times\mathbf{F}_2[\zeta|V^{\vee}]$. And
\[\mathrm{GL}_{\zeta}(E)\cap\mathrm{Sp}(E)=\mathbf{F}_2[\zeta]^{\times}=\langle\zeta\rangle.\]
That is, $\langle\zeta\rangle$ is its own centralizer in
$\mathrm{Sp}(E)$.
\smallskip

The normalizer of $\langle\zeta\rangle$ in $\mathrm{Sp}(E)$ admits $2$ generators
$\zeta,\sigma$, where
\[\sigma: \left\{ \begin{array}{cc}
           e_1\mapsto f_1,\ e_2\mapsto f_2,\ e_3\mapsto f_2+f_3 & \\
           f_1\mapsto e_1,\ f_2\mapsto e_2+e_3,\ f_3\mapsto e_3
           \end{array} \right. \]
And $\sigma, \zeta$ satisfy the relations :
\[\sigma^6=1,\ \sigma\zeta\sigma^{-1}=\zeta^{-2}.\] Notice that $|\langle\zeta,\sigma\rangle|=42$.
\smallskip

Let $\mathfrak{S}$ be the subgroup of $\mathrm{Sp}(E)$ consisting of all elements which act as the identity on $V$. By $g\mapsto (g-1)|V^{\vee}$, $\mathfrak{S}$ can be identified with
an $\mathbf{F}_2$-vector space of dimension $6$ which consists of all linear transformations
$A:V^{\vee}\to V$ such that the bilinear form
\[ u',v'\mapsto (u',Av')\] is symmetric in $u',v'\in V^{\vee}$.
\smallskip

For all $g\in\mathfrak{S}$, the function $v'\mapsto (v',(g-1)v')$ is
linear on $V^{\vee}$. Thus there is a unique vector
$v_g\in V$ satisfying
\[(v',(g-1)v')=(v_g, v'),\ \ \forall\ v'\in
V^{\vee}.\]

The function $\mathfrak{S}\to V$, $g\mapsto v_g$, is linear whose
kernel $\mathfrak{S}^1$ consists of all those
$g\in\mathfrak{S}$ such that the form
\[u',v'\mapsto (u',(g-1)v')\] is
alternating, i.e., that
\[(u',(g-1)v')=(u'\wedge v', \omega_g)\] for a
uniquely determined $2$-form $\omega_g\in\wedge^2V$.
\smallskip

The map
$g\mapsto\omega_g$ establishes a canonical bijection between
$\mathfrak{S}^1$ and $\wedge^2V$. The exact sequence
\[ 0\to\mathfrak{S}^1\to \mathfrak{S}\stackrel{g\mapsto v_g}{\longrightarrow} V\to 0\] is uniquely split as $\zeta$-modules. For, $\wedge^2V=\mathfrak{S}^1$ and $V$ are non-isomorphic
$\zeta$-modules. Let $\mathfrak{S}^2$ denote this complement of $\mathfrak{S}^1$ in $\mathfrak{S}$. So
$\mathfrak{S}=\mathfrak{S}^1\oplus\mathfrak{S}^2$.
\smallskip

In terms of matrices, every element $g\in\mathfrak{S}$ is of the form
\[ g: \left\{ \begin{array}{ll}
      e_i\mapsto e_i\ , \ i=1,2,3 & \\
      f_i\mapsto f_i+\sum_{j=1,2,3}A_{ji}e_j \end{array} \right.\]
where $A_{ij}$ is a symmetric matrix with coefficients in
$\mathbf{F}_2$.
\smallskip

The element $g$ belongs to $\mathfrak{S}^1$ if and only if
$A_{11}=A_{22}=A_{33}=0$. The $\zeta$-module $\mathfrak{S}^1$ is generated by $g_1$, where
\[g_1: \left\{ \begin{array}{ll}
        e_i\mapsto e_i,\ i=1,2,3 & \\
        f_1\mapsto f_1+e_2+e_3,\ f_2\mapsto f_2+e_1+e_3,\ f_3\mapsto
        f_3+e_1+e_2 \end{array} \right.\]

The $\zeta$-module $\mathfrak{S}^2$ is generated by $g_2$, where
\[g_2\ : \left\{ \begin{array}{ll}
          e_i\mapsto e_i,\ i=1,2,3 & \\
          f_1\mapsto f_1+e_2+e_3,\ f_2\mapsto f_2+e_1+e_3,\ f_3\mapsto f_3+e_1+e_2+e_3
          \end{array} \right. \]

The element $g\in\mathfrak{S}$ preserves the quadratic form $q$ if and
only if $A_{12}=A_{23}=A_{13}$. One has
$\mathfrak{S}^1\cap\mathrm{O}(q)=\{1, g_1\}$ and $\mathfrak{S}^2\cap\mathrm{O}(q)=\{1,g_2\}$.

\smallskip

\smallskip

{\bf Proposition 12.1.} \emph{Up to conjugation all solvable subgroups of $\mathrm{Sp}(E)$ that act transitively on $X$ are enumerated as follows }:
\smallskip

--- \emph{$\langle\zeta\rangle\mathfrak{S}$.}
\smallskip

--- \emph{$\langle\zeta,\sigma^2\rangle\mathfrak{S}$.}
\smallskip

--- \emph{$\langle\zeta\rangle\mathfrak{S}^1$, $\langle\zeta\rangle\mathfrak{S}^2$.}
\smallskip

--- \emph{$\langle\zeta,\sigma^2\rangle\mathfrak{S}^1$, $\langle\zeta, \sigma^2\rangle\mathfrak{S}^2$.}

\begin{proof} Suppose that $\mathfrak{G}$ is a solvable subgroup of $\mathrm{Sp}(E)$ which acts
transitively on $X$. Up to conjugation in $\mathrm{Sp}(E)$ one may assume that $\zeta\in\mathfrak{G}$. Recall that $|\mathrm{Sp}(E)|=2^9.3^4.5.7$, $|\mathrm{O}(q)|=2^7.3^4.5$.
\smallskip

\smallskip

--- \emph{Then $5\nmid |\mathfrak{G}|$ } :

\smallskip

\smallskip

Otherwise, as it is solvable, $\mathfrak{G}$ has a Hall subgroup of order $35$, say $\mathfrak{Q}$, which is cyclic. But $\mathbf{Z}/35\mathbf{Z}$ admits no faithful
$6$-dimensional representations over $\mathbf{F}_2$.
\smallskip

\smallskip

Thus $|\mathfrak{G}|=2^a.3^b.7$, for an integer $a\geq 2$ and an integer $0\leq b\leq 4$.
\smallskip

Let $\mathfrak{L}$ be a Hall subgroup of $\mathfrak{G}$ which is a product of $\langle\zeta\rangle$ and a $3$-Sylow subgroup of $\mathfrak{G}$. As $b\leq 4$, $\langle\zeta\rangle$ is normal in $\mathfrak{L}$.
So $\mathfrak{L}$ is a subgroup of $\langle\zeta,\sigma\rangle$. So
$\mathfrak{L}=\langle\zeta\rangle$ or
$\langle\zeta,\sigma^2\rangle$. In particular, $b=0$ or $1$.
\smallskip

Let $\mathfrak{H}$ be a Hall subgroup of $\mathfrak{G}$ which is a product of $\langle\zeta\rangle$ and a $2$-Sylow subgroup of $\mathfrak{G}$. As $a\geq 2$, $\mathfrak{H}$ is not a subgroup of
$\langle\zeta,\sigma\rangle$. That is to say, $\langle\zeta\rangle$ is not
normal in $\mathfrak{H}$. Let $\mathfrak{A}$ be a maximal abelian normal subgroup of the solvable group $\mathfrak{H}$.
\smallskip

\smallskip

--- \emph{The group $\mathfrak{A}$ is a $2$-group } :
\smallskip

\smallskip

For otherwise the unique $7$-Sylow subgroup of $\mathfrak{A}$ would be normal in
$\mathfrak{H}$.
\smallskip

\smallskip

--- \emph{The group $\mathfrak{A}$ is the unique $2$-Sylow subgroup of $\mathfrak{H}$ }:
\smallskip

\smallskip

As $\mathfrak{A}$ is a $2$-group, the subspace $E^{\mathfrak{A}}$ of $E$ consisting of all vectors fixed by
$\mathfrak{A}$ is a non-zero $\mathfrak{H}$-module. So $E^{\mathfrak{A}}$ is either $V$ or $V^{\vee}$.
Replacing $\mathfrak{G}$ by $\sigma\mathfrak{G}\sigma^{-1}$ if
necessary, we suppose $E^{\mathfrak{A}}=V$. Thus $\mathfrak{A}$ is a subgroup of $\mathfrak{S}$. Notice that $\sigma^3$ does not
normalize $V$. So
$\mathfrak{H}\cap\langle\zeta,\sigma\rangle=\langle\zeta\rangle$ and so $\mathfrak{H}$ has
$2^a=|\mathfrak{H}/\langle\zeta\rangle|$
$7$-Sylow subgroups. Then $\mathfrak{H}$ has a unique $2$-Sylow subgroup, say $\mathfrak{a}$, because
$2^a.7-2^a(7-1)=2^a$. Then
$E^{\mathfrak{a}}$ is a non-zero sub-$\mathfrak{H}$-module of
$E^{\mathfrak{A}}=V$. So $E^{\mathfrak{a}}=V$. So
$\mathfrak{a}\leq\mathfrak{S}$. Thus $\mathfrak{a}$ is abelian.
One concludes that $\mathfrak{A}=\mathfrak{a}$.
\smallskip

In particular, $\mathfrak{H}\leq\langle\zeta\rangle\mathfrak{S}$
and $\mathfrak{G}=\mathfrak{L}\mathfrak{H}\leq
\langle\zeta,\sigma^2\rangle\mathfrak{S}$.
\smallskip

To finish, it suffices to show that both
$\langle\zeta\rangle\mathfrak{S}^1$ and $\langle\zeta\rangle\mathfrak{S}^2$ act transitively on
$X$. Both have $56$ elements. And it is immediate to verify that each intersects $\mathrm{O}(q)$ in two elements.

\end{proof}
\smallskip

{\bf Proposition 12.2.} \emph{Let $(S,\eta, s)$, $\mathrm{char}(s)=\ell$, be as in $\S 4$.
Then $(E_7,\alpha_7)$ is elliptic over $\eta$ if and only if $\ell=2$.}

\begin{proof} By (12.1) all solvable subgroup of
$\{1,-1\}\times\mathrm{Sp}(E)$ that act transitively on
\[\{1,-1\}\times (\mathrm{Sp}(E)/\mathrm{O}(q))=\{1,-1\}\times X\] contain elementary $2$-groups of $2$-rank $\geq 3$. So $(E_7, \alpha_7)$ is not elliptic if $\ell>2$ (3.1), 7).
\smallskip

Suppose $\ell=2$. Then $\mathfrak{G}:=\{1,-1\}\times
\langle\zeta\rangle\mathfrak{S}$ is a quotient of
$\pi_1(\eta,\overline{\eta})$ by (4.1) and because
$\langle\zeta\rangle\mathfrak{S}$ has no index $2$ subgroups. Moreover, $\mathfrak{G}$ acts transitively on $\{1,-1\}\times X$ (12.1). So $(E_7,\alpha_7)$ is elliptic when
$\ell=2$ (3.1), 7).

\end{proof}

\bibliographystyle{amsplain}

\end{document}